\documentclass[preprint,authoryear,11pt]{elsarticle}
\usepackage{amssymb,amsmath,latexsym}
\usepackage[T1]{fontenc}
\usepackage{amsmath,amsfonts,amssymb}
\usepackage[english]{babel}
\usepackage{latexsym}
\usepackage{epsfig}
\usepackage{color}
\usepackage{graphics}
\usepackage{enumerate}
\usepackage{multirow}
\usepackage{tabularx}
\usepackage{natbib}
\usepackage[left=2.5cm,right=2cm,top=2cm,bottom=2cm]{geometry}

\makeatletter

\@addtoreset{equation}{section}
\makeatother

\newtheorem{thm}{Theorem}
\newtheorem{remi}{Remark}
\newtheorem{prop}{Proposition}
\newtheorem{cor}{Corollary}

\def\zak{\null\hfill{$\Box$}\par\vspace*{0.2cm}}
\def\Vec{\mathrm{vec}}

\usepackage{indic}
\usepackage{here}
\usepackage{float}
\RequirePackage[colorlinks,citecolor=blue,urlcolor=blue]{hyperref}
\hypersetup{pdfpagemode=FullScreen}
\hypersetup{
    pdftoolbar=false,        
    pdfmenubar=false,        
    pdffitwindow=false,     
    pdfstartview={FitH},    
    pdftitle={title},    
    pdfnewwindow=true,      
    colorlinks=true,       
    linkcolor=blue,          
    citecolor=blue,        
    filecolor=magenta,      
    urlcolor=cyan,           
    pdftoolbar=false,
   pdfmenubar=false
}
5

\begin{document}


\begin{frontmatter}

\title{Estimation of multivariate asymmetric power GARCH models}

\author[Yacouba]{Y. Boubacar Ma{\"\i}nassara}
\address[Yacouba]{Universit\'e Bourgogne Franche-Comt\'e, \\
Laboratoire de math\'{e}matiques de Besan\c{c}on, \\ UMR CNRS 6623, \\
16 route de Gray, \\ 25030 Besan\c{c}on, France.}
\ead{mailto:yacouba.boubacar\_mainassara@univ-fcomte.fr}

\author[Yacouba]{O. Kadmiri}
\ead{mailto:othman.kadmiri@univ-fcomte.fr}

\author[Yacouba]{B. Saussereau}
\ead{mailto:bruno.saussereau@univ-fcomte.fr}

\begin{abstract}
It is now widely accepted that volatility models have to incorporate the so-called leverage effect in order to  model the dynamics of daily financial returns. 
We suggest a new class of multivariate power transformed asymmetric models. It includes several functional forms of multivariate GARCH models which are of great interest in  financial modeling and time series literature. We provide an explicit necessary and sufficient condition to establish the strict stationarity of the model. We derive the asymptotic properties of the quasi-maximum likelihood estimator of the parameters. These properties are established both when the power of the transformation is known or is unknown. The asymptotic results are illustrated by Monte Carlo experiments. An application to real financial data is also proposed.
\end{abstract}

\end{frontmatter}


\section{Introduction}\label{Introduction}

The ARCH (AutoRegressive Conditional Heteroscedastic) model has been introduced by \cite{E-ARCH} in an univariate context. Since this work a lot of extensions have been proposed. A first one has been suggested four years later, namely the GARCH (Generalised ARCH) model by \cite{B-GARCH}. This model had for goal to improve modeling by considering the past conditional variance (volatility). Their concept are based on the past conditional heteroscedasticity which depends on the past values of the return. A consequence is the volatility has the same magnitude for a negative or positive return.

Financial series have their own characteristics which are usually difficult to reproduce artificially. An important characteristic is the leverage effect which considers negative returns differently than the positive returns. This is in contradiction with the construction of the GARCH model, because it cannot consider the asymmetry. The TGARCH (Threshold GARCH) model introduced by \cite{RZ-TGARCH} improves the modeling because it considers the asymmetry since the volatility is determined by the past negative observations and the past positive observations with different weights. Various asymmetric GARCH processes are introduced in the econometrics literature, for instance the EGARCH (Exponential GARCH) and the $\log-$GARCH models (see \cite{FWZ-EGARCH-2013} who studied the asymptotic properties of the estimators for an EGARCH$(1,1)$ models).

The standard GARCH model is based on the concept that the conditional variance is a linear function of the squared past innovations whereas the TGARCH model is based on the concept of conditional standard deviation (also called volatility). In reality,  as mentioned by \cite{E-ARCH}, other formulations (or several functional forms) of the volatility can be more appropriate. Motivated by Box-Cox power transformations,  \cite{HB-PARCH} proposed a general functional form for the ARCH model, the so-called PARCH (power-transformed ARCH) models in which the conditional variance is modeled by a power $\delta$ (see also \cite{DGE}). \cite{HK-2004} extend the PARCH models to the class of asymmetric models: the power transformed asymmetric (threshold) ARCH (APARCH for short) model. \cite{pan} generalized the APARCH (APGARCH)  model by adding the past realizations of the volatility. \cite{HZ-APGARCH} studied the asymptotic properties of the APGARCH models.

When one uses an  APGARCH model on real data, we often obtain interesting results on the power. 
In fact, as we can see on the Table~\ref{table1}, the power is not necessary equal to $1$ or $2$ and is different for each series.
\begin{table}[H]
 \caption{\small{Estimation of APGARH$(1,1)$ model for real dataset. We consider the series of the returns of the daily exchange rates of the Dollar (USD), the Yen (JPY), the Pounds (GBP) and the Canadian Dollar (CAD) with respect to the Euro. The observations cover the period from January 4,
1999 to July 13, 2017 which correspond to 4746 observations. The data were obtained from the web site of the National Bank of Belgium.}}
\begin{center}
\begin{tabular}{lccccc}
\hline\hline
Exchange rates & $\phantom{\sum\limits_{i}^n}\hat\omega\phantom{\sum\limits_{i}^n}$ & $\phantom{\sum\limits_{i}^n}\hat\alpha^+\phantom{\sum\limits_{i}^n}$ & $\phantom{\sum\limits_{i}^n}\hat\alpha^-\phantom{\sum\limits_{i}^n}$ &$\phantom{\sum\limits_{i}^n}\hat\beta\phantom{\sum\limits_{i}^n}$ & $\phantom{\sum\limits_{i}^n}\hat\delta\phantom{\sum\limits_{i}^n}$ \\
\hline
USD & 0.00279 & 0.02618 & 0.04063 & 0.96978 & 1.04728 \\
JPY & 0.00740 & 0.05331 & 0.08616 & 0.93580 & 1.12923 \\
GBP & 0.00240 & 0.06078 & 0.06337 & 0.94330 & 1.41851 \\
CAD & 0.00416 & 0.04054 & 0.03111 & 0.96114 & 1.56085 \\
\hline\hline \\
\end{tabular}
\end{center}
\label{table1}
\end{table}
In the econometric applications, the univariate APGARCH framework is very restrictive.
Despite the fact that  the volatility of univariate series has been widely studied in the literature, modeling the realizations of several series is of great practical importance. When several series displaying temporal  dependencies are available, it is useful to analyze them jointly, by viewing them as the components of a vector-valued (multivariate) process. Contrarily to the class of Vector Auto Regressive Moving Average (VARMA), there is no natural extension of GARCH models for vector series, and many MGARCH (multivariate GARCH) formulations are presented in the literature (see for instance \cite{Nel-MGARCH}, \cite{EK-MGARCH}, \cite{ Engle-DCC} and  \cite{ MFS-MGARCH}). See also \cite{BLR-MGARCH}, \cite{Silvennoinen2009} and \cite{BHL-VEC-2012handbook} for recent surveys on MGARCH processes. These extensions present numerous specific problems such that identifiabily conditions, estimation.  Among the numerous specifications of MGARCH models, the most popular seems to be the
Constant Conditional Correlations (CCC) model introduced by \cite{B-MGARCH}  and extended
by \cite{jeantheau} (denoted CCC-GARCH, in the sequel).

As mentioned before, to model the dynamics of daily financial returns, we need to incorporate the leverage effect to volatility models. Many asymmetric univariate GARCH models have been considered in the literature to capture the leverage effect, extensions to the multivariate setting have not been much developed. To our knowledge, notable exceptions are \cite{MHC-MGARCH} who extend the GJR model (see \cite{GJR-GJR}) to the CCC Asymmetric GARCH. Another extension is the Generalized Autoregressive Conditional Correlation (GARCC) model proposed by \cite{MFS-MGARCH}. Recently \cite{FZ-MAPGARCH} proposed an asymmetric CCC-GARCH (CCC-AGARCH) model that includes the CCC-GARCH introduced by \cite{B-MGARCH}  and its generalization by \cite{jeantheau}. The attractiveness of the CCC-AGARCH models follows from their tractability (see \cite{FZ-MAPGARCH}). \cite{SGE-logMGARCH} proposed a general framework for the estimation and inference in univariate and multivariate $\log-$GARCH-X models (when covariates or other conditioning variables ''$X$'' are added to the volatility equation) with Dynamic Conditional Correlations of unknown form via the VARMA-X representation (Vector Auto Regressive Moving Average).

The main purpose of this paper is to introduce and study a multivariate version of the APGARCH models. In view of the results summarized in Table \ref{table1}, it appears to be inadequate to consider a unique power for all the $m$ series. Hence we propose the CCC power transformed asymmetric (threshold) GARCH (denoted CCC-APGARCH or $\underline{\delta}_0-$CCC-GARCH, where  $\underline{\delta}_0$ is  a $m-$vector of powers and $m$ is the number of series considered).  Our model includes for examples the CCC-AGARCH developed by \cite{FZ-MAPGARCH} and of course, the most classical MGARCH model, the CCC-GARCH introduced by \cite{B-MGARCH}. An important feature in this family of models is that the interpretation of the coefficients and the conditional variance is simpler and explicit. We shall give a necessary and sufficient condition for the existence of a strictly stationary solution of the proposed model and we study the problem of estimation of the CCC-APGARCH.

For the estimation of GARCH and MGARCH models, the commonly used estimation method
is the quasi-maximum likelihood estimation (QMLE for short). The asymptotic distribution of the Gaussian QMLE is obtained for a wide class of asymmetric GARCH models with exogenous covariates
by \cite{francq_thieu_2018}.  The asymptotic theory of MGARCH models are well-known in the literature. For instance  \cite{jeantheau} gave general conditions for the strong consistency of the QMLE for multivariate GARCH models. \cite{CL-NA-MGARCH} (see also \cite{HP-MGARCH}) have proved the consistency and the asymptotic normality of the QMLE for the BEKK  formulation (the acronym comes from synthesized work on multivariate models by Baba, Engle, Kraft and Kroner).  Asymptotic results were established by \cite{LM-VARMA-GARCH} for the CCC formulation of an ARMA-GARCH. See also \cite{BW-NA-multi} who studied the asymptotic behavior of the QMLE in a general class of multidimensional causal processes allowing asymmetries.
Recently, the quasi-maximum likelihood (QML) results have been established for a MGARCH with stochastic correlations by \cite{FZ-EBE-multi} under the assumption that the system is estimable equation-by-equation. See also \cite{DFL}, who proved the asymptotic properties of the QML equation-by-equation estimator of Cholesky GARCH models and time-varying conditional
betas. 
\cite{FS-EBE-logMGARCH} prove the consistency and asymptotic normality of a least squares equation-by-equation estimator of a
multivariate $\log-$GARCH-X model with Dynamic Conditional Correlations by using the VARMA-X representation.
Strong consistency and asymptotic normality of CCC-Periodic-GARCH models are established by \cite{Bibi}. The asymptotic normality of 
maximum-likelihood estimator of Dynamic Conditional Beta is proved by \cite{Engle2016}. In our context, we use the  quasi-maximum likelihood estimation. The proofs of our results are quite technical. These are adaptations of the arguments used in \cite{FZ-MAPGARCH} when the power is known, \cite{HZ-APGARCH} and \cite{pan} when the power is unknown. We strength the fact that new techniques and arguments are applied in the proof of of the identifiability (see Subsection \ref{new-proof1} and in the proof of the invertibility of the Fisher information  matrix (see Subsection \ref{new-proof2}). 

 This paper is organized as follows. In Section \ref{model} we introduce the CCC-APGARCH model and show that it
includes some class of (M)GARCH models. We established the strict stationarity condition and we give an identifiability condition. Section \ref{QMLEdeltaConnu} is devoted to the asymptotic properties of the quasi-maximum likelihood estimation when the power $\underline{\delta}_0$ is known. In Section \ref{QMLEdeltaInconnu}, we consider the estimation of $\underline{\delta}_0$.
Wald test is developed in Section \ref{lineartest} in order to test the classical MGARCH model against a class of asymetric MGARCH models. The test can also be used to test the equality between the components of $\underline{\delta}_0$. Simulation studies and an illustrative application on real data are presented in  Section~\ref{NumIllust} and we provide a conclusion in Section \ref{conclusion}.  The proofs of the main results are collected in the Appendix \ref{app}.
\section{Model and strict stationarity condition}\label{model}
In all this work, we use the following notation $\underline{u}^{\underline{d}} := (u_1^{d_1},\ldots, u_m^{d_m})'$ for $\underline{u},\underline{d}\in\mathbb R^m$.
\subsection{Model presentation}

The $m$-dimensional process $\underline{\varepsilon}_t = (\varepsilon_{1,t}, \ldots, \varepsilon_{m,t})'$ is called a CCC-APGARCH$(p,q)$ if it verifies
\begin{equation}\label{MAPGARCH}
\left\{
\begin{aligned}
&\underline{\varepsilon}_t = H_t^{1/2}\eta_t, \\
&H_t = D_t R_0 D_t,\qquad D_t = \mbox{diag}(\sqrt{h_{1,t}},\dots, \sqrt{h_{m,t}}),\\
&\underline{h}_t^{\underline{\delta}_0/2} = \underline{\omega}_0 + \sum\limits_{i=1}^q\left\{ A_{0i}^+(\underline{\varepsilon}_{t-i}^+)^{\underline{\delta}_0/2} + A_{0i}^-(\underline{\varepsilon}_{t-i}^-)^{\underline{\delta}_0/2}\right\} + \sum\limits_{j=1}^p B_{0j} \underline{h}_{t-j}^{\underline{\delta}_0/2},
\end{aligned}
\right.	
\end{equation}
where $\underline{h}_t = (h_{1,t}, \ldots, h_{m,t})'$ and with $x^+ = \max(0, x)$ and $x^- = \min(0, x)$
\[ \underline{\varepsilon}_t^+ = \left(\{\varepsilon_{1,t}^+\}^2, \ldots, \{\varepsilon_{m,t}^+\}^2\right)'\qquad \underline{\varepsilon}_t^- = \left(\{-\varepsilon_{1,t}^-\}^2, \ldots, \{-\varepsilon_{m,t}^-\}^2\right)',\]
$\underline{\omega_0}$ and $\underline{\delta_0}$ are vectors of size $m \times 1$ with strictly positive coefficients, $A_{0i}^+, A_{0i}^-$ and $B_{0j}$ are matrices of size $m\times m$ with positive coefficients and $R_0$ is a correlation matrix. The parameters of the model are the coefficients of the vectors $\underline{\omega}_0$, $\underline{\delta}_0$, the coefficients of the matrices $A_{0i}^+, A_{0i}^-, B_{0j}$ and the coefficients in the lower triangular part excluding the diagonal of the matrix $R_0$. The number of unknown parameters is
\[s_0 = 2m + m^2(2q+p) + \dfrac{m(m-1)}{2}.\]
The innovation process $(\eta_t)_t$ is a vector of size $m\times 1$ and satisfies the assumption:\\[2mm]
\indent $\textbf{A0} : (\eta_t)$ is an  independent and identically distributed (iid for short)  sequence of variables on $\mathbb{R}^m$ with identity covariance matrix and $\mathbb{E}[\eta_t] = 0$.\\[2mm]
With this assumption, the  matrix $H_t$ is interpreted as the conditional variance (volatility) of $\underline{\varepsilon}_t$.
The representation \eqref{MAPGARCH} includes various MGARCH models.
For instance, if we assume that $A_{0i}^+ = A_{0i}^-$, the model \eqref{MAPGARCH} can be viewed as a multivariate extension and generalization of the
 the so-called PARCH (power-transformed ARCH) models introduced by \cite{HB-PARCH} (denoted CCC-PGARCH, in the sequel). Moreover, if we also fixed the power vector $\underline{\delta}_0=(2,\dots,2)'$ we obtain the CCC-GARCH$(p,q)$ model proposed by \cite{jeantheau}. Now, if we assume that $A_{0i}^+\neq A_{0i}^-$, the model CCC-AGARCH$(p,q)$ of \cite{FZ-MAPGARCH} is retrieved. If we also fixed $m=1$ (the univariate case, we retrieve the APGARCH introduced by \cite{pan}.
%
%

As remarked by \cite{FZ-MAPGARCH}, the interest of the APGARCH model \eqref{MAPGARCH} is to
consider the leverage effect observed for most of financial series. It is well known that a negative return has more effect on the volatility as a positive return. In fact, for the same magnitude,  the volatility  $\underline{h}^{\underline{\delta}_0/2}_t$ increases more if the return is negative than if it is positive.
In general, the estimation of the coefficients of the model suggests that one may think that $A_{0i}^+ < A_{0i}^-$ for some $i>0$ component by component.

\cite{B-MGARCH} introduced the CCC-GARCH model with the assumption that the coefficients matrices  $A_{0i}^+ = A_{0i}^-$ and $B_{0j}$ are diagonal. This assumption implies that the conditional variance $h_{k,t}$ of the $k$-th component of $\underline{\varepsilon}_t$ depends only on their own past values and not on the past values of the other components. By contrast, in the model \eqref{MAPGARCH}  (including the CCC-GARCH) the conditional variance $h_{k,t}$ of the $k$-th component of $\underline{\varepsilon}_t$ depends not only on its past values but also on the past values of the other components. For this reason, Model \eqref{MAPGARCH}  is referred to as the Extended CCC model by \cite{HT-2004-extendedCCC-GARCH}.
\subsection{Strict stationarity condition}

A sufficient condition for strict stationarity of the CCC-AGARCH($p,q$) is given by \cite{FZ-MAPGARCH}. For the model \eqref{MAPGARCH} the strict stationarity condition is established in the same way. For that sake we rewrite the  first equation of \eqref{MAPGARCH} as
\begin{equation}\label{eta-tilde}
 \underline{\varepsilon}_t = D_t\tilde{\eta}_t,\qquad \mbox{ where } \tilde{\eta}_t = (\tilde{\eta}_{1,t},\ldots, \tilde{\eta}_{m,t}) = R_0^{1/2}\eta_t.
\end{equation}
Using the third equation of model \eqref{MAPGARCH}, we may write
\begin{equation}\label{upsilon}
(\underline{\varepsilon}_t^\pm)^{\underline{\delta_0}/2} = (\Upsilon_t^{\pm,(\delta_0)})\underline{h}_t^{\underline{\delta_0}/2}, \mbox{ with } \Upsilon_t^{\pm,(\delta_0)} = \mbox{diag}\left((\pm\tilde{\eta}_{1,t}^\pm)^{\delta_{0,1}},\ldots, (\pm\tilde{\eta}_{m,t}^\pm)^{\delta_{0,m}}\right),
\end{equation}
where $\delta_0=(\delta_{0,1},\dots,\delta_{0,m})'$. To study the strict stationarity condition, we introduce the matrix expression for the model \eqref{MAPGARCH}
\begin{equation*}
\underline{z}_t = \underline{b}_t + C_t \underline{z}_{t-1},\label{matricielle}
\end{equation*}
where
\begin{equation*}
\underline{z}_t = \left(\left\{\left(\underline{\varepsilon}_t^{+}\right)^{\underline{\delta_0}/2}\right\}', \dots, \left\{\left(\underline{\varepsilon}_{t-q+1}^{+}\right)^{\underline{\delta_0}/2}\right\}', \left\{\left(\underline{\varepsilon}_t^{-}\right)^{\underline{\delta_0}/2}\right\}', \dots, \left\{\left(\underline{\varepsilon}_{t-q+1}^{-}\right)^{\underline{\delta_0}/2}\right\}', \left\{\underline{h}_t^{\underline{\delta_0}/2}\right\}', \dots, \left\{\underline{h}_{t-p+1}^{\underline{\delta_0}/2}\right\}'\right)',\\
\end{equation*}
\begin{equation*}
\underline{b}_t = \left(\left\{\left(\Upsilon_t^{+,(\delta_0)}\right)\underline{\omega}_0\right\}', 0_{m(q-1)}', \left\{\left(\Upsilon_t^{-,(\delta_0)}\right)\underline{\omega}_0\right\}', 0_{m(q-1)}', \left\{\underline{\omega}_0\right\}', 0_{m(p-1)}'\right)'
\end{equation*}
and
\begin{equation}
C_t = \left(
\begin{tabular}{ccc}
$\Upsilon_t^{+,(\delta_0)}A_{01 : q}^+$ & $\Upsilon_t^{+,(\delta_0)}A_{01 : q}^-$ & $\Upsilon_t^{+,(\delta_0)}B_{01:p}$\\
$I_{m(q-1)}$ & \multicolumn{2}{c}{$0_{m(q-1) \times m(p+q+1)}$}\\
$\Upsilon_t^{-,(\delta_0)}A_{01 : q}^+$ & $\Upsilon_t^{-,(\delta_0)}A_{01 : q}^-$ & $\Upsilon_t^{-,(\delta_0)}B_{01:p}$\\
$0_{m(q-1)\times mq}$ & $I_{m(q-1)}$ & $0_{m(q-1) \times m(p-1)}$ \\
$A_{01:q}^+$ & $A_{01:q}^-$ & $B_{01:p}$\\
\multicolumn{2}{c}{$0_{m(p+q+1)\times m(q-1)}$} & $I_{m(p-1)}$ \\
\end{tabular}
\right).
\end{equation}
We have denoted  $A_{01:q}^+= (A_{01}^+\ldots A_{0q}^+)$, $A_{01:q}^- = (A_{01}^-\ldots A_{0q}^-)$ and $B_{01:p} = (B_{01}\ldots B_{0p})$ (they are $m\times qm$ and  $m\times pm$ matrices). The matrix $C_t$ is of size $(p+2q)m \times (p+2q)m$.\\
Let $\gamma(\textbf{C}_0)$ the top Lyapunov exponent of the sequence $\textbf{C}_0 = \{C_t, t\in \mathbb{Z}\}$. It is defined by
\begin{equation*}
\gamma(\textbf{C}_0) := \lim_{t\to+\infty} \dfrac1t \mathbb{E}\left[\log\Vert C_{t}C_{t-1}\dots C_{1}\Vert\right] = \inf_{t \geq 1}\dfrac1t \mathbb{E}\left[\log\Vert C_{t}C_{t-1}\dots C_{1}\Vert\right].
\end{equation*}
Now we can state the following results. Their proofs are the same than the one in \cite{FZ-MAPGARCH} so they are omitted.
\begin{thm}\label{Strict}\  {(Strict stationarity)}\\
A necessary and sufficient condition for the existence of a strictly stationary and non anticipative solution process to model \eqref{MAPGARCH} is $\gamma(\textbf{C}_0) < 0$. 

When $\gamma(\textbf{C}_0) < 0$, the stationary and non anticipative solution is unique and ergodic.
\end{thm}

The two following corollaries are consequences of the necessary condition for strict stationarity. For  $A$ a square matrix, $\rho(A)$ denotes its spectral radius (i.e. the greatest modulus of its eigenvalues).
\begin{cor}\ \\
\label{cor1}
Let $\mathcal{B}_0$ be  the matrix polynomial defined by $\mathcal{B}_0(z) = I_m - zB_{01} - \cdots - z^pB_{0p}$ for $z \in \mathbb{C}$ and define
\[
\textbf{B}_0 = \left(
\begin{tabular}{cc}
\multicolumn{2}{c}{$B_{01:p}$}\\
$I_{(p-1)m}$ & $0_{(p-1)m\times 1}$\\
\end{tabular}
\right).
\]
Then, if $\gamma(\textbf{C}_0) < 0$ the following equivalent properties hold:
\begin{enumerate}[$\qquad (i)$]
	\item The roots of $\det(\mathcal{B}_0(z))$ are outside the unit disk,
	\item $\rho(\textbf{B}_0) < 1$.
\end{enumerate}
\end{cor}

\begin{cor}\label{cor2}\  \\
Suppose $\gamma(\textbf{C}_0) < 0$. Let $\underline{\varepsilon}_t$ be the strictly stationary and non anticipative solution of model \eqref{MAPGARCH}. There exists $s>0$ such that $\mathbb{E} \Vert \underline{h}_t^{\underline{\delta_0}/2} \Vert^s < \infty$ and $\mathbb{E} \Vert \underline{\varepsilon}_t^{\underline{\delta_0}/2} \Vert^s < \infty$.\\
\end{cor}

\subsection{Identifiability condition}

In this part, we are  interested in the identifiability condition to ensure the uniqueness of the parameters in the CCC-APGARCH representation. This is a crucial step before the estimation.

The parameter $\nu$ is defined by
\begin{equation*}
\nu := (\underline{\omega}',  {\alpha_{1}^{+}} ', \ldots, {\alpha_q^+}', {\alpha_1^-}', \ldots, {\alpha_q^-}', \beta'_{1}, \ldots, \beta'_{p}, \underline{\tau}',\rho')',
\end{equation*}
where $\alpha_i^+$ and $\alpha_i^-$ are defined by $\alpha_i^\pm = \mbox{vec}(A_i^\pm)$ for $i=1, \ldots, q$, $\beta_j = \mbox{vec}(B_j)$ for $j = 1,\ldots, p$,  $\underline{\tau}'=(\tau_1,...,\tau_m)$ is the vector of powers and $\rho = (\rho_{21},\ldots \rho_{m1}, \rho_{32}, \ldots, \rho_{m2}, \ldots, \rho_{mm-1})'$ such that the $\rho_{ij}$'s are the components of the matrix $R$.  
The  parameter $\nu$ belongs to the parameter space
\begin{equation*}
\Delta \subset ]0,+\infty[^{m} \times [0,\infty[^{m^2(2q+p)}\times]0,+\infty[^{m} \times ]-1,1[^{m(m-1)/2}.
\end{equation*}
The unknown true parameter value is denoted by
\begin{equation*}
\nu_0 := (\underline{\omega}_0', {\alpha_{01}^+}', \ldots, {\alpha_{0q}^+}', {\alpha_{01}^-}', \ldots, {\alpha_{0q}^-}', {\beta_{01}}',\ldots, {\beta_{0p}}', \underline{\delta}_0',\rho_0')'.
\end{equation*}
We adopt the following notation. For a matrix $A$ (which has to be seen as a parameter of the model), we write $A_0$ when the coefficients of the matrix are evaluated in the true value $\nu = \nu_0$. \\
Let $\mathcal{A}^+(L) = \sum_{i=1}^q A_i^+L^i$, $\mathcal{A}^-(L) = \sum_{i=1}^q A_i^-L^i$ and $\mathcal{B}(L) = I_m - \sum_{j=1}^pB_jL^j$ where $L$ is the backshift operator. By convention $\mathcal{A}^\pm = 0$ if $q=0$ and $\mathcal{B}(L) = I_m$ if $p=0$.\\
If the roots of $\mbox{det}(\mathcal{B}_0(L)) = 0$ are outside the unit disk, we have from $\mathcal{B}_0(L) \underline{h}_t^{\underline{\delta_0}/2} = \underline{\omega}_0 + \mathcal{A}_0^+(L) (\underline{\varepsilon}_t^+)^{\underline{\delta_0}/2} + \mathcal{A}_0^-(L) (\underline{\varepsilon}_t^-)^{\underline{\delta_0}/2}$ the compact expression:
\begin{equation}\label{IC}
\underline{h}_t^{\underline{\delta_0}/2} = \mathcal{B}_0(1)^{-1}\underline{\omega}_0 + \mathcal{B}_0(L)^{-1}\mathcal{A}_0^+(L) (\underline{\varepsilon}_t^+)^{\underline{\delta_0}/2} + \mathcal{B}_0(L)^{-1}\mathcal{A}_0^-(L) (\underline{\varepsilon}_t^-)^{\underline{\delta_0}/2}.
\end{equation}
The parameter $\nu_0$ is said to be identifiable if \eqref{IC} does not hold true when $\nu_0$ is replaced by $\nu \neq \nu_0$ belonging to $\Delta$.

The assumption that  the polynomials $\mathcal{A}_0^+, \mathcal{A}_0^+$ and $\mathcal{B}_0$ have no common roots is not sufficient to consider that there is not another triple $(\mathcal{A}_0^+, \mathcal{A}_0^+, \mathcal{B}_0)$ such that
\begin{equation}\label{identifiability}
\mathcal{B}^{-1}\mathcal{A}^+ = \mathcal{B}_0^{-1}\mathcal{A}_0^+
\quad \mbox{ and }\quad  \mathcal{B}^{-1}\mathcal{A}^- = \mathcal{B}_0^{-1}\mathcal{A}_0^-.
\end{equation}
This condition is equivalent as the existence of an operator $U(B)$ such that
\[\mathcal{A}^+(L) = U(L)\mathcal{A}_0^+(L),\quad \mathcal{A}^-(L) = U(L)\mathcal{A}_0^-(L)\quad \mbox{ and }\quad  \mathcal{B}(L) = U(L)\mathcal{B}_0(L).\]
The matrix $U(L)$ is unimodular if $\det(U(L))$ is a constant not equal to zero. If the common factor to both polynomials is unimodular,
\[P(L) = U(L)P_1(L), \quad Q(L) = U(L)Q_1(L) \Rightarrow \det(U(L)) = c,\]
the polynomials $P(L)$ and $Q(L)$ are left-coprimes.\\
But in the vectorial case, suppose that $\mathcal{A}_0^+$, $\mathcal{A}_0^-$ and $\mathcal{B}_0$ are left-coprimes is not sufficient to consider that \eqref{identifiability} have no solution for $\nu \neq \nu_0$ (see \cite{FZ-MAPGARCH}).\\
To obtain a mild condition, for any column $i$ of the matrix operators $\mathcal{A}_0^+$, $\mathcal{A}_0^-$ and $\mathcal{B}_0$, we denote by  $q_i^+(\nu)$, $q_i^-(\nu)$, and $p_i(\nu)$ their maximal degrees. We suppose that the maximal values of the orders are imposed:
\begin{equation}\label{condition}
\forall \nu \in \Delta, \forall i =1, \ldots, m, \quad q_i^+(\nu) \leq q_i^+,\quad q_i^-(\nu) \leq q_i^-,\quad \mbox{ and }\quad p_i(\nu) \leq p_i
\end{equation}
where $q_i^+\leq q, q_i^-\leq q$ and $p_i\leq p$ are fixed integers. \\
We denote $a_{q_i^+}^+(i)$ the column vector of the coefficients $L^{q_i^+}$ , $a_{q_i^-}^-(i)$ the column vector of the coefficients $L^{q_i^-}$ in the column $i$ of $\mathcal{A}_0^+$, respectively $\mathcal{A}_0^-$ and $b_{p_i}(i)$ the column vector of the coefficients $L^{p_i}$ in the column $i$ of $\mathcal{B}_0$.\\
\begin{prop}\label{identif}\  {(Identifiability condition)}\\
If the matrix polynomials $\mathcal{A}_0^+(L), \mathcal{A}_0^-(L)$ and $\mathcal{B}_0(L)$ are left-coprime, $\mathcal{A}_0^+(1) + \mathcal{A}_0^-(1) \neq 0$ and if the matrix
\[M(\mathcal{A}_0^+(L), \mathcal{A}_0^-(L), \mathcal{B}_0(L)) = \left[a_{q_1^+}^+(1) \ldots a_{q_m^+}^+(m) a_{q_1^-}^-(1) \ldots a_{q_m^-}^-(m) b_{p_1}(1) \ldots b_{p_m}(m)\right]\]
has full rank $m$, under the constraints \eqref{condition} with $q_i^+ = q_i^+(\nu_0), q_i^- = q_i^-(\nu_0)$ and $p_i = p_i(\nu_0)$ for any value of $i$, then
\begin{equation*}
\left\{
\begin{aligned}
\mathcal{B}(L)^{-1}\mathcal{A}^+(L) = \mathcal{B}_0(L)^{-1}\mathcal{A}_0^+(L)\\
\mathcal{B}(L)^{-1}\mathcal{A}^-(L) = \mathcal{B}_0(L)^{-1}\mathcal{A}_0^-(L)
\end{aligned}
\right. \Rightarrow (\mathcal{A}^+, \mathcal{A}^-, \mathcal{B}) = (\mathcal{A}_0^+, \mathcal{A}_0^-, \mathcal{B}_0).
\end{equation*}
\end{prop}
\  \\
\textbf{Proof of Proposition \ref{identif}}\\
The proof of the identifiability condition is identical as \cite{FZ-MAPGARCH} in the case of the CCC-AGARCH model.
\section{Estimation when the power is known}\label{QMLEdeltaConnu}
In this section, we assume that $\underline{\delta}_0$ is known. We write $\underline{\delta}_0=\underline{\delta}$ in order to simplify the writings.
For the estimation of GARCH and MGARCH models, the commonly used
estimation method is the QMLE, which can also be viewed as a
nonlinear least squares estimation (LSE).
The QML method is particularly relevant for GARCH models because
it provides consistent and asymptotically normal estimators for strictly stationary GARCH processes
under mild regularity conditions. For example, no moment assumptions on the observed process are required (see for instance \cite{FZ-bernoulli} or \cite{FZ-2010}).

As remarked in Section~\ref{model}, some particular cases of Model \eqref{MAPGARCH} are obtained for $\underline{\delta}=(2,\dots,2)'$: the CCC-AGARCH introduced by \cite{FZ-MAPGARCH} and the CCC model introduced by \cite{HT-2004-extendedCCC-GARCH} with $A_{0i}^+ = A_{0i}^-$. This section provides asymptotic results which can, in particular, be applied to those models.
\subsection{QML estimation}
The procedure of estimation and the asymptotic properties are similar to those of the model CCC-AGARCH introduced by \cite{FZ-MAPGARCH}.

The parameters are the coefficients of the vector $\underline{\omega}_0$, the matrices $A_{0i}^+, A_{0i}^-$ and $B_{0j}$, and the coefficients of the lower triangular part without the diagonal of the correlation matrix $R_0$. The number of parameters is
\begin{equation*}
s_0 = m + m^2(p+2q) + \dfrac{m(m-1)}{2}.
\end{equation*}
The goal is to estimate the $s_0$ coefficients of the model \eqref{MAPGARCH}. In this section, we note the parameter
\[\theta := (\underline{\omega}', {\alpha_1^+}', \ldots, {\alpha_q^+}', {\alpha_1^-}', \ldots, {\alpha_q^-}', {\beta_1}',\ldots, {\beta_p}',\rho')',\]
where $\alpha_i^+$ and $\alpha_i^-$ are define by $\alpha_i^\pm = \mbox{vec}(A_i^\pm)$ for $i=1, \ldots, q,
\beta_j = \mbox{vec}(B_j)$ for $j = 1,\ldots, p$ and $\rho = (\rho_{21},\ldots \rho_{m1}, \rho_{32}, \ldots, \rho_{m2}, \ldots, \rho_{mm-1})'$.
The parameter $\theta$ belongs to the parameter space
\begin{equation*}
\Theta \subset ]0,+\infty[^{m} \times [0,\infty[^{m^2(2q+p)}\times ]-1,1[^{m(m-1)/2}.
\end{equation*}
The unknown true value of the parameter  is denoted by
\begin{equation*}
\theta_0 := (\underline{\omega}_0', {\alpha_{01}^+}', \dots, {\alpha_{0q}^+}', {\alpha_{01}^-}', \dots, {\alpha_{0q}^-}', {\beta_{01}}',\dots, {\beta_{0p}}',\rho_0')'.
\end{equation*}
The determinant of a square matrix $A$ is denoted by $\det(A)$ or $\vert A \vert$.

Let $(\underline{\varepsilon}_1, \ldots, \underline{\varepsilon}_n)$ be a realization of length $n$ of the unique non-anticipative strictly stationary solution $(\underline{\varepsilon}_t)$ of Model \eqref{MAPGARCH}. Conditionally to nonnegative initial values $\underline{\varepsilon}_0, \ldots, \underline{\varepsilon}_{1-q}, \underline{\tilde{h}}_0, \ldots, \underline{\tilde{h}}_{1-p}$, the Gaussian quasi-likelihood writes
\begin{equation*}
L_n(\theta) = L_n(\theta ; \underline{\varepsilon}_1,\dots, \underline{\varepsilon}_n) = \prod_{t=1}^n \dfrac{1}{(2\pi)^{m/2}\vert\tilde{H}_t\vert^{1/2}}\exp\left(-\dfrac12\underline{\varepsilon}_t' \tilde{H}_t^{-1}\underline{\varepsilon}_t\right),
\end{equation*}
where the $\tilde{H}_t$ are recursively defined, for $t \geq 1$, by
\begin{equation*}
\left\{
\begin{aligned}
\tilde{H}_t &= \tilde{D}_tR\tilde{D}_t, \qquad \tilde{D}_t = \mbox{diag}\left(\sqrt{\tilde{h}_{1,t}},\dots, \sqrt{\tilde{h}_{m,t}}\right)\\
\underline{\tilde{h}}_t^{\underline{\delta}/2} &:=\underline{\tilde{h}}_t^{\underline{\delta}/2}(\theta)=\underline{\omega} + \sum\limits_{i=1}^q A_i^+(\underline{\varepsilon}_{t-i}^+)^{\underline{\delta}/2} + A_i^-(\underline{\varepsilon}_{t-i}^-)^{\underline{\delta}/2} + \sum\limits_{j=1}^p B_j \underline{\tilde{h}}_{t-j}^{\underline{\delta}/2}.
\end{aligned}
\right.
\end{equation*}
A quasi-likelihood estimator of $\theta$ is defined as any measurable solution $\hat{\theta}_n$ of
\begin{equation}\label{QMLE}
\hat{\theta}_n = \underset{\theta \in \Theta}{\arg\max}\  L_n(\theta) = \underset{\theta \in \Theta}{\arg\min}\  \tilde{\mathcal{L}}_n(\theta),
\end{equation}
where
\begin{equation*}
\tilde{\mathcal{L}}_n(\theta) = \dfrac1n \sum\limits_{t=1}^n \tilde{l}_t, \qquad \tilde{l}_t = \tilde{l}_t(\theta) = \underline{\varepsilon}_t' \tilde{H}_t^{-1}\underline{\varepsilon}_t + \log \vert \tilde{H}_t\vert.
\end{equation*}
\subsection{Asymptotic properties}
To establish the strong consistency, we need the following assumptions borrowed from \cite{FZ-MAPGARCH}:\\
	\indent $\textbf{A1}:$ $\theta_0 \in \Theta$ and $\Theta$ is compact,\\
	\indent $\textbf{A2}:$ $\gamma(\textbf{C}_0) < 0$ and $\forall \theta \in \Theta, \det(\mathcal{B}(z)) = 0 \Rightarrow \vert z\vert > 1$,\\
	\indent $\textbf{A3}:$ For $i=1,\ldots, m$ the distribution of $\tilde{\eta}_{it}$ is not concentrated on 2 points and $\mathbb{P}(\tilde{\eta}_{it}>0)\in  (0,1)$.\\
	\indent $\textbf{A4}:$ if $p>0, \mathcal{A}_0^+(1) + \mathcal{A}_0^-(1) \neq 0,  \mathcal{A}_0^+(z), \mathcal{A}_0^-(z)$ and $\mathcal{B}_0(z)$ are left-coprime and the matrix $M(\mathcal{A}_0^+, \mathcal{A}_0^-, \mathcal{B}_0)$ has full rank $m$.\\
	If the space $\Theta$ is constrained by \eqref{condition}, Assumption $\textbf{A4}$ can be replaced by the more general condition
	\indent $\textbf{A4'}:$ $\textbf{A4}$ with $M(\mathcal{A}_0^+, \mathcal{A}_0^-, \mathcal{B}_0)$ replaced by $[A_{0q}^+A_{0q}^-B_{0p}]$.\\
	\indent $\textbf{A5}:$ $R$ is a positive-definite correlation matrix for all $\theta \in \Theta$.\\
To ensure the strong consistency of the QMLE, a compactness assumption is required (i.e \textbf{A1}).	
The assumption \textbf{A2} makes reference to the condition of strict stationarity for the model (\ref{MAPGARCH}). This assumption implies that for the true parameter $\theta_0$, Model (\ref{MAPGARCH}) admits a  strictly stationary solution but is less restrictive concerning the other values $\theta\in\Theta$. The second part of Assumption \textbf{A2} implies that the roots of $\det(\mathcal{B}(z))$ are outside the unit disk. Assumptions \textbf{A3} and \textbf{A4} or \textbf{A4'}
are made for identifiability reasons. In particular $\mathbb{P}(\tilde{\eta}_{it}>0)\in  (0,1)$ ensures that the process $(\underline{\varepsilon}_{it})$ (for $i=1,\dots, m$) takes positive and negative values with a positive probability (if, for instance, the $(\underline{\varepsilon}_{it})$ were
a.s. positive, the parameters $\alpha_{0j}^-$ for $j=1,\dots,q$ could not be identified).

We are now able to state the following strong consistency theorem.
\begin{thm}
\label{le-joli-label-de-yacouba}
Let $(\hat{\theta}_n)$ be a sequence of QMLE satisfying \eqref{QMLE}. Then, under $\textbf{A0},\dots,\textbf{A5}$ or $\textbf{A0},\dots,\textbf{A4'}$ and $\textbf{A5}$, we have
$\hat{\theta}_n \longrightarrow \theta_0,$ almost surely as $n\to \infty.$
\end{thm}
The proof of this result is postponed to Subsection \ref{A2} in the Appendix  \ref{app}. 

It will be useful to approximate the sequence $(\tilde{l}_t)$ by an ergodic and stationary
sequence. Under Assumption $\textbf{A2}$ there exists a strictly stationary, non anticipative and ergodic solution $(h_t)_t = (h_t(\theta))_t$ of
\begin{align}\label{theorique}
\underline{h}_t^{\underline{\delta}/2} &:=\underline{h}_t^{\underline{\delta}/2}(\theta)=\underline{\omega} + \sum\limits_{i=1}^q A_i^+(\underline{\varepsilon}_{t-i}^+)^{\underline{\delta}/2} + A_i^-(\underline{\varepsilon}_{t-i}^-)^{\underline{\delta}/2} + \sum\limits_{j=1}^p B_j \underline{h}_{t-j}^{\underline{\delta}/2}.
\end{align}
We denote $D_t=D_t(\theta) = \mbox{diag}\left(\sqrt{{h}_{1,t}},\dots, \sqrt{{h}_{m,t}}\right)$ and $H_t={H}_t(\theta) = {D}_t(\theta)R{D}_t(\theta)$ and we define
\begin{equation*}
{\mathcal{L}}_n(\theta) = \dfrac1n \sum\limits_{t=1}^n {l}_t, \qquad {l}_t = {l}_t(\theta) = \underline{\varepsilon}_t' {H}_t^{-1}\underline{\varepsilon}_t + \log \vert {H}_t\vert.
\end{equation*}
To establish the asymptotic normality, the following additional assumptions are required:\\
	\indent $\textbf{A6}:$ $\theta_0\in \stackrel{\circ}{\Theta}$, where $\stackrel{\circ}{\Theta}$ is the interior of $\Theta$.\\
	\indent $\textbf{A7}:$ $\mathbb{E}\Vert \eta_t\eta_t'\Vert^2 < \infty$.\\	
Assumption \textbf{A6} prevents the situation where certain components of $\theta_0$ are equal to zero (more precisely the coefficients of the matrices in our model). One refers to Section 8.2 of \cite{FZ-2010} and \cite{pedersen} for a discussion on this topic. 

The second main result of this section is the following asymptotic normality theorem.
\begin{thm}
\label{AN-connu}
Under the assumptions of Theorem~\ref{le-joli-label-de-yacouba} and \textbf{A6}--\textbf{A7}, when $n\to\infty$, we have
\[\sqrt{n}(\hat{\theta}_n - \theta_0) \overset{\mathcal{L}}{\longrightarrow}\mathcal{N}(0, J^{-1}IJ^{-1}),\]
where $J$ is a positive-definite matrix and $I$ is a positive semi-definite matrix, defined by
\[I :=I(\theta_0)= \mathbb{E} \left[ \dfrac{\partial l_t(\theta_0)}{\partial \theta} \dfrac{\partial l_t(\theta_0)}{\partial \theta'}\right], \qquad
J:=J(\theta_0) = \mathbb{E} \left[ \dfrac{\partial^2l_t(\theta_0)}{\partial\theta\partial\theta'}\right].\]
\end{thm}
The proof of this result is postponed to Subsection \ref{anconnu} in the Appendix  \ref{app}. 
\begin{remi}\label{rem-i}
In the one dimensional case (when $m=1$), we have $I=(\mathbb E \eta_t^4 -1 )J$ with $$J=
\frac{4}{\underline{\delta}^2_0}\mathbb E \left ( \frac{1}{\underline{h}_t^{\underline{\delta}_0}(\theta_0)} \frac{\partial \underline{h}_t^{\underline{\delta}_0/2}(\theta_0)}{\partial \theta}
\frac{\partial \underline{h}_t^{\underline{\delta}_0/2}(\theta_0)}{\partial \theta'} \right )
 .$$
This expression is in accordance with the one of Theorem 2.2 in \cite{HZ-APGARCH}.
\end{remi}

%
\section{Estimation when the power is unknown}\label{QMLEdeltaInconnu}
In this section, it is assumed that $\underline{\delta}_0$ is unknown.
\subsection{QML estimation}
Now we consider the case when the power is unknown. Thus we consider joint estimation of $\theta$ and $\underline{\delta}$. In practice $\underline{\delta}_0$ is difficult to identified, as it was remarked in \cite{HZ-APGARCH} for the APGARCH model. In line with \cite{HZ-APGARCH} we make the following assumption which is used to ensure that $ \underline{\delta}_0$ is identified.  One refers to Remark 3.2 of \cite{HZ-APGARCH} for a discussion on how the following assumption differs from the one given in \cite{pan}.  \\[2mm]
\indent $\textbf{A8}:$ $\eta_t$ has a positive density on some neighbourhood of zero.\\[2mm]
%
To define the QML estimator of $\nu$, we replace $H_t$ by $\mathcal{H}_t$ in the expression of the criterion defined in \eqref{QMLE} and we obtain recursively $\tilde{\mathcal{H}}_t$, for $t \geq 1$,
\begin{equation*}
\left\{
\begin{aligned}
\tilde{\mathcal{H}}_t &= \tilde{D}_tR\tilde{D}_t, \qquad \tilde{D}_t = \mbox{diag}\left(\sqrt{\tilde{h}_{1,t}},\dots, \sqrt{\tilde{h}_{m,t}}\right)\\
\underline{\tilde{h}}_t &:=\underline{\tilde{h}}_t(\nu)= \left(\underline{\omega} + \sum\limits_{i=1}^q A_i^+(\underline{\varepsilon}_{t-i}^+)^{\underline{\tau}/2} + A_i^-(\underline{\varepsilon}_{t-i}^-)^{\underline{\tau}/2} + \sum\limits_{j=1}^p B_j \underline{\tilde{h}}_{t-j}^{\underline{\tau}/2}\right)^{2/\underline{\tau}}.
\end{aligned}
\right.
\end{equation*}
A quasi-maximum likelihood estimator of $\nu$ is defined as any mesurable solution $\hat{\nu}_n$ of
\begin{equation}\label{QMLE un}
\hat{\nu}_n = \underset{\nu \in \Delta}{\arg\min}\  \tilde{\mathcal{L}}_n(\nu),
\end{equation}
where
\[\tilde{\mathcal{L}}_n(\nu) = \dfrac1n \sum\limits_{t=1}^n \tilde{l}_t, \qquad \tilde{l}_t = \tilde{l}_t(\nu) = \underline{\varepsilon}_t' \tilde{\mathcal{H}}_t^{-1}\underline{\varepsilon}_t + \log \vert \tilde{\mathcal{H}}_t\vert.\]
\subsection{Asymptotic properties}
To establish the consistency and the asymptotic normality, we need some assumptions similar to those we assumed when the power is known. We will assume $\textbf{A1},\dots,\textbf{A6}$ with the parameter $\theta$ which is replaced by $\nu$ and the space parameter $\Theta$ is replaced by $\Delta$.
%
%
\begin{thm}\  {(Strong consistency)}\label{theorem31}\\
Let $(\hat{\nu}_n)$ be a sequence of QMLE satisfying \eqref{QMLE}. Then, under $\textbf{A0},\dots,\textbf{A5}$
or $\textbf{A0},\dots,\textbf{A4'}$, $\textbf{A5}$ and $\textbf{A8}$, we have
$\hat{\nu}_n \longrightarrow \nu_0,$ almost surely when $n\to +\infty.$
\end{thm}
The proof of this result is postponed to Subsection \ref{A5} in the Appendix \ref{app}.
%
\begin{thm}\  {(Asymptotic normality)}\label{AN-inconnu-c-est-le-label-de-bruno}\\
Under the assumptions of Theorem~\ref{theorem31}, \textbf{A6}, \textbf{A7} and \textbf{A8}, when $n\to\infty$, we have
\[\sqrt{n}(\hat{\nu}_n - \nu_0) \overset{\mathcal{L}}{\longrightarrow}\mathcal{N}(0, J^{-1}IJ^{-1}),\]
where $J$ is positive-definite matrix and $I$ is a positive semi-definite matrix, defined by
\[I :=I(\nu_0)= \mathbb{E} \left[ \dfrac{\partial l_t(\nu_0)}{\partial \nu} \dfrac{\partial l_t(\nu_0)}{\partial \nu'}\right], \qquad
J:=J(\nu_0) = \mathbb{E} \left[ \dfrac{\partial^2l_t(\nu_0)}{\partial\nu\partial\nu'}\right].\]
\end{thm}
The proof of this result is postponed to Subsection \ref{aninconnu} in the Appendix \ref{app}.
\begin{remi}
As in Remark \ref{rem-i}, in the one dimensional case, we have $I=(\mathbb E \eta_t^4 -1 )J$ with $$J=
\mathbb E \left ( \frac{\partial \log \underline{h}^2_t(\nu_0)}{\partial \nu}
\frac{\partial \log \underline{h}^2_t(\nu_0)}{\partial \nu'} \right )
 .$$
This expression is again  in accordance with the one of Theorem 3.1 in \cite{HZ-APGARCH}.
\end{remi}
\section{Linear tests}\label{lineartest}
The asymptotic normality results from Theorem \ref{AN-connu} and Theorem \ref{AN-inconnu-c-est-le-label-de-bruno} are used to test linear constraints on the parameter.
We thus consider a null hypothesis of the form
\begin{equation}\label{Wald}
\textit{H}_0 : C\nu_0 = c,
\end{equation}
where $C$ is a known $s \times s_0$ matrix of rank $s$ and $c$ is a known  vector of size $s \times 1$. The Wald test is a standard parametric test for testing $\textit{H}_0$ and it is particularly appropriate in the context of financial series. Let $\hat{I}$ and $\hat{J}$ be weakly consistent estimators of
$I$ and $J$ involved in the asymptotic normality of the QMLE. For instance, $I$ and $J$ can be estimated by their empirical or observable counterparts
given by
\begin{equation*}
\hat{I} = \dfrac{1}{n} \sum\limits_{t=1}^n \dfrac{\partial \tilde{l}_t(\hat{\nu}_n)}{\partial \nu} \dfrac{\partial \tilde{l}_t(\hat{\nu}_n)}{\partial \nu'} \qquad \mbox{and}\qquad \hat{J} = \dfrac{1}{n}\sum\limits_{t=1}^n \dfrac{\partial^2 \tilde{l}_t(\hat{\nu}_n)}{\partial \nu\partial \nu'}.
\end{equation*}
Under the assumptions of Theorem~\ref{AN-inconnu-c-est-le-label-de-bruno}, and the assumption that the matrix $I$ is
invertible, the Wald test statistic is defined as follows
\begin{equation*}
W_n(s) = (C\hat{\nu}_n - c)' (C \hat{J}^{-1}\hat{I}\hat{J}^{-1}C')^{-1}(C\hat{\nu}_n - c).
\end{equation*}
We reject the null hypothesis for large values of $W_n(s)$. If $W_n(s) > \chi^2_s(\alpha)$ (here $\chi^2_s(\alpha)$ correspond to the quantile of the Khi$^2$ distribution of asymptotic level $\alpha$), we could conclude that the powers of the model \eqref{MAPGARCH} are equals. Thus we can consider a model with equal powers and so we are able to reduce the dimension of the parameter space.

We can also test the equality between the components of the power $\underline{\nu}$. For instance, in the bivariate case with the orders fixed at $p=q=1$, we take $C=(0_{14}, 1, 0, 0)$ and $c=\tau_2$ in \eqref{Wald} (remind that  $\nu_0  = (\underline{\omega}', \alpha^{+'},\alpha^{-'},\beta',{\tau}_1,{\tau}_2,\rho)'$).
Via this kind of test, we can also test the asymmetric property. If we reject the asymmetric assumption, we can consider the standard CCC-PGARCH model and reduce the dimension of the parameter space.\\

We strength the fact that our normality results can not be used in order to determine the orders $(p{,}q)$ of the model \eqref{MAPGARCH}. Indeed, this can be usually done when one tests the nullity of the coefficients of the matrix $B_p$. If they are all equal to zero, we can consider a smaller order $p$ (and similarly with the matrices $A_q^+$ and $A_q^-$ with the order $q$). Unfortunately, the assumption $\textbf{A6}$ excludes the case of vanishing coefficients so we can not apply our results in this situation. One refers to Section 8.2 of \cite{FZ-2010} and \cite{pedersen} for a discussion on this topic. 
\section{Numerical illustrations}
\label{NumIllust}
In this section, we make some simulations with $\eta_t \sim \mathcal{N}(0, I_2)$ and we compute the QML estimator to estimate the coefficients of the model \eqref{MAPGARCH}. The simulations are made with the open source statistical software R (see R Development Core Team, 2017) or (see http://cran.r-project.org/). We use the function {\tt nlminb()} to minimize the $\log-$quasi-likelihood.

\subsection{Estimation when the power is known}

 We fixed the orders of the model \eqref{MAPGARCH} at $q=1$ and $p=0$. 
We computed the estimator on $N=100$ independent simulated trajectories to assess the performance of the QMLE in finite sample. The trajectories are made in dimension 2 with a length $n=500$ and $n=5,000$. 

\begin{figure}[H]
\centering
\includegraphics[width=0.75\textwidth]{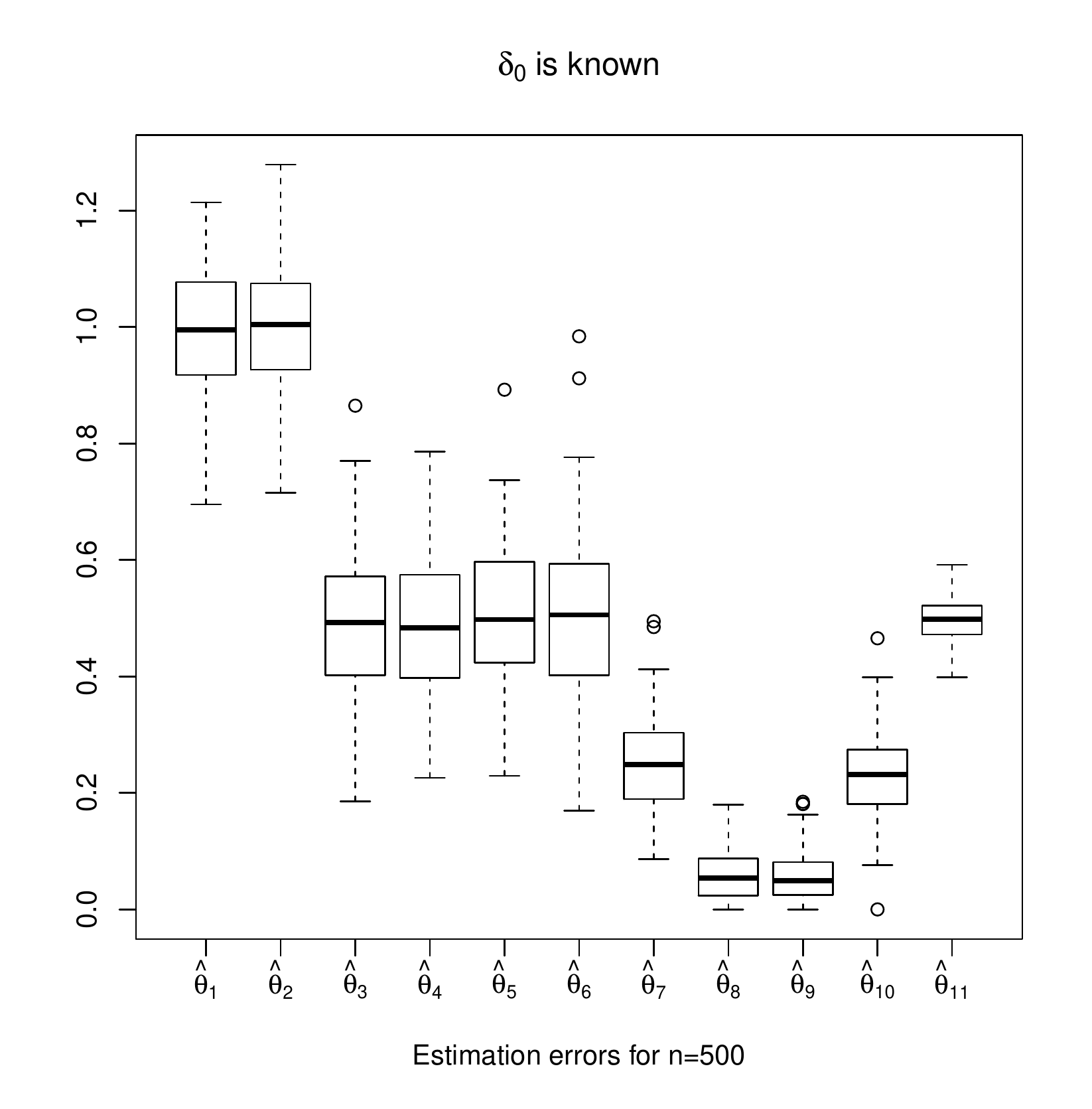}
\caption{\label{figure1} {\footnotesize Box-plot of the QMLE distribution for $n=500$}}
\end{figure}

\begin{figure}[H]
\centering
\includegraphics[width=0.75\textwidth]{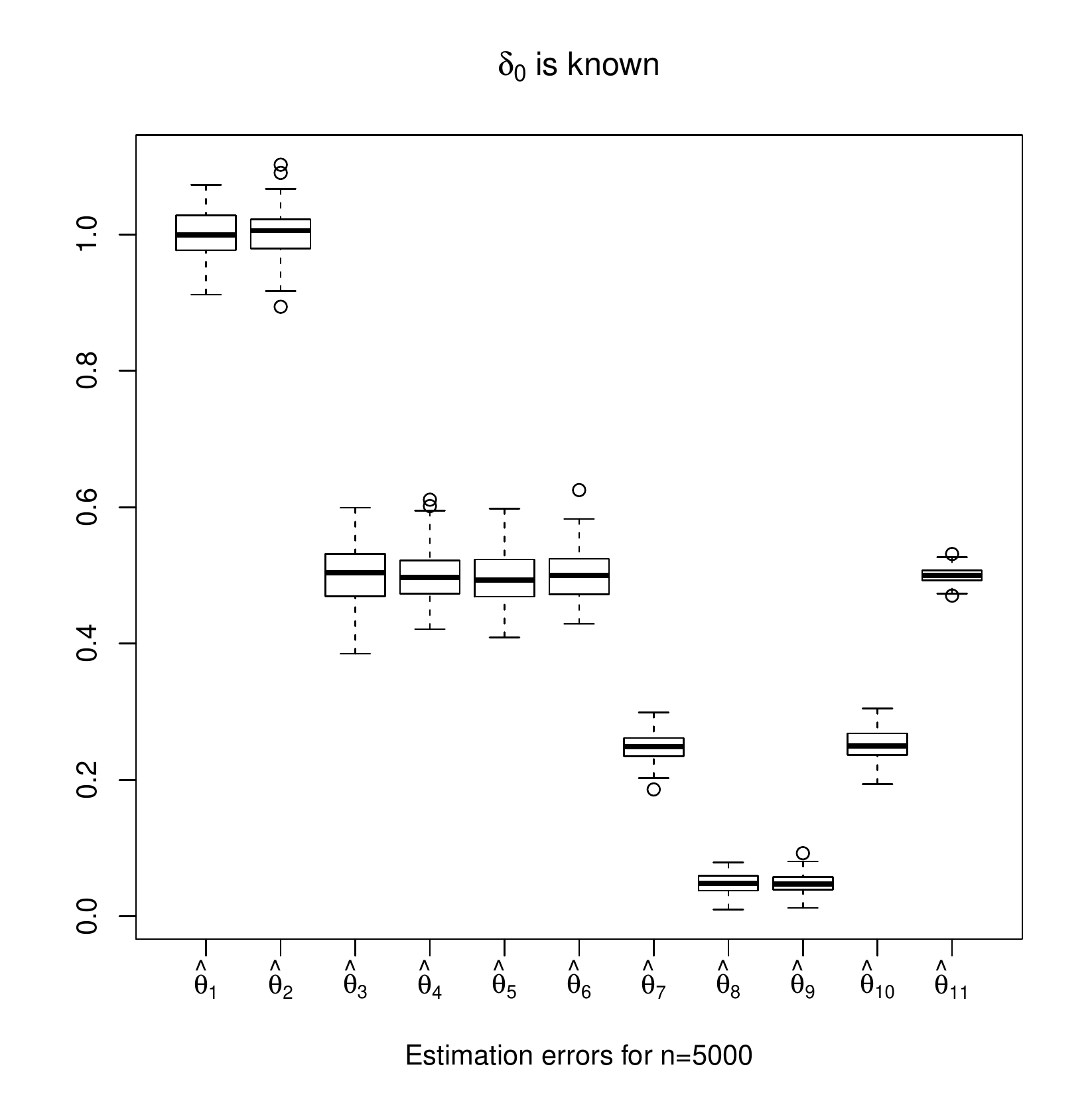}
\caption{\label{figure2} {\footnotesize Box-plot of the QMLE distribution for $n=5,000$}}
\end{figure}

The parameter used to simulate the trajectories are given in first row of the Table \ref{table2} 
and the space $\Theta$ associated is chosen to satisfied the assumptions of Theorem~\ref{AN-connu}.
As expected, Table \ref{table2} shows that the bias and the RMSE decrease when the size of the sample increases. Figures~\ref{figure1} and \ref{figure2} summarize via box-plot, the distribution of the QMLE for these simulations. Of course the precision around the estimated coefficients is better when the size of the sample increases (see Figures~\ref{figure1} and \ref{figure2}).
\begin{table}[H]
 \caption{\small{Sampling distribution of the QMLE of $\theta_0$ for the CCC-APGARCH(0,1) model \eqref{MAPGARCH} for $\underline{\delta}=(2,2)'$.}}
\begin{center}
\begin{tabular}{ll|c|cc|cc|c}
\hline \hline
& & $\underline{\omega}_0$ &  \multicolumn{2}{c|}{$A_0^+$} &  \multicolumn{2}{c|}{$A_0^-$} & $\rho_0$\\
Length & \multirow{2}{*}{True val.} & 1 & 0.25 & 0.05 & 0.5 & 0.5  & 0.5 \\
& & 1 & 0.05 & 0.25 & 0.5 & 0.5\\
\hline
\multirow{4}{*}{$n=500$} & \multirow{2}{*}{Bias} & -0.00498  & 0.00180 & 0.00789 &  -0.00525 & 0.00485 &  -0.00471 \\
& & -0.00090 & 0.00683 & -0.02083 & -0.01505 & 0.00122 & \\
\cline{2-8}
& \multirow{2}{*}{RMSE} & 0.10714& 0.08290 & 0.04520 & 0.12666 & 0.12379& 0.03774  \\
& & 0.12073 & 0.04308 & 0.07730 & 0.12738 & 0.15095 & \\
\hline
\multirow{4}{*}{$n=5,000$} & \multirow{2}{*}{Bias} & 0.00105  & -0.00063 & -0.00024  & 0.00346 & -0.00324  & 0.00010 \\
& & 0.00044 & -0.00175 & 0.00188 & 0.00175 & 0.00121& \\
\cline{2-8}
& \multirow{2}{*}{RMSE} & 0.03526  & 0.02129 & 0.01464 & 0.04219 & 0.03950 & 0.01143\\
& & 0.03745 & 0.01437 & 0.02220 & 0.03883 &0.03890& \\
\hline\hline \\
\end{tabular}
\end{center}
\label{table2}
\end{table}
%
\subsection{Estimation when the power is unknown}
\begin{figure}[H]
\centering
\includegraphics[width=0.75\textwidth]{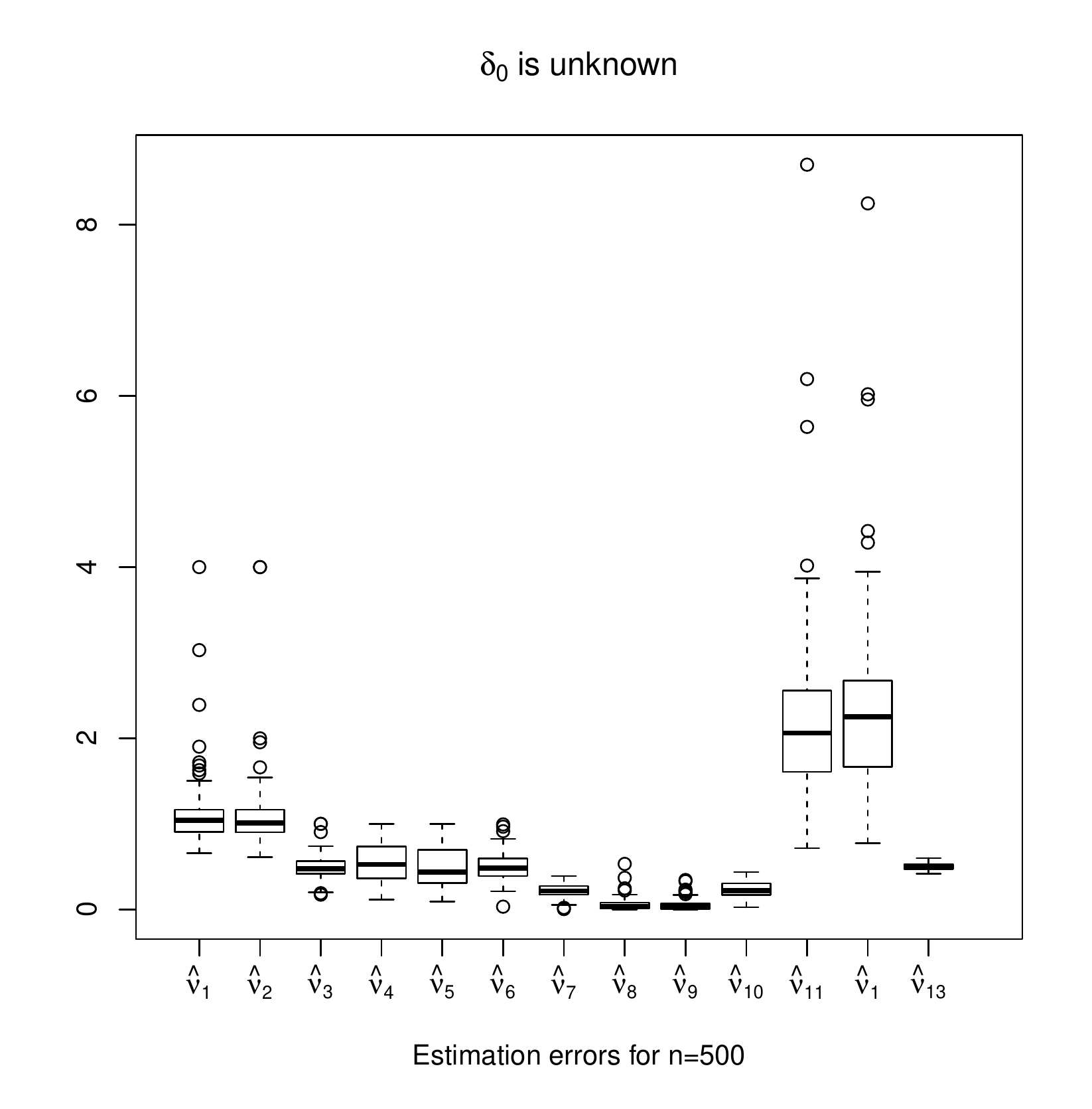}
\caption{\label{figure3} {\footnotesize Box-plot of the QMLE distribution for these simulations experiments, with $n=500$. }}
\end{figure}

\begin{figure}[H]
\centering
\includegraphics[width=0.75\textwidth]{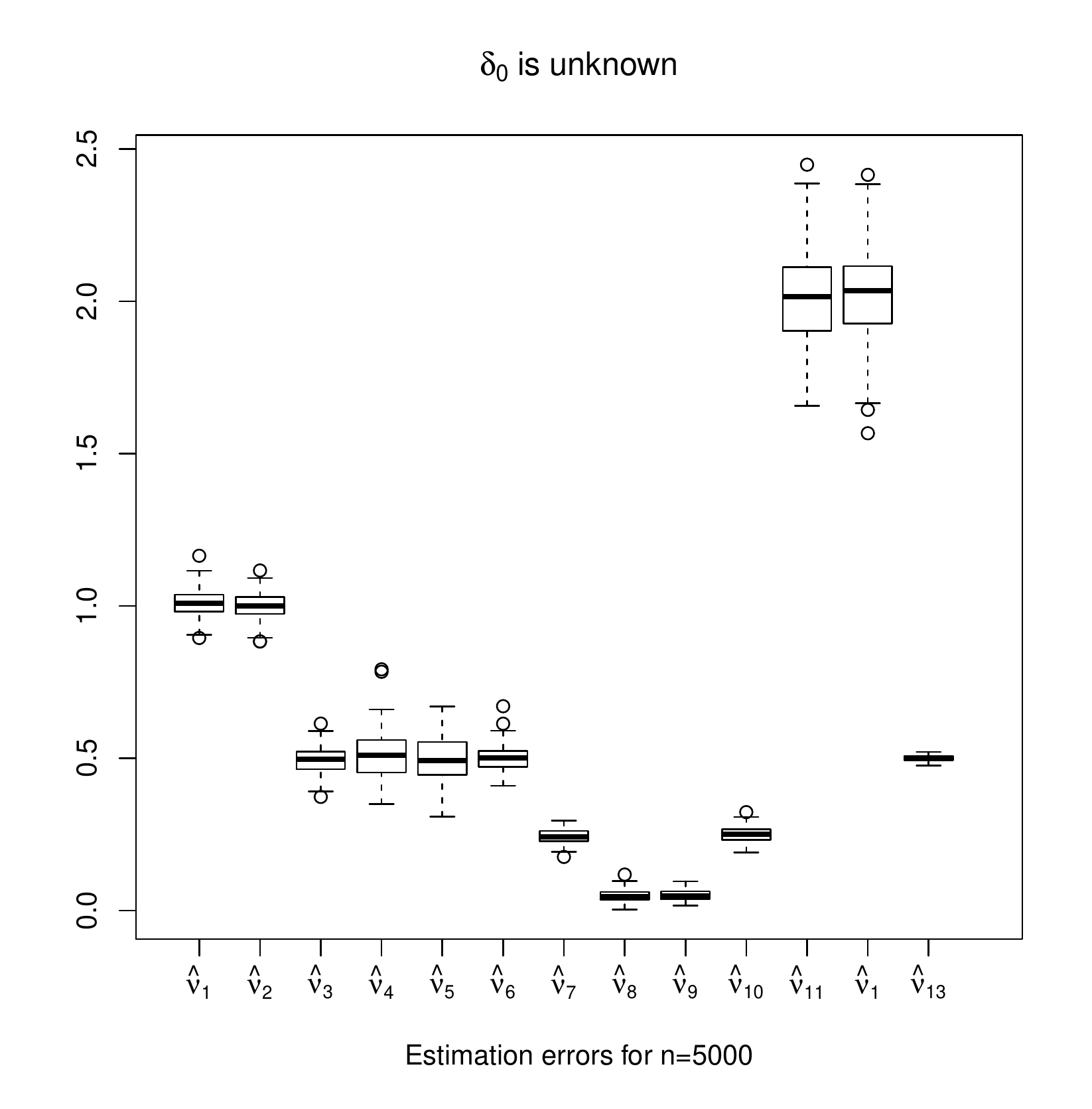}
\caption{\label{figure4} {\footnotesize Box-plot of the QMLE distribution for these simulations experiments, with $n=5,000$. }}
\end{figure}

In this section, we also fixed the orders $q=1$ and $p=0$ of Model \eqref{MAPGARCH} and we estimate the parameter $\nu_0$. We compute, as the previous section, the estimator on $N=100$ independent simulated trajectories in dimension 2 with sizes $n=500$ and $n=5,000$. Table \ref{table3} presents the results of QMLE of the parameters including the powers $\underline{\delta}=(\delta_1,\delta_2)'$. The results are satisfactory even for the parameter $\underline{\delta}$ which is difficult to identify in practice (see Section 4.1 of \cite{HZ-APGARCH}) and which is more difficult to estimate than the others parameters. The conclusions are similar as in the case where $\underline{\delta}$ is known. The dispersion around the parameter $\underline{\delta}$ is much more precise when the sample size is $n=5,000$. The bias and the RMSE decrease when the size of sample increases. Figure \ref{figure3} shows the distribution of the parameter $\hat\nu_n$ for the size $n=500$. We remark that the estimation of the parameter $\underline{\delta}$ has an important dispersion as regard to the other parameters. When the size of the sample is $n=5,000$ (see Figure \ref{figure4}), the estimation of the power is more accurate.

Table \ref{tabWald} displays the relative percentages of rejection of the Wald test  proposed in Section \ref{lineartest}  for testing the null hypothesis $\textit{H}_0: C\nu_0=(1,1)'$ in the case of bivariate CCC-APGARCH$(0,1)$, where $C$ is a $2\times 13$ matrix with 1 in positions $(1,11)$ and $(2,12)$ and 0 elsewhere and $\nu_0=(0.2,0.3,0.25,0.05,0.05,0.25,0.5,0.5,0.5,0.5,1,1,0.5)'$.  We simulated $N=100$ independent trajectories of size $n=1,000$ and $n=10,000$ of  Model \eqref{MAPGARCH} with $m=2$, $q=1$ and $p=0$. The nominal asymptotic level of the tests is $\alpha=5\%$. The values in bold correspond to the null hypothesis $\textit{H}_0$. We remark that the relative rejection frequencies of the Wald test are close to the nominal 5\% level under the null, and are close to 100\% under the alternative. We draw the conclusion that the proposed test well controls the error of first kind.
%

\begin{table}[H]
 \caption{\small{Sampling distribution of the QMLE of $\nu_0$ for the CCC-APGARCH$(0,1)$ model \eqref{MAPGARCH} with $\underline{\delta}$ unknown.}}
\begin{center}
\begin{tabular}{ll|c|cc|cc|c|c}
\hline \hline
& & $\underline{\omega}_0$ &  \multicolumn{2}{c|}{$A_0^+$} &  \multicolumn{2}{c|}{$A_0^-$} & $\underline{\delta}_0$ & $\rho_0$\\
Length & \multirow{2}{*}{True val.} & 1 & 0.25 & 0.05 & 0.5 & 0.5  & 2 & 0.5 \\
& & 1 & 0.05 & 0.25 & 0.5 & 0.5 & 2 & \\
\hline
\multirow{4}{*}{$n=500$} & \multirow{2}{*}{Bias} & 0.12600 & -0.03211 & 0.00370 & -0.00853 & -0.00155 & 0.27015 & -0.00112 \\
& & 0.11629 & 0.00937 & -0.02054 & 0.04720 & 0.00177 & 0.33354 & \\
\cline{2-9}
& \multirow{2}{*}{RMSE} & 0.45385 & 0.08449 & 0.06808 & 0.14385 & 0.27142 & 1.13662 & 0.03772 \\
& & 0.48839 & 0.07917 & 0.09275 & 0.26265 & 0.16417 & 1.14144 & \\
\hline
\multirow{4}{*}{$n=5,000$} & \multirow{2}{*}{Bias} & 0.00921 & -0.00490 & 0.00082 & -0.00526 & -0.00574 & 0.00000 & -0.00013 \\
& & 0.00072& -0.00052 & 0.00183  & 0.01204 & 0.00196 & 0.00000 & \\
\cline{2-9}
& \multirow{2}{*}{RMSE} & 0.05023 & 0.02434 & 0.01783 & 0.04278 & 0.07725 & 0.16761 & 0.01067  \\
& & 0.04825 & 0.01935 & 0.02661 & 0.08433 & 0.04402 & 0.16004 & \\
\hline\hline \\
\end{tabular}
\end{center}
\label{table3}
\end{table}

\begin{table}[H]
 \caption{\small{Empirical rejection frequency (in \%) of $W_n(s)$ for testing the null hypothesis $H_0: C\nu_0=(1,1)'$ in the case of a CCC-APGARCH$(0,1)$ model
(\ref{MAPGARCH}), where $C$ is a $2\times 13$ matrix with 1 in positions $(1,11)$ and $(2,12)$ and 0 elsewhere. 
The number of replications is $N=100$ and the nominal asymptotic level is $\alpha=5\%$. }}
\begin{center}
{\small
\begin{tabular}{l rrrr rrr}
\hline\hline
$(\tau_1,\tau_2)$& & \multicolumn{2}{c}{$(0.5,0.5)$} & \multicolumn{2}{c}{$(0.5,1)$} & \multicolumn{2}{c}{$(1,1)$}\\
Length $n$& &$1,000$&$10,000$ & $1,000$&$10,000$& $1,000$&$10,000$\\
\cline{3-8}  $W_n(2)$ &&100.0&100.0 &100.0 &100.0&{\bf 6.0}&{\bf 5.0}
\\
\hline
$(\tau_1,\tau_2)$& & \multicolumn{2}{c}{$(0.5,1.5)$} & \multicolumn{2}{c}{$(1,1.5)$} & \multicolumn{2}{c}{$(2,2)$}\\
Length $n$& &$1,000$&$10,000$& $1,000$&$10,000$& $1,000$&$10,000$\\
\cline{3-8} $W_n(2)$ &&100.0&100.0&83.0  &100.0 &80.0&100.0
\\
\hline\hline
\end{tabular}
}
\end{center}
\label{tabWald}
\end{table}
%

\subsection{Estimation on a real dataset}
In this section, we propose to estimate the bivariate series returns of the daily exchange rates between the Dollar (USD) and the Yen (JPY) against the Euro (EUR). The observations cover the period from January 4, 1999 to July 13, 2017 which correspond to 4746 observations. The data were obtained from the website of the National Bank of Belgium (https://www.nbb.be).
 We divide the full period in three subperiods with equal length ($n = 1581$): the first period runs from 1999-01-04 to 2005-03-02, the second from 2005-03-03 to 2011-05-05 and the last one from 2011-05-06 to 2017-07-13.
 
First, we consider the univariate case and an APGARCH$(1,1)$ model is estimated by QML on the full period and the three subperiods. Tables \ref{table3univdeltafixe} and \ref{table3univ} show the univariate analysis of the two series (when $\underline{\delta}_0$ is known and estimated, respectively) and give an idea of the behaviour of the modeling. We observe:
 \begin{enumerate}
 \item  a strong leverage effect (resp. no leverage effect) of the JPY with respect to its own past return for the full, first and second periods (resp. for the third periods) confirmed by the Wald test proposed in Section \ref{lineartest}.
%
 \item  almost no leverage effect is also confirmed by Wald test in the volatility of the USD for some periods 
 \item  a strong persistence of past volatility for the USD and JPY
 \item  different values of the power $\hat\tau$ for the USD in the three subperiods.
\end{enumerate}

To give an idea of the reliability of the parameter estimates (in the bivariate case) obtained from the full period (4746 observations), the period has been divided into 3 subperiods as in the univariate analysis (see Tables \ref{table3univdeltafixe} and \ref{table3univ}). 
Tables \ref{table4bivconnu} and \ref{table4bivinconnu} present the estimation of the model \eqref{MAPGARCH} for the bivariate series $(USD_t,JPY_t)'$ when $\underline{\delta}_0$ is known and estimated, respectively. The conclusion is the same as in Tables~\ref{table3univdeltafixe} and \ref{table3univ}. But we also remark a strong correlation between the two exchange rates. The power $\underline{\delta}_0$ is equivalent to the univariate case and differs for each subperiod. We also note that none is equal to the others. The persistence matrix $B_{01}$ has a diagonal form and the persistence of past volatility is strong as in the univariate case. The coefficients corresponding to the negative returns are generally bigger than those for the positive returns and some of them are equal to zero. Since we obtain null estimated coefficients, we think that the assumption $\textbf{A6}$ is not satisfied and thus our asymptotic normality theorems do not apply.
In future works we intent to study how the identification (orders $p$ and $q$ selection) procedure should be adapted in the CCC-APGARCH$(p,q)$ framework considered in the present paper by extending Section 8.2 of \cite{FZ-2010}.

\begin{table}[H]
 \caption{\small{Estimation of daily exchange rates of the Dollar and Yen against the Euro with the APGARCH$(1,1)$ model in dimension 1.}}
\begin{center}
\begin{tabular}{cc|cc|cc|cc|cc}
\hline \hline
$\underline{\delta}_0$&Parameters & \multicolumn{2}{c|}{Full Period} & \multicolumn{2}{c|}{First Period} & \multicolumn{2}{c|}{Second Period} & \multicolumn{2}{c}{Third Period} \\
& & USD & JPY & USD & JPY & USD & JPY & USD & JPY \\
\hline
&$\hat\omega$ & 0.00492 & 0.01107  & 0.01098 & 0.01245 & 0.00245 & 0.01192 & 0.00363 & 0.01115\\
\cline{2-10}
0.5&$\hat\alpha^+_1$ & 0.02573 & 0.05211  & 0.02132 & 0.05832 & 0.03822 & 0.03087 & 0.00984 & 0.05656\\
\cline{2-10}
&$\hat\alpha^-_1$ & 0.03512 & 0.06995  & 0.02919 & 0.07707 & 0.04447 & 0.06687 & 0.02780 & 0.06050\\
\cline{2-10}
&$\hat\beta_1$ & 0.96960 & 0.93866  & 0.96626 & 0.93106 & 0.96405 & 0.94787 & 0.98016 & 0.94095\\
\hline
&$\hat\omega$ & 0.00293 & 0.00803  & 0.00803 & 0.00993 & 0.00184 & 0.00826 & 0.00176 & 0.00826\\
\cline{2-10}
1&$\hat\alpha^+_1$ & 0.02628 & 0.05377  & 0.02205 & 0.06009 & 0.04306 & 0.02422 & 0.00127 & 0.06211\\
\cline{2-10}
&$\hat\alpha^-_1$ & 0.04059 & 0.08504  & 0.03498 & 0.09537 & 0.04933 & 0.08989 & 0.03315 & 0.06808\\
\cline{2-10}
&$\hat\beta_1$ & 0.96972 & 0.93620  & 0.96588 & 0.92635 & 0.96180 & 0.94543 & 0.98318 & 0.93977\\
\hline
&$\hat\omega$ & 0.00177 & 0.00597  & 0.00600 & 0.00815 & 0.00149 & 0.00564 & 0.00103 & 0.00637\\
\cline{2-10}
1.5&$\hat\alpha^+_1$ & 0.02447 & 0.05075  & 0.02048 & 0.05620 & 0.04050 & 0.01979 & 0.00205 & 0.05886\\
\cline{2-10}
&$\hat\alpha^-_1$ & 0.03803 & 0.08468  & 0.03274 & 0.09880 & 0.04617 & 0.09510 & 0.03102 & 0.06150\\
\cline{2-10}
&$\hat\beta_1$ & 0.97041 & 0.93463  & 0.96675 & 0.92254 & 0.96131 & 0.94293 & 0.98273 & 0.94044\\
\hline
&$\hat\omega$ & 0.00113 & 0.00466  & 0.00456 & 0.00693 & 0.00120 & 0.00396 & 0.00064 & 0.00509\\
\cline{2-10}
2&$\hat\alpha^+_1$ & 0.02100 & 0.04517  & 0.01672 & 0.04936 & 0.03408 & 0.01669 & 0.00406 & 0.05209\\
\cline{2-10}
&$\hat\alpha^-_1$ & 0.03156 & 0.07586  & 0.02646 & 0.08730 & 0.03970 & 0.08730  & 0.02402 & 0.04845\\
\cline{2-10}
&$\hat\beta_1$ & 0.97110 & 0.93242  & 0.96856 & 0.94046 & 0.96131 & 0.94046 & 0.98233 & 0.94082\\
\hline
&$\hat\omega$ & 0.00078 & 0.00378  & 0.00347 & 0.00597 & 0.00096 & 0.00285 & 0.00044 & 0.00425\\
\cline{2-10}
2.5&$\hat\alpha^+_1$ & 0.01652 & 0.03838  & 0.01164 & 0.04137 & 0.02680 & 0.01340 & 0.00386 & 0.04405\\
\cline{2-10}
&$\hat\alpha^-_1$ & 0.02438 & 0.06435  & 0.01945 & 0.08541 & 0.03246 & 0.07338 & 0.01624 & 0.03406\\
\cline{2-10}
&$\hat\beta_1$ & 0.97161 & 0.92895  & 0.97126 & 0.91149 & 0.96114 & 0.93840 & 0.98284 & 0.94060\\
\hline \hline
\end{tabular}
\end{center}
\label{table3univdeltafixe}
\end{table}

\begin{table}[H]
 \caption{\small{Estimation of daily exchange rates of the Dollar and Yen against the Euro with the APGARCH$(1,1)$ model in dimension 1.}}
\begin{center}
\begin{tabular}{c|cc|cc|cc|cc}
\hline \hline
Parameters & \multicolumn{2}{c|}{Full Period} & \multicolumn{2}{c|}{First Period} & \multicolumn{2}{c|}{Second Period} & \multicolumn{2}{c}{Third Period} \\
 & USD & JPY & USD & JPY & USD & JPY & USD & JPY \\
\hline
$\hat\omega$ & 0.00279 & 0.00740 & 0.01070 & 0.00896 & 0.00127 & 0.00810 & 0.00185 & 0.00777 \\
\hline
$\hat\alpha^+_1$ & 0.02618 & 0.05331 & 0.02150 & 0.05864 & 0.03599 & 0.02394 & 0.00161 & 0.06189 \\
\hline
$\hat\alpha^-_1$ & 0.04063 & 0.08616 & 0.02992 & 0.09854 & 0.04158 & 0.09056 & 0.03302 & 0.06754 \\
\hline
$\hat\beta_1$ & 0.96978 & 0.93580 & 0.96620 & 0.92450 & 0.96130 & 0.94531 & 0.98309 & 0.93990 \\
\hline
$\hat\tau$ & 1.04728 & 1.12923 & 0.53666 & 1.24883 & 1.86581 & 1.02499 & 0.96029 & 1.10951 \\\hline \hline
\end{tabular}
\end{center}
\label{table3univ}
\end{table}

\begin{table}[H]
 \caption{\small{Estimation of daily exchange rates of the (Dollar,Yen) against the Euro with the CCC-APGARCH$(1,1)$ model in dimension 2
 when $\underline{\delta}_0$ is fixed. }}
\begin{center}
\begin{tabular}{cccccc}
\hline \hline
$\underline{\delta}_0$&Parameters & Full Period & First Period & Second Period & Third Period \\
\hline
&$\underline{\hat\omega}$ & 0.00535 & 0.33950 & 0.01813 & 0.06369 \\
&& 0.14680 & 0.33068 & 0.08347 & 0.09406 \\
\cline{2-6}
&$\hat\alpha^+_1$ & 0.04265 & 0.01283 & 0.10252 & 0.14253 \\
&& 0.00000 & 0.00000 & 0.04520 & 0.00361 \\
&& 0.00484 & 0.11529 & 0.06218 & 0.02385 \\
&& 0.04232 & 0.03573 & 0.00000 & 0.07928 \\
\cline{2-6}
&$\hat\alpha^-_1$ & 0.02825 & 0.00000 & 0.05721 & 0.02225 \\
&& 0.00000 & 0.00000 & 0.00000 & 0.00000 \\
$(2,1)$&& 0.02203 & 0.03164 & 0.00000 & 0.06317 \\
&& 0.11239 & 0.17432 & 0.08536 & 0.04928 \\
\cline{2-6}
&$\hat\beta_1$ & 0.94718 & 0.37445 & 0.86761 & 0.73823 \\
&& 0.04727 & 0.00000 & 0.02670 & 0.00000 \\
&& 0.00000 & 0.00000 & 0.00000 & 0.00000 \\
&& 0.79038 & 0.63165 & 0.86991 & 0.87176 \\
\cline{2-6}
&$\hat\rho_{21}$ & 0.68316 & 0.74252 & 0.66079 & 0.69230 \\
\hline
&$\underline{\hat\omega}$ & 0.00600 & 0.06562 & 0.01928 & 0.04892 \\
&& 0.04538 & 0.12019 & 0.01354 & 0.02421 \\
\cline{2-6}
&$\hat\alpha^+_1$ & 0.03022 & 0.02592 & 0.07476 & 0.09216 \\
&& 0.00526 & 0.00000 & 0.06385 & 0.00262 \\
&& 0.00000 & 0.00000 & 0.02459 & 0.00734 \\
&& 0.06315 & 0.04730 & 0.00443 & 0.07986 \\
\cline{2-6}
&$\hat\alpha^-_1$ & 0.02421 & 0.02062 & 0.03516 & 0.03387 \\
&& 0.00000 & 0.0000 & 0.00000 & 0.00000 \\
$(2,2)$&& 0.00808 & 0.00986 & 0.00000 & 0.01241 \\
&& 0.14923 & 0.25040 & 0.13401 & 0.05892 \\
\cline{2-6}
&$\hat\beta_1$ & 0.95080 & 0.82370  & 0.87099 & 0.77251 \\
&& 0.00000 & 0.00000 & 0.00000 & 0.00000 \\
&& 0.00000 & 0.00000 & 0.00528 & 0.00000 \\
&& 0.81186 & 0.65577 & 0.88146 & 0.88318 \\
\cline{2-6}
&$\hat\rho_{21}$ & 0.55335 & 0.61313 & 0.52596 & 0.56145 \\
\hline\hline
\end{tabular}
\end{center}
\label{table4bivconnu}
\end{table}

\begin{table}[H]
 \caption{\small{Estimation of daily exchange rates of the (Dollar,Yen) against the Euro with the CCC-APGARCH$(1,1)$ model in dimension 2. }}
\begin{center}
\begin{tabular}{ccccc}
\hline \hline
Parameters & Full Period & First Period & Second Period & Third Period \\
\hline
$\hat{\underline{\omega}}$ & 0.00136 & 0.00117 & 0.02407 & 0.11245 \\
& 0.06124 & 0.15275 & 0.01453 & 0.07930 \\
\hline
$\hat\alpha^+_1$ & 0.03050 & 0.00719 & 0.06434 & 0.04457 \\
& 0.00000 & 0.00000 & 0.06129 & 0.01818 \\
& 0.00000 & 0.05530 & 0.01975 & 0.00000 \\
& 0.05368 & 0.04459 & 0.00292 & 0.05135 \\
\hline
$\hat\alpha^-_1$ & 0.02351 & 0.00000 & 0.03112 & 0.02393 \\
& 0.00000 & 0.00613 & 0.00000 & 0.00000 \\
& 0.01072 & 0.00994 & 0.00000 & 0.00000 \\
& 0.12207 & 0.15270 & 0.12752 & 0.01499 \\
\hline
$\hat\beta_1$ & 0.95326 & 0.41630 & 0.89072 & 0.90158 \\
& 0.03182 & 0.00000 & 0.00000 & 0.00000 \\
& 0.00000 & 0.00000 & 0.00382 & 0.00410 \\
& 0.80512 & 0.70082 & 0.88657 & 0.87580 \\
\hline
$\hat\rho_{21}$ & 0.55106 & 0.61076 & 0.52380 & 0.57556 \\
\hline
$\hat\tau_1$ & 2.01916 & 3.93035 & 1.95333 & 1.53699 \\
$\hat\tau_2$& 1.88965 & 1.72984  &  1.98917 & 1.77340 \\
\hline\hline
\end{tabular}
\end{center}
\label{table4bivinconnu}
\end{table}
\section{Conclusion}
\label{conclusion}
In this paper we propose a class of multivariate asymmetric GARCH models which includes numerous functional forms of MGARCH. We provide an explicit necessary and sufficient condition to the strict stationary of the proposed model. In addition the asymptotic properties of the QMLE are investigated in the two cases of $\underline{\delta}$ known and unknown). We remark  that moment conditions on the observed process are not needed. A Wald test is proposed to test $s$ linear constraints on the parameter. In Monte Carlo experiments we demonstrated that the QMLE of the CCC-APGARCH models provide some satisfactory results, at least for the models considered in our study.
%
%
%
\appendix
\section{Appendix : Proofs of the main results}\label{app}
The proofs of our results are quite technical. These are adaptations of the arguments used in \cite{FZ-MAPGARCH} when the power is known, \cite{HZ-APGARCH} and \cite{pan} when the power is unknown. We strength the fact that new techniques and arguments are applied in the proof of of the identifiability (see Subsection \ref{new-proof1} and in the proof of the invertibility of the Fisher information  matrix (see Subsection \ref{new-proof2}). 
\subsection{Preliminaries}
In the following technical proofs we will use the following notations. 
We will use the multiplicative norm defined as:
\[\Vert A \Vert := \sup\limits_{\Vert x \Vert \leq 1}\Vert Ax\Vert = \rho^{1/2}(A'A),\]
where $A$ is a matrix of size $d_1 \times d_2$ and $\Vert x \Vert$ is the Euclidian norm of vector $x \in \mathbb{R}^{d_2}$ and $\rho(\cdot)$ is the spectral radius. We recall that this norm satisfies
\[\Vert A \Vert^2 \leq \sum\limits_{i,j}a_{i,j}^2 = Tr(A'A) \leq d_2\Vert A\Vert^2, \quad \vert A'A\vert \leq \Vert A\Vert^{2d_2}.\]
Moreover we have the following relation
\[\vert Tr(AB)\vert \leq (d_1d_2)^{1/2}\Vert A \Vert \Vert B \Vert \]
as long as the matrix product is well defined (actually $B$ is a $d_2\times d_1$ matrix).

We recall some useful derivation rules for matrix valued functions. If we consider $f(A)$ a scalar function of a matrix $A$, where all the $a_{ij}$ are considered as a function of an one real variable $x$, we have
\begin{align}
\label{DfA}
\dfrac{\partial f(A)}{\partial x} &  = \sum\limits_{i,j}\dfrac{\partial f(A)}{\partial a_{ij}}\dfrac{\partial a_{ij}}{\partial x} = Tr\left(\dfrac{\partial f(A)}{\partial A'}\dfrac{\partial A}{\partial x}\right).
\end{align}
For a non singular matrix $A$, we have the following relations:
\begin{equation}
\dfrac{\partial c' A c}{\partial A'} = cc'
\end{equation}
\begin{equation}
\dfrac{\partial Tr(CA'BA')}{\partial A'} = C'AB' + B'AC'
\end{equation}
\begin{equation}
\dfrac{\partial \log\vert \det(A)\vert}{\partial A'} = A^{-1}
\end{equation}
\begin{equation}
\dfrac{\partial A^{-1}}{\partial x} = -A^{-1}\dfrac{\partial A}{\partial x}A^{-1}
\end{equation}
\begin{equation}\label{Dtrace}
\dfrac{\partial Tr(CA^{-1}B)}{\partial A'} = -A^{-1}BCA^{-1}
\end{equation}
\begin{equation}
\dfrac{\partial Tr(CAB)}{\partial A'} = BC
\end{equation}
%
%

\subsection{Proof of Theorem \ref{le-joli-label-de-yacouba}}\label{A2}
We prove the consistency of the QMLE result when $\underline{\delta}=\underline{\delta_0}$ known following the same lines as in \cite{FZ-MAPGARCH}.

We first rewrite Equation \eqref{theorique} in the matrix form as follows:
\begin{equation}\label{EcMat1}
\mathbb{H}_t(\theta) = \underline{c}_t(\theta) + \mathbb{B}(\theta)\cdot \mathbb{H}_{t-1}(\theta),
\end{equation}
with
\[ \mathbb{H}_t(\theta) = \begin{pmatrix} \underline{h}_t^{\underline{\delta}/2}(\theta) \\ \underline{h}_{t-1}^{\underline{\delta}/2}(\theta) \\ \vdots \\ \underline{h}_{t-p+1}^{\underline{\delta}/2}(\theta) \end{pmatrix},\quad \underline{c}_t(\theta) = \begin{pmatrix} \underline{\omega} + \sum\limits_{i=1}^q A_i^+(\underline{\varepsilon}_{t-i}^+)^{\underline{\delta}/2} + A_i^- (\underline{\varepsilon}_{t-1}^-)^{\underline{\delta}/2} \\ 0 \\ \vdots \\ 0 \end{pmatrix}\]
and
\[\mathbb{B}(\theta) = \begin{pmatrix} B_1 & B_2 & & \ldots & & B_p\\ I_m & & & & & 0 \\ 0 & I_m & & & & \vdots\\ \vdots & & & \ddots & & \vdots \\ 0 & & & \ldots & I_m & 0 \end{pmatrix}. \]
For simplicity, one writes $\mathbb{H}_t$ instead of $\mathbb{H}_t(\theta)$ when there is no possible confusion (and analogously one writes $ \underline{c}_t$ and $\mathbb{B}$).
We iterate the expression \eqref{EcMat1} and we obtain
\begin{equation}\label{Hgrandinfini}
\mathbb{H}_t = \underline{c}_t + \mathbb{B}\underline{c}_{t-1} + \mathbb{B}^2\underline{c}_{t-2} + \ldots + \mathbb{B}^{t-1}\underline{c}_1 + \mathbb{B}^t\mathbb{H}_0 = \sum\limits_{k=0}^\infty \mathbb{B}^k \underline{c}_{t-k}\ .
\end{equation}

The proof is decomposed in the four following points which will be treated in separate subsections.
\begin{enumerate}
	\item[$(i)$] Initial values do not influence quasi-likelihood: $\lim_{n\to\infty}\sup_{\theta \in \Theta} \vert \mathcal{L}_n(\theta) - \tilde{\mathcal{L}}_n(\theta) \vert  = 0$, almost surely.
	\item[$(ii)$] Identifiability: if there exists $t \in \mathbb{Z}$ such that $\underline{h}_t^{\underline{\delta}_0/2}(\theta) = \underline{h}_t^{\underline{\delta}_0/2}(\theta_0)$ almost surely and $R = R_0$, then $\theta = \theta_0$
	\item[$(iii)$] Minimization of the quasi-likelihood on the true value: $\mathbb{E}_{\theta_0} \vert l_t(\theta_0) \vert < \infty,$ and if $\theta \neq \theta_0$, $\mathbb{E}_{\theta_0}[l_t(\theta)] > \mathbb{E}_{\theta_0} [l_t(\theta_0)]$.
	\item[$(iv)$] For any $\theta \neq \theta_0$ there exists a neighborhood $V(\theta)$ such that
\begin{align}\label{iv}
& \underset{n\to \infty}{\lim\inf}\inf\limits_{\theta^\ast \in V(\theta)} \tilde{\mathcal{L}}_n(\theta^\ast) > \mathbb{E}_{\theta_0}l_1(\theta_0),\quad a.s.
\end{align}
\end{enumerate}
\subsubsection{Initial values do not influence quasi-likelihood}
\label{A21}
We define the vectors $\tilde{\mathbb{H}}_t$ by replacing the variables $\underline{h}_{t-i}^{\underline{\delta}/2}$, $i=0,\dots,p-1$,  in $\mathbb{H}_t$ by $\underline{\tilde h}_{t-i}^{\underline{\delta}/2}$ and we have
\begin{align}\label{ec2}
\tilde{\mathbb{H}}_t &= \underline{c}_t + \mathbb{B}\underline{c}_{t-1} +  \ldots + \mathbb{B}^{t-q-1}\underline{c}_{q+1} + \mathbb{B}^{t-q}\tilde{\underline{c}}_{q} + \ldots + \mathbb{B}^{t-1}\tilde{\underline{c}}_1 + \mathbb{B}^t\mathbb{H}_0\ ,
\end{align}
where the vector $\tilde{\underline{c}}_t$ is obtained by replacing $\underline{\varepsilon}_0,\ldots, \underline{\varepsilon}_{1-q}$ in $\underline{c}_t$ by the initial values.\\

From the Assumption $\textbf{A2}$ and Corollary \ref{cor1},  we have $\rho(\mathbb{B})<1$ and we deduce from the compactness of $\Theta$ that we have $\sup\limits_{\theta\in \Theta} \rho(\mathbb{B})<1$.
Using the two iterative equations \eqref{EcMat1} and \eqref{ec2}, we obtain almost surely that for any $t$:
\begin{equation}\label{MajH1}
\begin{aligned}
\sup\limits_{\theta\in\Theta} \left\Vert \mathbb{H}_t - \tilde{\mathbb{H}}_t\right\Vert &= \sup\limits_{\theta\in\Theta} \left\Vert \sum\limits_{k=1}^q \mathbb{B}^{t-k}(\underline{c}_k - \tilde{\underline{c}}_k) + \mathbb{B}^t(\mathbb{H}_0 - \tilde{\mathbb{H}}_0)\right\Vert 
\leq K\rho^t \ ,
\end{aligned}
\end{equation}
where $K$ is a random constant that depends on the past values of $\{\varepsilon_t, t\le 0\}$. We may write \eqref{MajH1} as 
\begin{equation}\label{Majh.delta}
\begin{aligned}
\sup\limits_{\theta\in\Theta} \left\Vert \underline{ h}_{t}^{\underline{\delta}/2}(\theta) - \underline{\tilde h}_{t}^{\underline{\delta}/2}(\theta)\right\Vert
\leq K\rho^t \ .
\end{aligned}
\end{equation}
Thus, for $i_1=1,\dots,m$, since $\min\left( { h}_{i_1,t}^{\delta_{i_1}/2}(\theta),{\tilde h}_{i_1,t}^{{\delta_{i_1}}/2}(\theta)\right)\geq\omega=\inf\limits_{1\leq i\leq m}\underline{\omega}(i)$, the mean-value theorem implies that
\begin{align}\label{MajH2}
\sup\limits_{\theta\in\Theta} \left| { h}_{i_1,t}(\theta)-{\tilde h}_{i_1,t}(\theta)\right| &\leq\frac{2}{\delta_{i_1}} \sup\limits_{\theta\in\Theta}\max\left( { h}_{i_1,t}^{1-\delta_{i_1}/2}(\theta),{\tilde h}_{i_1,t}^{1-{\delta_{i_1}/2}}(\theta)\right) \sup\limits_{\theta\in\Theta}
\left| { h}_{i_1,t}^{\delta_{i_1}/2}(\theta)-{\tilde h}_{i_1,t}^{{\delta_{i_1}}/2}(\theta)\right| \nonumber \\
&\leq\frac{2K}{\delta_{i_1}} \left(\sup\limits_{\theta\in\Theta}\frac{1}{\omega}\right)\sup\limits_{\theta\in\Theta}\max\left( { h}_{i_1,t}(\theta),{\tilde h}_{i_1,t}(\theta)\right) \rho^t
\leq K\rho^t \ ,
\end{align}
and similarly
\begin{align}\label{MajH2.1}
\sup\limits_{\theta\in\Theta} \left| { h}_{i_1,t}^{1/2}(\theta)-{\tilde h}_{i_1,t}^{1/2}(\theta)\right| &\leq\frac{1}{\delta_{i_1}} \sup\limits_{\theta\in\Theta}\max\left( { h}_{i_1,t}^{(1-\delta_{i_1})/2}(\theta),{\tilde h}_{i_1,t}^{(1-{\delta_{i_1}})/2}(\theta)\right) \sup\limits_{\theta\in\Theta}
\left| { h}_{i_1,t}^{\delta_{i_1}/2}(\theta)-{\tilde h}_{i_1,t}^{{\delta_{i_1}}/2}(\theta)\right|  \nonumber \\
&\leq\frac{K}{\delta_{i_1}} \left(\sup\limits_{\theta\in\Theta}\frac{1}{\omega}\right)\sup\limits_{\theta\in\Theta}\max\left( { h}_{i_1,t}^{1/2}(\theta),{\tilde h}_{i_1,t}^{1/2}(\theta)\right) \rho^t
\leq K\rho^t \ .
\end{align}
From \eqref{MajH2} we can deduce that, almost surely, we have
\begin{equation}\label{MajH-H1}
\begin{aligned}
\sup\limits_{\theta\in\Theta} \left\Vert H_t - \tilde{H}_t\right\Vert \leq K\rho^t,\qquad \forall t.
\end{aligned}
\end{equation}
Since $\Vert R^{-1} \Vert$ is the inverse of the eigenvalue of smaller module of $R$ and $\Vert \tilde{D}_t^{-1}\Vert = [\min_i(h_{ii,t}^{1/\delta_i})]^{-1}$, we have
\begin{align}\label{MajDR1}
\sup\limits_{\theta\in\Theta} \Vert \tilde{H}_t^{-1}\Vert&  \leq \sup\limits_{\theta\in\Theta}\Vert \tilde{D}^{-1}\Vert^2\Vert R^{-1}\Vert \leq \sup\limits_{\theta\in\Theta}\left [\min_i(\underline{\omega}(i)^{1/\delta_i})\right ]^{-1}\Vert R^{-1}\Vert \leq K,
\end{align}
by using the fact that $R$ is a positive-definite matrix (see assumption \textbf{A5}), the compactness of $\Theta$ and the strict positivity of the components of $\underline{\omega}$. Similarly, we have 
\begin{equation}\label{MajH1.1}
\sup\limits_{\theta\in\Theta} \Vert H_t^{-1}\Vert \leq K.
\end{equation}
One may writes
\begin{equation*}
\begin{aligned}
\sup\limits_{\theta\in \Theta} \vert \mathcal L_t(\theta) - \tilde{\mathcal L}_t(\theta)\vert &\leq \dfrac{1}{n} \sum\limits_{t=1}^n\sup\limits_{\theta\in\Theta}\left\vert \underline{\varepsilon}_t'(H_t^{-1} - \tilde{H}_{t}^{-1})\underline{\varepsilon}_t + \log\vert H_t \vert - \log \vert \tilde{H}_t\vert\right\vert \\
&\leq \dfrac{1}{n} \sum\limits_{t=1}^n\sup\limits_{\theta\in\Theta}\left\vert \underline{\varepsilon}_t'(H_t^{-1} - \tilde{H}_{t}^{-1})\underline{\varepsilon}_t \right\vert + \dfrac{1}{n} \sum\limits_{t=1}^n\sup\limits_{\theta\in\Theta}\left\vert \log\vert H_t \vert - \log \vert \tilde{H}_t\vert\right\vert \\
&\le S_1+S_2 \ .
\end{aligned}
\end{equation*}
We can rewrite the first  term $S_1$ in the right hand side of the above inequality as
\begin{equation*}
\begin{aligned}
S_1 &=  \dfrac{1}{n} \sum\limits_{t=1}^n\sup\limits_{\theta\in\Theta}\left\vert \underline{\varepsilon}_t' H_t^{-1}(H_t - \tilde{H}_{t})\tilde{H}_t^{-1}\underline{\varepsilon}_t \right\vert \\ %
 &=  \dfrac{1}{n} \sum\limits_{t=1}^n\sup\limits_{\theta\in\Theta}\left\vert Tr(H_t^{-1}(H_t - \tilde{H}_{t})\tilde{H}_t^{-1}\underline{\varepsilon}_t \underline{\varepsilon}_t')\right\vert \\
 &\leq K\dfrac{1}{n} \sum\limits_{t=1}^n\sup\limits_{\theta\in\Theta}\Vert H_t^{-1}\Vert\Vert H_t - \tilde{H}_{t}\Vert\Vert\tilde{H}_t^{-1}\Vert\Vert\underline{\varepsilon}_t \underline{\varepsilon}_t'\Vert \\
 &\leq K \dfrac{1}{n} \sum\limits_{t=1}^n \rho^t \Vert\underline{\varepsilon}_t \underline{\varepsilon}_t'\Vert
\end{aligned}
\end{equation*}
where we have used \eqref{MajH-H1}, \eqref{MajDR1} and \eqref{MajH1.1}. Using the Borel-Cantelli lemma and Corollary \ref{cor2}, we deduce that $\rho^t \Vert\underline{\varepsilon}_t\underline{\varepsilon}_t'\Vert = \rho^t \underline{\varepsilon}_t'\underline{\varepsilon}_t$ goes to zero almost surely. Consequently, the Ces\'aro lemma implies that $n^{-1} \sum_{t=1}^n \rho^t\Vert \underline{\varepsilon}_t\underline{\varepsilon}_t'\Vert \to 0$ when $n$ goes to infinity.\\
For the second term $S_2$ we use $\log(1+x)\leq x$ for $ x> -1$ and the inequality $|\det(A)|\le \rho(A)^m\le  \|A\|^m$ and we obtain
\begin{equation*}
\begin{aligned}
\log\vert H_t\vert - \log\vert \tilde{H}_t\vert &= \log\vert H_t \tilde{H}_t^{-1}\vert \\
& = \log\vert I_m + (H_t - \tilde{H}_t)\tilde{H}_t^{-1}\vert \\
& \leq m \log\Vert I_m + (H_t - \tilde{H}_t)\tilde{H}_t^{-1}\Vert \\
& \leq m \Vert H_t - \tilde{H}_t\Vert \Vert \tilde{H}_t^{-1}\Vert \ ,
\end{aligned}
\end{equation*}
and, by symmetry, we have
\begin{equation*}
\begin{aligned}
\log\vert \tilde{H}_t\vert - \log\vert H_t\vert & \leq m \Vert H_t - \tilde{H}_t\Vert \Vert H_t^{-1}\Vert.
\end{aligned}
\end{equation*}
Using \eqref{MajH-H1}, \eqref{MajDR1} and \eqref{MajH1.1}, we have
$$
\left\vert \log\vert H_t \vert - \log \vert \tilde{H}_t\vert\right\vert \le
m \Vert H_t - \tilde{H}_t\Vert \left ( \Vert \tilde{H}_t^{-1}\Vert + \Vert {H}_t^{-1}\Vert\right )\leq K\rho^t . $$
Using the same arguments as for $S_1$, we conclude that $S_2$ goes to $0$.
We have shown that $\sup_{\theta\in\Theta} \vert \mathcal L_t(\theta) - \tilde{\mathcal L}_t(\theta)\vert  \underset{n\to \infty}{\longrightarrow} 0$ almost surely and thus $(i)$ is proved.
\subsubsection{Identifiability}\label{A22}
We suppose that for some $\theta \neq \theta_0$ we have
\begin{equation*}
\underline{h}_t(\theta) = \underline{h}_t(\theta_0),\quad \mathbb{P}_{\theta_0} -p.s., \quad  \mbox{and}\quad R(\theta) = R(\theta_0).
\end{equation*}
So we have $\rho = \rho_0$. Remind that by \eqref{IC} we have
\begin{equation*}
\underline{h}_{t}^{\underline{\delta}/2}(\theta_0) = \mathcal{B}^{-1}_0(1)\underline{\omega} + \mathcal{B}^{-1}_0(L)\mathcal{A}^+_0(L)(\underline{\varepsilon}_t^+)^{\underline{\delta}/2} + \mathcal{B}^{-1}_0(L)\mathcal{A}^-_0(L)(\underline{\varepsilon}_t^-)^{\underline{\delta}/2}.
\end{equation*}
We have a similar expression with the parameter $\theta$:
\begin{equation*}
\begin{aligned}
\underline{h}_t^{\underline{\delta}/2}(\theta) &= \mathcal{B}^{-1}(1)\underline{\omega} + \mathcal{B}^{-1}(L)\mathcal{A}^+(L)(\underline{\varepsilon}_t^+)^{\underline{\delta}/2} + \mathcal{B}^{-1}(L)\mathcal{A}^-(L)(\underline{\varepsilon}_t^-)^{\underline{\delta}/2}
\end{aligned}
\end{equation*}
Consequently we have
\begin{align*}
0 &= \mathcal{B}^{-1}(1)\underline{\omega} - \mathcal{B}^{-1}_0(1)\underline{\omega}_0  + \mathcal{B}^{-1}(L)\mathcal{A}^+(L)(\underline{\varepsilon}_t^+)^{\underline{\delta}/2} - \mathcal{B}^{-1}_0(L)\mathcal{A}^+_0(L)(\underline{\varepsilon}_t^+)^{\underline{\delta}/2}\\
& \hspace{2cm} +\mathcal{B}^{-1}(L)\mathcal{A}^-(L)(\underline{\varepsilon}_t^-)^{\underline{\delta}/2} - \mathcal{B}^{-1}_0(L)\mathcal{A}^-_0(L)(\underline{\varepsilon}_t^-)^{\underline{\delta}/2}
\end{align*}
and
\begin{align*}
\mathcal{B}^{-1}(1)\underline{\omega} - \mathcal{B}^{-1}_0(1)\underline{\omega}_0 &= \left[\mathcal{B}^{-1}(L)\mathcal{A}^+(L)- \mathcal{B}^{-1}_0(L)\mathcal{A}^+_0(L)\right](\underline{\varepsilon}_t^+)^{\underline{\delta}/2}+ \left[\mathcal{B}^{-1}(L)\mathcal{A}^-(L) - \mathcal{B}^{-1}_0(L)\mathcal{A}^-_0(L)\right](\underline{\varepsilon}_t^-)^{\underline{\delta}/2} \\
& = \mathcal{P}^+(L)(\underline{\varepsilon}_t^+)^{\underline{\delta}/2} + \mathcal{P}^-(L)(\underline{\varepsilon}_t^-)^{\underline{\delta}/2},
\end{align*}
where
\begin{equation}\label{ppm}
\mathcal{P}^\pm(L) = \mathcal{B}^{-1}_0(L)\mathcal{A}^\pm_0(L) - \mathcal{B}^{-1}(L)\mathcal{A}^\pm(L) = \sum\limits_{i=0}^\infty P_i^\pm L^i .
\end{equation}
We remark that $P_0^\pm = \mathcal{P}^\pm(0) = 0$ by the identifiability conditions. Using \eqref{eta-tilde} and \eqref{upsilon} we can write
\begin{equation*}
\begin{aligned}
\mathcal{P}^+(L)(\underline{\varepsilon}_t^+)^{\underline{\delta}/2} + \mathcal{P}^-(L)(\underline{\varepsilon}_t^-)^{\underline{\delta}/2} &= P_1^+(\underline{\varepsilon}_{t-1}^+)^{\underline{\delta}/2} + P_1^-(\underline{\varepsilon}_{t-1}^-)^{\underline{\delta}/2}+Z_{t-2}\\
&= P_1^+ \Upsilon_{t-1}^{+,(\delta)} \underline{h}_{t-1}^{\underline{\delta}/2} + P_1^- \Upsilon_{t-1}^{-,(\delta)} \underline{h}_{t-1}^{\underline{\delta}/2} + Z_{t-2}\quad a.s.,
\end{aligned}
\end{equation*}
where $Z_{t-2}$ is measurable with respect to the $\sigma$-field $\mathcal{F}_{t-2}$ generated by $\{\eta_{t-2},\eta_{t-3},\ldots\}$.
Hence we have
$$
 P_1^+ \Upsilon_{t-1}^{+,(\delta)} \underline{h}_{t-1}^{\underline{\delta}/2} + P_1^- \Upsilon_{t-1}^{-,(\delta)} \underline{h}_{t-1}^{\underline{\delta}/2} =  \mathcal{B}^{-1}(1)\underline{\omega} - \mathcal{B}^{-1}_0(1)\underline{\omega}_0-Z_{t-2}= \tilde Z_{t-2}$$
 where $\tilde Z_{t-2}$ is another $\mathcal{F}_{t-2}$-measurable random matrix.
Since $P_1^+\Upsilon_{t-1}^{+,(\delta)} + P_1^-\Upsilon_{t-1}^{-,(\delta)}$ is independent from $(Z_{t-2},\underline{h}_{t-1})$ and since $\underline{h}_{t-1}>0$, $P_1^+\Upsilon_{t-1}^{+,(\delta)} + P_1^-\Upsilon_{t-1}^{-,(\delta)} = C$ for some constant matrix $C$. Since the matrices-$\Upsilon_{t-1}^{\pm,(\delta)}$ are diagonal (see \eqref{upsilon}), the element $(i,j)$ of the matrix $C$ satisfies
$$ C(i,j) = P_1^{+}(i,j)(\tilde{\eta}_{j,t-1}^+)^{\delta_j} + P_1^{-}(i,j)(-\tilde\eta_{j,t-1}^-)^{\delta_j} . $$
If $P_1^+(i,j)P_1^-(i,j) \neq 0$, then $\tilde{\eta}_{j,t-1}$ takes at most two different values, which is in contradiction with \textbf{A3}.
If $P_1^+(i,j) \neq 0$ and $P_1^-(i,j) = 0$, then $P_1^+(i,j)(\tilde{\eta}_{j,t-1}^+)^{\delta_j} = C(i,j)$ which entails $C(i,j) = 0$, since $\mathbb{P}(\tilde{\eta}_{j,t-1}> 0) \neq 0$, and then $\tilde{\eta}_{j,t-1} = 0$, a.s., which is also in contradiction with \textbf{A3}. We thus have $P_1^+ = P_1^- = 0$.
We argue similarly for $P_i^{\pm}$ for $i\ge 2$ and by \eqref{ppm} we obtain that $\mathcal{P}^+(L) = \mathcal{P}^-(L) = 0$.
Therefore, in view of \eqref{ppm}, we may apply Proposition \ref{identif} because we assumed \textbf{A4} (or \textbf{A4'}). Thus we have  $\theta = \theta_0$.
We have thus established $(ii)$.

The proof of $(iii)$ and $(iv)$ is strictly identical to the one given in \cite{FZ-MAPGARCH}. Therefore it is  omitted and the proof of Theorem \ref{le-joli-label-de-yacouba} is
complete.
\zak


\subsection{Proof of Theorem \ref{AN-connu}}\label{anconnu}
Here, we prove the asymptotic normality result when $\underline{\delta}$ is known. The prove is based on the standard Taylor expansion. We have
\[\sqrt{n}(\hat{\theta}_n - \theta_0) = - \left [\dfrac1n \sum\limits_{t=1}^n \dfrac{\partial^2 \tilde{l}_t(\theta_{ij}^\ast)}{\partial \theta_i \partial\theta_j}\right ]^{-1}\left(\dfrac{1}{\sqrt{n}}\sum\limits_{t=1}^n \dfrac{\partial \tilde{l}_t(\theta_0)}{\partial\theta}\right),\]
where the parameter $\theta_{ij}^\ast$ is between $\hat{\theta}_n$ and $\theta_0$.
To establish the asymptotic normality result when the power is known, we will decomposed the proof in six intermediate points as in \cite{FZ-MAPGARCH}.
\begin{enumerate}
	\item[$(a)$] First derivative of the  quasi log-likelihood.
	\item[$(b)$] Existence of moments at any order of the score.
	\item[$(c)$] Asymptotic normality of the score vector:
	\begin{equation}\label{NA1}
\dfrac{1}{\sqrt{n}}\sum\limits_{t=1}^n \dfrac{\partial {l}_t(\theta_0)}{\partial\theta} \overset{\mathcal{L}}{ \underset{n\to \infty}{\longrightarrow}} \mathcal{N}(0,I).
\end{equation}
	\item[$(d)$] Convergence to $J$:
	\begin{equation}\label{ConvJ1}
\dfrac1n \sum\limits_{t=1}^n \dfrac{\partial^2 {l}_t(\theta_{ij}^\ast)}{\partial \theta_i \partial\theta_j}  \underset{n\to \infty}{\longrightarrow} J(i,j) \mbox{ in probability.}
\end{equation}
	\item[$(e)$] Invertibility of the matrix $J$.
	\item[$(f)$] Asymptotic irrelevance of the initial values.
\end{enumerate}
We shall need the following notations.
\begin{enumerate}[$\qquad \ast$]
	\item $s_0 = m + (p+2q)m^2 + m(m-1)/2$,
	\item $s_1 = m + (p+2q)m^2$,
	\item $s_2 = m + 2qm^2$,
	\item $s_3 = m + qm^2$.
\end{enumerate}
These notations will help us to point out the different entries of our parameter $\theta$.
\subsubsection{First derivative of log-likelihood}\label{A31}
The aim of this subsection is to establish the expressions of the first order derivatives of the quasi log-likelihood. We shall use the following notations:  $D_{0t} = D_t(\theta_0), R_0 = R(\theta_0)$,
\[D_{0t}^{(i)} = \dfrac{\partial D_t}{\partial \theta_i}(\theta_0), \quad R_0^{(i)} = \dfrac{\partial R}{\partial\theta_i}(\theta_0),\]
	\[D_{0t}^{(i,j)} = \dfrac{\partial^2D_t}{\partial\theta_i\partial\theta_j}(\theta_0),\quad R_0^{(i,j)} = \dfrac{\partial^2 R}{\partial\theta_i\partial\theta_j}(\theta_0),\]
and $\underline{\varepsilon}_t = D_{0t}\tilde{\eta}_t$, where $\tilde{\eta}_t(\theta) = R^{1/2}\eta_t(\theta)$ with $\tilde{\eta}_t=\tilde{\eta}_t(\theta_0) = R_0^{1/2}\eta_t$.

We recall the expression:
\begin{equation*}
\begin{aligned}
l_t(\theta) &= \underline{\varepsilon}_t' H_t^{-1} \underline{\varepsilon}_t + \log(\mathrm{det}(H_t)) \\
&= \underline{\varepsilon}_t'(D_tRD_t)^{-1}\underline{\varepsilon}_t + \log(\mathrm{det}(D_tRD_t))\\
	&= \underline{\varepsilon}_t' D_t^{-1} R^{-1}D_t^{-1} \underline{\varepsilon}_t + \log(\mathrm{det}(D_t) \mathrm{det}(R) \mathrm{det}(D_t))\\
	&= \underline{\varepsilon}_t' D_t^{-1} R^{-1}D_t^{-1} \underline{\varepsilon}_t + 2\log(\mathrm{det}(D_t)) + \log(\mathrm{det}(R)).
\end{aligned}
\end{equation*}
\noindent$\bullet$ We differentiate with respect to $\theta_i$ for $i = 1, \ldots, s_1$ (that is with respect to $(\underline{\omega}',  {\alpha_{1}^{+}} ', \ldots, {\alpha_q^+}', {\alpha_1^-}', \ldots, {\alpha_q^-}', \beta'_{1}, \ldots, \beta'_{p})'$). We have
\begin{align}\label{DL1}
\dfrac{\partial l_t(\theta)}{\partial \theta_i} & = -Tr\left((\underline{\varepsilon}_t\underline{\varepsilon}_t'D_t^{-1}R^{-1} + R^{-1}D_t^{-1}\underline{\varepsilon}_t\underline{\varepsilon}_t') D_t^{-1}\dfrac{\partial D_t}{\partial \theta_i}D_t^{-1}\right) + 2Tr\left(D_t^{-1} \dfrac{\partial D_t}{\partial\theta_i}\right) \\
\dfrac{\partial l_t(\theta_0)}{\partial \theta_i}& = Tr\left[ \left(I_m - R_0^{-1}\tilde{\eta}_t\tilde{\eta}_t'\right)D_{0t}^{(i)} D_{0t}^{-1} + \left(I_m - \tilde{\eta}_t\tilde{\eta}_t'R_0^{-1}\right)D_{0t}^{-1}D_{0t}^{(i)}\right]. \label{DL1bis}
\end{align}
\noindent$\bullet$  We differentiate with respect to $\theta_i$  for $i = s_1+1, \ldots, s_0$ (that is with respect to $\rho'$). We have
\begin{align}
\label{DL2}
\dfrac{\partial l_t(\theta)}{\partial \theta_i} &= -Tr\left( R^{-1} D_t^{-1} \underline{\varepsilon}_t\underline{\varepsilon}_t' D_t^{-1}R^{-1}\dfrac{\partial R}{\partial\theta_i} \right) + Tr\left(R^{-1} \dfrac{\partial R}{\partial\theta_i}\right) \\
\label{DL2bis}
\dfrac{\partial l_t(\theta_0)}{\partial \theta_i}& = Tr\left[\left(I_m - R^{-1}\tilde{\eta}_t\tilde{\eta}_t'\right)R^{-1}R^{(i)}\right].
\end{align}
\subsubsection{Existence of moments of any order of the score}\label{A32}
\begin{enumerate}[(i)]
\item For $i=1, \ldots, s_1$, in view of \eqref{DL1bis} we have 
\begin{equation*}
\begin{aligned}
 \left |\dfrac{\partial l_t(\theta_0)}{\partial \theta_i}\right | & = \left\vert Tr\left[\left(I_m - R_0^{-1}\tilde{\eta}_t\tilde{\eta}_t'\right)D_{0t}^{(i)} D_{0t}^{-1} + \left(I_m - \tilde{\eta}_t\tilde{\eta}_t'R_0^{-1}\right)D_{0t}^{-1}D_{0t}^{(i)}\right] \right\vert \\
 &\leq 2m \underbrace{\left\Vert \left(I_m - R_0^{-1}\tilde{\eta}_t\tilde{\eta}_t'\right) \right\Vert}_{\text{constant}} \left\Vert D_{0t}^{(i)} D_{0t}^{-1} \right\Vert  \leq K \left\Vert D_{0t}^{(i)} D_{0t}^{-1} \right\Vert .
\end{aligned}
\end{equation*}
\item For $i=s_1+1,\ldots, s_0$, in view of \eqref{DL2bis}  we have
\begin{equation*}
\begin{aligned}
 \left |\dfrac{\partial l_t(\theta_0)}{\partial \theta_i}\right | & = \left\vert Tr\left[\left(I_m - R_0^{-1}\tilde{\eta}_t\tilde{\eta}_t'\right)R_{0}^{-1} R_{0}^{(i)} \right] \right\vert \\
  &\leq m \underbrace{\left\Vert \left(I_m - R_0^{-1}\tilde{\eta}_t\tilde{\eta}_t'\right) \right\Vert}_{\text{constant}} \underbrace{\left\Vert R_{0}^{-1} R_{0}^{(i)} \right\Vert}_{\text{constant}} \leq K.
\end{aligned}
\end{equation*}
\item For $i,j = 1,\ldots, s_1$,  in view of \eqref{DL1bis}  we have
\begin{equation*}
\begin{aligned}
\mathbb{E}\left\vert\dfrac{\partial l_t(\theta_0)}{\partial\theta_i} \dfrac{\partial l_t(\theta_0)}{\partial\theta_j}\right\vert  &= \mathbb{E}\left\vert Tr\left[\left(I_m - R_0^{-1}\tilde{\eta}_t\tilde{\eta}_t'\right)D_{0t}^{(i)} D_{0t}^{-1} + \left(I_m - \tilde{\eta}_t\tilde{\eta}_t'R_0^{-1}\right)D_{0t}^{-1}D_{0t}^{(i)}\right] \right. \\
	&\qquad \left.\times Tr\left[\left(I_m - R_0^{-1}\tilde{\eta}_t\tilde{\eta}_t'\right)D_{0t}^{(j)} D_{0t}^{-1} + \left(I_m - \tilde{\eta}_t\tilde{\eta}_t'R_0^{-1}\right)D_{0t}^{-1}D_{0t}^{(j)}\right]\right\vert \\
	&\leq \mathbb{E}\left[ K_1 \left\Vert D_{0t}^{(i)} D_{0t}^{-1} \right\Vert \times K_2 \left\Vert D_{0t}^{(j)} D_{0t}^{-1} \right\Vert \right]\\
	&\leq K\left(\mathbb{E}\Vert D_{0t}^{(i)}D_{0t}^{-1} \Vert^2\mathbb{E}\Vert D_{0t}^{(j)}D_{0t}^{-1} \Vert^2\right)^{1/2}.
\end{aligned}
\end{equation*}
 by using the Cauchy-Schwarz inequality.
\item In view of \eqref{DL1bis} and \eqref{DL2bis}, we also have for $i= 1, \ldots, s_1$ and $j = s_1+1,\ldots, s_0$,
\begin{equation*}
\begin{aligned}
\mathbb{E} \left\vert\dfrac{\partial l_t(\theta_0)}{\partial\theta_i} \dfrac{\partial l_t(\theta_0)}{\partial\theta_j}\right\vert  &= \mathbb{E}\left\vert Tr\left[\left(I_m - R_0^{-1}\tilde{\eta}_t\tilde{\eta}_t'\right)D_{0t}^{(i)} D_{0t}^{-1} + \left(I_m - \tilde{\eta}_t\tilde{\eta}_t'R_0^{-1}\right)D_{0t}^{-1}D_{0t}^{(i)}\right] \right. \\
	& \qquad\times \left.Tr\left[ \left( I_m - R_0^{-1} \tilde{\eta}_t \tilde{\eta}_t' \right) R_0^{-1} R_0^{(j)} \right] \right \vert \\
	&\leq \mathbb{E}\left[ K_1\left\Vert D_{0t}^{(i)} D_{0t}^{-1} \right\Vert  \times K_2 \right] \leq K \mathbb{E}\left\Vert D_{0t}^{(i)} D_{0t}^{-1} \right\Vert.
\end{aligned}
\end{equation*}
\item For $i,j = s_1+1, \ldots, s_0$, , in view of \eqref{DL2bis}  we have,
\begin{equation*}
\begin{aligned}
\mathbb{E}\left\vert\dfrac{\partial l_t(\theta_0)}{\partial\theta_i} \dfrac{\partial l_t(\theta_0)}{\partial\theta_j}\right\vert  &= \mathbb{E}\left\vert Tr\left[ \left( I_m - R_0^{-1} \tilde{\eta}_t \tilde{\eta}_t' \right) R_0^{-1} R_0^{(i)} \right] \right. \\
	&\qquad \left.\times Tr\left[ \left( I_m - R_0^{-1} \tilde{\eta}_t \tilde{\eta}_t' \right) R_0^{-1} R_0^{(j)} \right] \right \vert  \leq \mathbb{E}[K_1 \times K_2]\leq K.
\end{aligned}
\end{equation*}
\end{enumerate}
To have the finiteness of the moments of the first derivative of the log-likelihood, it remains to treat the cases   (i), (iii) and (iv) above. Thus, we  have to control the term $\Vert D_{0t}^{(i)}D_{0t}^{-1}\Vert$.
Since
\[ D_{0t} = \mathrm{Diag}(\underline{h}_{t}^{1/2}(\theta_0)) = \mathrm{Diag}\left ( (\underline{h}_{t}^{\underline{\delta}/2}(\theta_0))^{1/\underline{\delta}}\right ) , \]
we can work component by component. We have for $i_1 = 1,\ldots, m$ and $i=1,\ldots,s_1$
\begin{align}\label{DD}
\dfrac{\partial D_{0t}(i_1,i_1)}{\partial \theta_{i}} =\dfrac{\partial \left(h_{i_1,t}^{\delta_{i_1}/2}\right)^{1/\delta_{i_1}}}{\partial \theta_{i}} = \dfrac{1}{\delta_{i_1}} h_{i_1,t}^{1/2} \times \dfrac{1}{h_{i_1,t}^{\delta_{i_1}/2}} \dfrac{\partial h_{i_1,t}^{\delta_{i_1}/2}}{\partial\theta_i}(\theta_0).
\end{align}
Control the term $\Vert D_{0t}^{(i)}D_{0t}^{-1}\Vert$ is equivalent to control $1/ h_{i_1,t}^{\delta_{i_1}/2}\partial h_{i_1,t}^{\delta_{i_1}/2} / \partial \theta_i$ in $\theta_0$. So it is  sufficient to prove that for any $r_0\ge 1$
\begin{equation} \label{score.moment.r0}
\mathbb{E}\left\vert \dfrac{1}{h_{i_1,t}^{\delta_{i_1}/2}} \dfrac{\partial h_{i_1,t}^{\delta_{i_1}/2}}{\partial\theta_i}(\theta_0)\right\vert^{r_0} < \infty\ .
\end{equation}
For this purpose, we shall use the matrix expression \eqref{Hgrandinfini}. Three kinds   of computations (listed (a), (b) and (c) below) are necessary according to the parameter with respect to which we differentiate.
\begin{enumerate}[(a)]
\item  We first differentiate with respect to $\underline{\omega}$ and we obtain
\[\dfrac{\partial\mathbb{H}_t}{\partial\theta_i} = \sum\limits_{k = 0}^{\infty} \mathbb{B}^k\dfrac{\partial \underline{c}_{t-k}}{\partial\theta_i}, \text{ for }i=1,\ldots, m, \]
and since $\partial \underline{c}_{t-k}(j_1)/\partial\theta_i = (0, \ldots, 1, \ldots, 0)'$ (the vector composed with 0 and 1 at the $j_1-$th position for $j_1=1,\ldots, m$), we have
\[\theta_i \dfrac{\partial \mathbb{H}_t(i_1)}{\partial\theta_i} = \sum\limits_{k=0}^{\infty}\sum\limits_{j_1=1}^m \mathbb{B}^k(i_1,j_1) \theta_i \ \dfrac{\partial \underline{c}_{t-k}(j_1)}{\partial\theta_i} \leq \sum\limits_{k=0}^{\infty} \sum\limits_{j_1=1}^m\mathbb{B}^k(i_1,j_1)\underline{c}_{t-k}(j_1) = \mathbb{H}_t(i_1) \]
where $ \mathbb{H}_t(i_1)$ is the $i_1-$component of $ \mathbb{H}_t$.
So we have
\begin{align}\label{no1}
&\dfrac{1}{\mathbb{H}_t(i_1)} \dfrac{\partial\mathbb{H}_t(i_1)}{\partial\theta_i} \leq \dfrac{1}{\theta_i}.
\end{align}
%
%
\item We differentiate with respect to $({\alpha_{1}^{+}} ', \ldots, {\alpha_q^+}', {\alpha_1^-}', \ldots, {\alpha_q^-}')'$. We have
\[\dfrac{\partial\mathbb{H}_t}{\partial\theta_i} = \sum\limits_{k = 0}^{\infty} \mathbb{B}^k\dfrac{\partial \underline{c}_{t-k}}{\partial\theta_i}\text{ for } i=m+1,\ldots, s_2,\]
with
\[\dfrac{\partial \underline{c}_{t-k}(j_1)}{\partial\theta_i} = \sum\limits_{l=1}^q\sum\limits_{j_2=1}^m \dfrac{\partial A_{l}^\pm(j_1,j_2)}{\partial\theta_i} (\pm{\varepsilon}_{j_2,t-l}^\pm)^{\delta_{j_2}}, \text{ for }j_1=1,\dots,m\]
where $\partial A_{l}^\pm(j_1, j_2) / \partial \theta_i$ is a null matrix or a  matrix whose entries are all zero except the one (equal to $1$) which is located at the same place of $\theta_i$. Thus
\[\theta_i \dfrac{\partial \mathbb{H}_t(i_1)}{\partial\theta_i} = \sum\limits_{k=0}^{\infty} \sum\limits_{j_1=1}^m\mathbb{B}^k(i_1,j_1) \theta_i \dfrac{\partial\underline{c}_{t-k}(j_1)}{\partial\theta_i} \leq \sum\limits_{k=0}^\infty \sum\limits_{j_1=1}^m \mathbb{B}^k(i_1,j_1)\underline{c}_{t-k}(j_1)= \mathbb{H}_t(i_1),\]
and  we have
\begin{align}
\label{no2}
&\dfrac{1}{\mathbb{H}_t(i_1)} \dfrac{\partial \mathbb{H}_t(i_1)}{\partial\theta_i} \leq \dfrac{1}{\theta_i}.
\end{align}
\item We now differentiate with respect to $(\beta'_{1}, \ldots, \beta'_{p})'$ and we have
\begin{equation}\label{derive-beta-connu}
\dfrac{\partial \mathbb{H}_t}{\partial\theta_i} = \sum\limits_{k=1}^\infty \left\{ \sum\limits_{j=1}^k \mathbb{B}^{j-1} \dfrac{\partial \mathbb{B}}{\partial\theta_i} \mathbb{B}^{k-j}\right\} \underline{c}_{t-k} \text{ for }i=s_2+1,\ldots, s_1.
\end{equation}
The matrix $ \partial \mathbb{B}/\partial\theta_i$ is a matrix whose entries are all 0, apart from a 1 located at the same place
as $\theta_i$ in $\mathbb{B}$. 
Thus $\theta_i \partial \mathbb{B}/\partial\theta_i \le \mathbb B$   and using \eqref{derive-beta-connu}, for all $i = s_2+1, \ldots, s_1$ we obtain
\[\theta_i \dfrac{\partial \mathbb{H}_t}{\partial\theta_i} \leq \sum\limits_{k=1}^\infty \left\{ \sum\limits_{j=1}^k \mathbb{B}^{j-1} \mathbb{B} \mathbb{B}^{k-j}\right\} \underline{c}_{t-k} = \sum\limits_{k=1}^\infty k \mathbb{B}^k \underline{c}_{t-k}. \]
Reasoning component by component, we have
\begin{equation} \label{inegaliteHi1}
\theta_i \dfrac{\mathbb{H}_t(i_1)}{\partial\theta_i} \leq \sum\limits_{k=1}^\infty k \sum\limits_{j_1 = 1}^m \mathbb{B}^k (i_1, j_1)\underline{c}_{t-k}(j_1). 
\end{equation}
We use $\mathbb{H}_t(i_1) \geq \omega + \sum\limits_{j_1 = 1}^m \mathbb{B}^k (i_1,j_1)\underline{c}_{t-k}(j_1),\;\forall k$, with $\omega = \inf\limits_{1\leq i\leq m}\underline{\omega}(i)$ and the relation $x/(1+x)\leq x^s$ which is valid for all $x\geq 0$ and $s \in [0,1]$. In view of \eqref{inegaliteHi1}, we obtain
\begin{equation*}
\begin{aligned}
\dfrac{\theta_i}{\mathbb{H}_t(i_1)} \dfrac{\partial \mathbb{H}_t(i_1)}{\partial \theta_i}  &\leq \sum\limits_{j_1=1}^m \sum\limits_{k=1}^\infty \dfrac{ k \mathbb{B}^k(i_1, j_1)\underline{c}_{t-k}(j_1)}{\omega + \mathbb{B}^k(i_1,j_1)\underline{c}_{t-k}(j_1)} \\
	& \leq \sum\limits_{j_1 = 1}^m \sum\limits_{k=1}^\infty k \left( \dfrac{\mathbb{B}^k(i_1,j_1)\underline{c}_{t-k}(j_1)}{\omega}\right)^{s/r_0}\\
	& \leq K \sum\limits_{j_1 = 1}^m \sum\limits_{k=1}^\infty k \rho^{sk/r_0} \underline{c}_{t-k}^{s/r_0}(j_1), 
\end{aligned}
\end{equation*}
with the constant $\rho$ which belongs to the interval $[0,1)$. Since $\underline{c}_{t-k}(\theta_0)$ has moments of any order (see Corollary \ref{cor2}) we have proved that
\begin{align}\label{prout}
\mathbb E \left  | \dfrac{\theta_i}{\mathbb{H}_t(i_1)} \dfrac{\partial \mathbb{H}_t(i_1)}{\partial \theta_i}  \right |^{r_0} &<\infty \ .
\end{align}
\medskip
\end{enumerate}
Since \eqref{no1} (as well as \eqref{no2} and \eqref{prout}) is an equivalent writing of \eqref{score.moment.r0} (as well as \eqref{no2} and \eqref{prout}), we deduce that \eqref{score.moment.r0} is true. By continuity of the functions that are involved on our estimations, the above inequalities are uniform on a neighborhood $V(\theta_0)$ of $\theta_0\in\overset{\circ}{\Theta}$: for all $i_1=1,\dots,m$ and all $i=1,\dots,s_1$ we have
\begin{align}
\label{unif1}
\mathbb{E}\sup_{\theta\in V(\theta_0)}\left\vert \dfrac{1}{{h}_{i_1,t}^{{\delta}_{i_1}/2}} \dfrac{\partial {h}_{i_1,t}^{{\delta}_{i_1}/2}}{\partial\theta_i}(\theta)\right\vert^{r_0}  & < \infty .
\end{align}

 \subsubsection{Asymptotic normality of the score vector.}\label{A33}
The proof is the same than the one in \cite{FZ-MAPGARCH}. It follows the following arguments:
\begin{itemize}
\item the process $\{ \partial l_t(\theta_0)/\partial \theta\}_t$ is stationary,
\item $\partial l_t(\theta_0) / \partial \theta$ is measurable with respect to the $\sigma-$field generated by $\{\eta_u, u<t\}$,
\item $\mathbb{E}\left[\dfrac{\partial l_t(\theta_0)}{\partial \theta} \vert \mathcal{F}_{t-1}\right] = 0$ thus we have a martingale-difference sequence,
\item Subsection \ref{A32}  implies that the matrix $I:=\mathbb{E}\left[\dfrac{\partial l_t(\theta_0)}{\partial\theta}\dfrac{\partial l_t(\theta_0)}{\partial\theta'}\right]$ is well defined.
\end{itemize}
Using the central limit theorem from \cite{billing} we obtain \eqref{NA1}.

\subsubsection{Convergence to $J$}\label{A34}
\subsubsection*{$\rightsquigarrow$ Expression of the second order derivatives of the log-likelihood}
\addcontentsline{toc}{paragraph}{$\rightsquigarrow$\ Expression of the second order derivatives of the log-likelihood}
We start from the expression \eqref{DL1} and \eqref{DL2} in order to compute the second order derivatives of the log-likelihood.  According to the index of $\theta_i$, we have three cases:
\begin{enumerate}[(i)]
\item For all $i,j = 1,\ldots, s_1$ we deduce from \eqref{DL1} that
\begin{align}
\label{D2-1}
\dfrac{\partial^2 l_t(\theta_0)}{\partial\theta_i\partial\theta_j} & = Tr(c_1 + c_2 + c_3)
\end{align}
with
\begin{equation*}
\begin{aligned}
c_1 &= D_{0t}^{-1}\underline{\varepsilon}_t\underline{\varepsilon}_t'D_{0t}^{-1}D_{0t}^{(i)}D_{0t}^{-1}R^{-1}D_{0t}^{-1}D_{0t}^{(j)} + D_{0t}^{-1}\underline{\varepsilon}_t\underline{\varepsilon}_t' D_{0t}^{-1} R^{-1}D_{0t}^{-1}D_{0t}^{(i)}D_{0t}^{-1}D_{0t}^{(j)}  \\
	& \hspace{3cm} -D_{0t}^{-1}\underline{\varepsilon}_t\underline{\varepsilon}_t'D_{0t}^{-1}R^{-1}D_{0t}^{-1}D_{0t}^{(i,j)} + D_t^{-1}\underline{\varepsilon}_t\underline{\varepsilon}_t'D_{0t}^{-1}R^{-1}D_{0t}^{-1} D_{0t}^{(j)}D_{0t}^{-1} D_{0t}^{(i)}, \\
c_2 &= D_{0t}^{-1}R^{-1}D_{0t}^{-1}D_{0t}^{(i)}D_{0t}^{-1}\underline{\varepsilon}_t\underline{\varepsilon}_t'D_{0t}^{-1}D_{0t}^{(j)} + D_{0t}^{-1}R^{-1} D_{0t}^{-1} \underline{\varepsilon}_t\underline{\varepsilon}_t' D_{0t}^{-1}D_{0t}^{(i)}D_{0t}^{-1}D_{0t}^{(j)}  \\
	& \hspace{3cm} -D_{0t}^{-1}R^{-1}D_{0t}^{-1}\underline{\varepsilon}_t\underline{\varepsilon}_t'D_{0t}^{-1}D_{0t}^{(i,j)} + D_t^{-1}R^{-1}D_{0t}^{-1}\underline{\varepsilon}_t\underline{\varepsilon}_t'D_{0t}^{-1} D_{0t}^{(j)}D_{0t}^{-1} D_{0t}^{(i)}, \\
c_3 & = -2D_{0t}^{-1} D_{0t}^{(i)}D_{0t}^{-1} D_{0t}^{(j)} + 2D_{0t}^{-1} D_{0t}^{(i,j)}.
\end{aligned}
\end{equation*}
%
%
\item For all $i, = 1,\ldots, s_1$ and $j=s_1+1, \ldots, s_0$  we have
\begin{align}\label{D2-2}
\dfrac{\partial^2 l_t(\theta_0)}{\partial\theta_i\partial\theta_j} & = Tr(c_4 + c_5)
\end{align}
with
\begin{align*}
c_4 & = R^{-1} D_{0t}^{-1} D_{0t}^{(i)} D_{0t}^{-1}\underline{\varepsilon}_t\underline{\varepsilon}_t'D_{0t}^{-1} R^{-1} R^{(j)} \ \text{and}\\
c_5  & =  R^{-1} D_{0t}^{-1} \underline{\varepsilon}_t\underline{\varepsilon}_t' D_{0t}^{-1} D_{0t}^{(i)} D_{0t}^{-1} R^{-1} R^{(j)}.
\end{align*}
This comes from the following computations. Starting from \eqref{DL1}, we use \eqref{DfA} and \eqref{Dtrace} and we obtain
\begin{align*}
-\dfrac{\partial}{\partial\theta_j} Tr\left(\underline{\varepsilon}_t\underline{\varepsilon}_t'D_{0t}^{-1} R^{-1} D_{0t}^{-1}D_{0t}^{(i)}D_{0t}^{-1}\right) & = Tr\left(R^{-1}D_{0t}^{-1} D_{0t}^{(i)}D_{0t}^{-1} \underline{\varepsilon}_t\underline{\varepsilon}_t'D_{0t}^{-1}R^{-1}R^{(j)}\right) \\
-\dfrac{\partial}{\partial\theta_j} Tr\left(R^{-1} D_{0t}^{-1}\underline{\varepsilon}_t\underline{\varepsilon}_t'D_{0t}^{-1} D_{0t}^{(i)} D_{0t}^{-1}\right) & =  Tr\left(R^{-1}D_{0t}^{-1} \underline{\varepsilon}_t\underline{\varepsilon}_t'D_{0t}^{-1} D_{0t}^{(i)}D_{0t}^{-1}R^{-1}R^{(j)}\right) \ \text{and} \\
2\dfrac{\partial }{\partial\theta_j}Tr\left(D_{0t}^{-1} \dfrac{\partial D_{0t}}{\partial\theta_i}\right) & = 0.
\end{align*}
\item For all $i,j = s_1+1,\ldots, s_0$ we have from \eqref{DL2}
\begin{align}\label{D2-3}
\dfrac{\partial^2 l_t(\theta_0)}{\partial\theta_i\partial\theta_j} & = Tr(c_6)
\end{align}
with
\begin{equation*}
\begin{aligned}
c_6 &= R^{-1}R^{(i)}R^{-1}D_{0t}^{-1}\underline{\varepsilon}_t\underline{\varepsilon}_t'D_{0t}^{-1}R^{-1}R^{(j)} + R^{-1}D_{0t}^{-1}\underline{\varepsilon}_t\underline{\varepsilon}_t'D_{0t}^{-1}R^{-1} R^{(i)} R^{-1}R^{(j)} \\
	& \hspace{3cm}- R^{-1}D_{0t}^{-1}\underline{\varepsilon}_t\underline{\varepsilon}_t'D_{0t}^{-1} R^{-1} R^{(i,j)} - R^{-1}R^{(i)}R^{-1}R^{(j)} + R^{-1} R^{(i,j)}.
\end{aligned}
\end{equation*}
Remark that, with our parametrization, $R^{(i,j)}=\partial^2 R/\partial\theta_i\partial\theta_j = 0$.
\end{enumerate}
Thanks to the three above cases, we remark that to control the second order derivatives,  it is sufficient to control the new term $\Vert D^{-1}_{0t} D_{0t}^{(i,j)} \Vert$. Indeed, all the other terms in $c_1,\dots,c_6$ can be controlled thanks to the results from Subsection \ref{A33}.
\subsubsection*{$\rightsquigarrow$ Existence of the moments of the second order derivatives of the log-likelihood}
\addcontentsline{toc}{paragraph}{$\rightsquigarrow$ Existence of the moments of the second order derivatives of the log-likelihood}
We take derivatives in the expression \eqref{DD} and we obtain
\begin{equation*}
\begin{aligned}
\dfrac{\partial^2 D_t(i_1,i_1)}{\partial\theta_{i}\partial\theta_{j}} &= \dfrac{1}{\delta_{i_1}}h_{i_1,t}^{1/2}(\theta)\left[ \dfrac{1}{h_{i_1,t}^{\delta_{i_1}/2}}
\dfrac{\partial h_{i_1,t}^{\delta_{i_1}/2}}{\partial \theta_{j}} \dfrac{1}{h_{i_1,t}^{\delta_{i}/2}} \dfrac{\partial h_{i_1,t}^{\delta_{i}/2}}{\partial \theta_{i}} \left(\dfrac{1}{\delta_{i_1}} - 1 \right) + \dfrac{1}{h_{i_1,t}^{\delta_{i_1}/2}} \dfrac{\partial^2 h_{i_1,t}^{\delta_{i_1}/2}}{\partial \theta_{i} \partial \theta_{j}} \right](\theta) .
\end{aligned}
\end{equation*}
It only remains to control the last term in the right hand side of the above identity. We will prove that
\begin{align}\label{D2h}
\mathbb{E}\left\vert \dfrac{1}{\underline{h}_{i_1,t}^{\underline{\delta}_{i_1}/2}} \dfrac{\partial^2 \underline{h}_{i_1,t}^{\underline{\delta}_{i_1}/2}}{\partial\theta_i\partial\theta_j}(\theta_0)\right\vert^{r_0} < \infty.
\end{align}
By \eqref{EcMat1}, we have for $i_1=1,\dots,m$:
$$ \mathbb{H}_t(i_1)= \sum\limits_{k=0}^{\infty} \sum\limits_{j_1=1}^m\mathbb{B}^k(i_1,j_1)\underline{c}_{t-k}(j_1) $$
and we remark that
\begin{align*}
0 & = \dfrac{\partial^2 \mathbb{H}_t(i_1)}{\partial \omega_i^2} = \dfrac{\partial^2 \mathbb{H}_t(i_1)}{\partial \omega_i\partial\omega_j} = \dfrac{\partial^2 \mathbb{H}_t(i_1)}{\partial \omega_i\partial\alpha^\pm_j} \\
0 & = \dfrac{\partial^2\mathbb{H}_t(i_1)}{\partial \alpha^\pm_i \partial\alpha^\pm_j} \ \text{and}\\
0 & = \dfrac{\partial^2\mathbb{H}_t(i_1)}{\partial (\beta_i)^2 }  \ .
\end{align*}
We also have
\begin{align}
\dfrac{\partial^2\mathbb{H}_t(i_1)}{\partial \omega_i \partial\beta_j} &=  \sum\limits_{k=1}^\infty\sum_{j_1=1}^m \left\{\sum\limits_{h=1}^k \mathbb{B}^{h-1}(i_1,j_1)\dfrac{\partial \mathbb{B}(i_1,j_1)}{\partial\beta_j}\mathbb{B}^{k-h}(i_1,j_1)\right\}\dfrac{\partial \underline{c}_{t-k}(j_1)}{\partial\omega_i}  \label{Domegabeta}\\
\dfrac{\partial^2\mathbb{H}_t(i_1)}{\partial \alpha^\pm_i \partial\beta_j} &= \sum\limits_{l=1}^q \sum\limits_{k=1}^\infty\sum_{j_1=1}^m\sum_{j_2=1}^m \left\{\sum\limits_{h=1}^k \mathbb{B}^{h-1}(i_1,j_1)\dfrac{\partial \mathbb{B}(i_1,j_1)}{\partial\beta_j}\mathbb{B}^{k-h}(i_1,j_1)\right\}\dfrac{\partial A^\pm_l(j_1,j_2)}{\partial\alpha^\pm_i}(\pm{\varepsilon}^\pm_{j_2,t-l})^{\delta_{j_2}} \label{Dalphabeta}
\end{align}
and
\begin{align}\label{Dbetabeta}
\dfrac{\partial^2\mathbb{H}_t}{\partial\beta_i \partial\beta_j } &= \sum\limits_{k=2}^\infty\left[\sum\limits_{l=2}^k \left\{\left(\sum\limits_{r=1}^{l-1} \mathbb{B}^{r-1}\dfrac{\partial \mathbb{B}}{\partial\beta_j}\mathbb{B}^{l-1-r}\right)\dfrac{\partial \mathbb{B}}{\partial\beta_i}\mathbb{B}^{k-l}\right\}\right. \nonumber \\
&\hspace{3cm} \left.+ \sum\limits_{l = 1}^{k-1}\left\{\mathbb{B}^{l-1}\dfrac{\partial \mathbb{B}}{\partial\beta_i} \left(\sum\limits_{r=1}^{k-l}\mathbb{B}^{r-1}\dfrac{\partial \mathbb{B}}{\partial\beta_j}\mathbb{B}^{k-l-r}\right)\right\}\right]\underline{c}_{t-k} \  .
\end{align}
We recall that the matrix $ \partial \mathbb{B}/\partial\beta_j$  is a matrix whose entries are all 0, apart from a 1 located at the same place
as $\beta_j$ in $\mathbb{B}$. 

Now we treat separately the three different expressions of the derivatives.
\begin{enumerate}[(a)]
\item For $i=1, \ldots, m$ we deduce from \eqref{Domegabeta} that
\begin{align*}
\omega_i\beta_j\dfrac{\partial^2\mathbb{H}_t(i_1)}{\partial \omega_i \partial\beta_j} &=  \sum\limits_{k=1}^\infty\sum_{j_1=1}^m  \left\{\sum\limits_{h=1}^k \mathbb{B}^{j-1}(i_1,j_1)\beta_j\dfrac{\partial \mathbb{B}(i_1,j_1)}{\partial\beta_j}\mathbb{B}^{k-h}(i_1,j_1)\right\}\omega_i\dfrac{\partial \underline{c}_{t-k}(j_1)}{\partial\omega_i}  \\
& \le \sum\limits_{k=1}^\infty\sum\limits_{j_1 = 1}^m k \mathbb{B}^k(i_1,j_1)  \underline{c}_{t-k}(j_1).
\end{align*}
Using the same arguments as for \eqref{prout}, we obtain that
\begin{align}\label{prout2}
\mathbb E \left |\dfrac{\omega_i\beta_j}{\mathbb{H}_t(i_1)} \dfrac{\partial^2\mathbb{H}_t(i_1)}{\partial\omega_i\partial\beta_j} \right |^{r_0} <\infty\ .
\end{align}
\item From \eqref{Dalphabeta} we have
\begin{align*}
\alpha_i^\pm\beta_j\dfrac{\partial^2\mathbb{H}_t(i_1)}{\partial \alpha^\pm_i \partial\beta_j} &= \sum\limits_{l=1}^q \sum\limits_{k=1}^\infty\sum_{j_1=1}^m\sum_{j_2=1}^m \left\{\sum\limits_{h=1}^k \mathbb{B}^{h-1}(i_1,j_1)\beta_j\dfrac{\partial \mathbb{B}(i_1,j_1)}{\partial\beta_j}\mathbb{B}^{k-h}(i_1,j_1)\right\}\alpha_i^\pm\dfrac{\partial A^\pm_l(j_1,j_2)}{\partial\alpha^\pm_i}(\pm{\varepsilon}^\pm_{j_2,t-l})^{\delta_{j_2}} \\
& \le  \sum\limits_{l=1}^q \sum\limits_{k=1}^\infty\sum_{j_1=1}^m \sum_{j_2=1}^mk\mathbb{B}^{k}(i_1,j_1) A^\pm_l(j_1,j_2)(\pm{\varepsilon}^\pm_{j_2,t-l})^{\delta_{j_2}}\\
& \le \sum\limits_{k=1}^\infty\sum\limits_{j_1 = 1}^m k \mathbb{B}^k(i_1,j_1)  \underline{c}_{t-k}(j_1)
\end{align*}
and we may proceed as in the proof of \eqref{prout} and we obtain that
\begin{align}\label{prout3}
\mathbb E \left |\dfrac{\alpha_i^\pm\beta_j}{\mathbb{H}_t(i_1)} \dfrac{\partial^2\mathbb{H}_t(i_1)}{\partial\alpha_i^\pm\partial\beta_j} \right |^{r_0} <\infty\ .
\end{align}
\item Using \eqref{Dbetabeta} we write
\begin{align*}
\beta_i\beta_j\dfrac{\partial^2\mathbb{H}_t(i_1)}{\partial\beta_i \partial\beta_j } &= \sum_{j_1=1}^m\sum\limits_{k=2}^\infty\left[\sum\limits_{l=2}^k \left\{\left(\sum\limits_{r=1}^{l-1} \mathbb{B}^{r-1}(i_1,j_1)\beta_j\dfrac{\partial \mathbb{B}(i_1,j_1)}{\partial\beta_j}\mathbb{B}^{l-1-r}(i_1,j_1)\right)\beta_i\dfrac{\partial \mathbb{B}(i_1,j_1)}{\partial\beta_i}\mathbb{B}^{k-l}(i_1,j_1)\right\}\right. \nonumber \\
&\hspace{0.2cm} \left.+ \sum\limits_{l = 1}^{k-1}\left\{\mathbb{B}^{l-1}(i_1,j_1)\beta_i\dfrac{\partial \mathbb{B}(i_1,j_1)}{\partial\beta_i} \left(\sum\limits_{r=1}^{k-l}\mathbb{B}^{r-1}(i_1,j_1)\beta_j\dfrac{\partial \mathbb{B}(i_1,j_1)}{\partial\beta_j}\mathbb{B}^{k-l-r}(i_1,j_1)\right)\right\}\right]\underline{c}_{t-k}(i_1) \\
& \leq \sum\limits_{j_1=1}^m\sum\limits_{k=2}^\infty \left[\sum\limits_{l=2}^k(l-1)\mathbb{B}^k(i_1,j_1) + \sum\limits_{l=1}^{k-1}(k-i)\mathbb{B}^k(i_1,j_1)\right]\underline{c}_{t-k}(j_1)\\
&\le  \sum\limits_{j_1=1}^m\sum\limits_{k=2}^\infty k(k-1)\mathbb{B}^k(i_1,j_1)\underline{c}_{t-k}(j_1).
\end{align*}
Arguing as before we deduce that
\begin{align}\label{prout4}
\mathbb E \left |\dfrac{\beta_i\beta_j}{\mathbb{H}_t(i_1)} \dfrac{\partial^2\mathbb{H}_t(i_1)}{\partial\beta_i\partial\beta_j} \right |^{r_0} <\infty\ .
\end{align}
\end{enumerate}
We deduce from \eqref{prout2}, \eqref{prout3} and  \eqref{prout4} that \eqref{D2h} is true.

Once again, by continuity of the involved functions, the above inequalities are uniform on a neighborhood $V(\theta_0)$ of $\theta_0\in\overset{\circ}{\Theta}$: for all $i_1=1,\dots,m$ and all $i,j,k=1,\dots,s_1$ we have
\begin{align}
\label{unif2}
\mathbb{E}\sup_{\theta\in V(\theta_0)}\left\vert \dfrac{1}{\underline{h}_{i_1,t}^{\underline{\delta}_{i_1}/2}} \dfrac{\partial^2 \underline{h}_{i_1,t}^{\underline{\delta}_{i_1}/2}}{\partial\theta_i\partial\theta_j}(\theta)\right\vert^{r_0} &< \infty \ .
\end{align}

\subsubsection*{$\rightsquigarrow$ Existence of the moments of the third order derivatives of the log-likelihood}
\addcontentsline{toc}{paragraph}{$\rightsquigarrow$ Existence of the moments of the third order derivatives of the log-likelihood}

First we write the quite heavy expressions of the third order derivatives derivatives with respect to the different parameters.
\begin{enumerate}[(i)]
\item For $i, j, k = 1,\ldots, s_1$, the derivatives with respect to $\theta_i$, $\theta_j$ and $\theta_k$ will correspond to the derivatives with respect to the parameters  $\underline{\omega}'$, $\mbox{vec}(A_l^\pm)'$ et $\mbox{vec}(B_{l'})'$. We obtain from \eqref{D2-1} that
\begin{equation}\label{Der-3-D-ijk}
\dfrac{\partial^3l_t(\theta)}{\partial \theta_i\partial \theta_j\partial \theta_k} = Tr\left( \big ( \mbox{$\sum_{p=1}^4$} c_{1p} + c_{2p} \big )  + c_{31}+c_{32}\right )
\end{equation}
with
\begin{align*}
c_{11} &= -D_t^{-1}D_t^{(i)}D_t^{-1}\underline{\varepsilon}_t\underline{\varepsilon}_t'D_t^{-1}D_t^{(j)}D_t^{-1}R^{-1}D_t^{-1}D_t^{(k)} - D_t^{-1}\underline{\varepsilon}_t\underline{\varepsilon}_t' D_t^{-1}D_t^{(i)}D_t^{-1}D_t^{(j)}D_t^{-1}R^{-1}D_t^{-1}D_t^{(k)}  \\
	& \hspace{1.5cm} + D_t^{-1}\underline{\varepsilon}_t\underline{\varepsilon}_t'D_t^{-1}D_t^{(i,j)}D_t^{-1}R^{-1}D_t^{-1}D_t^{(k)} - D_t^{-1}\underline{\varepsilon}_t\underline{\varepsilon}_t'D_t^{-1}D_t^{(j)}D_t^{-1}D_t^{(i)}D_t^{-1}R^{-1}D_t^{-1} D_t^{(k)} \\
	& \hspace{1.5cm} - D_t^{-1}\underline{\varepsilon}_t\underline{\varepsilon}_t'D_t^{-1}D_t^{(j)}D_t^{-1}R^{-1}D_t^{-1}D_t^{(i)}D_t^{-1}D_t^{(k)} + D_t^{-1}\underline{\varepsilon}_t\underline{\varepsilon}_t'D_t^{-1}R^{-1}D_t^{-1}D_t^{(j)}D_t^{-1}D_t^{(i,k)},\\
c_{12} &= -D_t^{-1}D_t^{(i)}D_t^{-1}\underline{\varepsilon}_t\underline{\varepsilon}_t'D_t^{-1}R^{-1}D_t^{-1}D_t^{(j)}D_t^{-1}D_t^{(k)} - D_t^{-1}\underline{\varepsilon}_t\underline{\varepsilon}_t'D_t^{-1}D_t^{(i)}D_t^{-1}R^{-1}D_t^{-1}D_t^{(j)}D_t^{-1}D_t^{(k)} \\
	& \hspace{1.5cm} -D_t^{-1}\underline{\varepsilon}_t\underline{\varepsilon}_t'D_t^{-1}R^{-1}D_t^{-1}D_t^{(i)}D_t^{-1}D_t^{(j)}D_t^{-1}D_t^{(k)} + D_t^{-1}\underline{\varepsilon}_t\underline{\varepsilon}_t'D_t^{-1}R^{-1}D_t^{-1}D_t^{(i,j)}D_t^{-1}D_t^{(k)} \\
	& \hspace{1.5cm} -D_t^{-1}\underline{\varepsilon}_t\underline{\varepsilon}_t'D_t^{-1}R^{-1}D_t^{-1}D_t^{(j)}D_t^{-1}D_t^{(i)}D_t^{-1}D_t^{(k)} + D_t^{-1}\underline{\varepsilon}_t\underline{\varepsilon}_t'D_t^{-1}R^{-1}D_t^{-1}D_t^{(j)}D_t^{-1}D_t^{(i,k)}, \\
c_{13} & = D_t^{-1}D_t^{(i)}D_t^{-1}\underline{\varepsilon}_t\underline{\varepsilon}_t'D_t^{-1}R^{-1}D_t^{-1}D_t^{(j,k)} + D_t^{-1}\underline{\varepsilon}_t\underline{\varepsilon}_t'D_t^{-1}D_t^{(i)}D_t^{-1}R^{-1}D_t^{-1}D_t^{(j,k)}\\
& \hspace{1.5cm} + D_t^{-1}\underline{\varepsilon}_t\underline{\varepsilon}_t'D_t^{-1}R^{-1}D_t^{-1}D_t^{(i)}D_t^{-1}D_t^{(j,k)} - D_t^{-1}\underline{\varepsilon}_t\underline{\varepsilon}_t'D_t^{-1}R^{-1}D_t^{-1}D_t^{(i,j,k)}, \\
c_{14} &= -D_t^{-1}D_t^{(i)}D_t^{-1}\underline{\varepsilon}_t\underline{\varepsilon}_t'D_t^{-1}R^{-1}D_t^{-1}D_t^{(j)}D_t^{-1}D_t^{(k)} - D_t^{-1}\underline{\varepsilon}_t\underline{\varepsilon}_t'D_t^{-1}D_t^{(i)}D_t^{-1}R^{-1}D_t^{-1}D_t^{(j)}D_t^{-1}D_t^{(k)}\\
& \hspace{1.5cm} - D_t^{-1}\underline{\varepsilon}_t\underline{\varepsilon}_t'D_t^{-1}R^{-1}D_t^{-1}D_t^{(i)}D_t^{-1}D_t^{(j)}D_t^{-1}D_t^{(k)} + D_t^{-1}\underline{\varepsilon}_t\underline{\varepsilon}_t'D_t^{-1}R^{-1}D_t^{-1}D_t^{(i,j)}D_t^{-1}D_t^{(k)}\\
& \hspace{1.5cm} - D_t^{-1}\underline{\varepsilon}_t\underline{\varepsilon}_t'D_t^{-1}R^{-1}D_t^{-1}D_t^{(j)}D_t^{-1}D_t^{(i)}D_t^{-1}D_t^{(k)} + D_t^{-1}\underline{\varepsilon}_t\underline{\varepsilon}_t'D_t^{-1}R^{-1}D_t^{-1}D_t^{(j)}D_t^{-1}D_t^{(i,k)},\\
c_{21} &= -D_t^{-1}D_t^{(i)}D_t^{-1}R^{-1}D_t^{-1}D_t^{(j)}D_t^{-1}\underline{\varepsilon}_t\underline{\varepsilon}_t'D_t^{-1}D_t^{(k)} - D_t^{-1}R^{-1} D_t^{-1}D_t^{(i)}D_t^{-1}D_t^{(j)}D_t^{-1}\underline{\varepsilon}_t\underline{\varepsilon}_t'D_t^{-1}D_t^{(k)}  \\
	& \hspace{1.5cm} + D_t^{-1}R^{-1}D_t^{-1}D_t^{(i,j)}D_t^{-1}\underline{\varepsilon}_t\underline{\varepsilon}_t'D_t^{-1}D_t^{(k)} - D_t^{-1}R^{-1}D_t^{-1}D_t^{(j)}D_t^{-1}D_t^{(i)}D_t^{-1}\underline{\varepsilon}_t\underline{\varepsilon}_t'D_t^{-1} D_t^{(k)} \\
	& \hspace{1.5cm} - D_t^{-1}R^{-1}D_t^{-1}D_t^{(j)}D_t^{-1}\underline{\varepsilon}_t\underline{\varepsilon}_t'D_t^{-1}D_t^{(i)}D_t^{-1}D_t^{(k)} + D_t^{-1}R^{-1}D_t^{-1}\underline{\varepsilon}_t\underline{\varepsilon}_t'D_t^{-1}D_t^{(j)}D_t^{-1}D_t^{(i,k)},\\
c_{22} &= -D_t^{-1}D_t^{(i)}D_t^{-1}R^{-1}D_t^{-1}\underline{\varepsilon}_t\underline{\varepsilon}_t'D_t^{-1}D_t^{(j)}D_t^{-1}D_t^{(k)} - D_t^{-1}R^{-1}D_t^{-1}D_t^{(i)}D_t^{-1}\underline{\varepsilon}_t\underline{\varepsilon}_t'D_t^{-1}D_t^{(j)}D_t^{-1}D_t^{(k)} \\
	& \hspace{1.5cm} -D_t^{-1}R^{-1}D_t^{-1}\underline{\varepsilon}_t\underline{\varepsilon}_t'D_t^{-1}D_t^{(i)}D_t^{-1}D_t^{(j)}D_t^{-1}D_t^{(k)} + D_t^{-1}R^{-1}D_t^{-1}\underline{\varepsilon}_t\underline{\varepsilon}_t'D_t^{-1}D_t^{(i,j)}D_t^{-1}D_t^{(k)} \\
	& \hspace{1.5cm} -D_t^{-1}R^{-1}D_t^{-1}\underline{\varepsilon}_t\underline{\varepsilon}_t'D_t^{-1}D_t^{(j)}D_t^{-1}D_t^{(i)}D_t^{-1}D_t^{(k)} + D_t^{-1}R^{-1}D_t^{-1}\underline{\varepsilon}_t\underline{\varepsilon}_t'D_t^{-1}D_t^{(j)}D_t^{-1}D_t^{(i,k)},\\
c_{23} & = D_t^{-1}D_t^{(i)}D_t^{-1}R^{-1}D_t^{-1}\underline{\varepsilon}_t\underline{\varepsilon}_t'D_t^{-1}D_t^{(j,k)} + D_tR^{-1}D_t^{-1}D_t^{(i)}D_t^{-1}\underline{\varepsilon}_t\underline{\varepsilon}_t'D_t^{-1}D_t^{(j,k)}\\
& \hspace{1.5cm} + D_t^{-1}R^{-1}D_t^{-1}\underline{\varepsilon}_t\underline{\varepsilon}_t'D_t^{-1}D_t^{(i)}D_t^{-1}D_t^{(j,k)} - D_t^{-1}R^{-1}D_t^{-1}\underline{\varepsilon}_t\underline{\varepsilon}_t'D_t^{-1}D_t^{(i,j,k)}, \\
c_{24} &= -D_t^{-1}D_t^{(i)}D_t^{-1}R^{-1}D_t^{-1}\underline{\varepsilon}_t\underline{\varepsilon}_t'D_t^{-1}D_t^{(j)}D_t^{-1}D_t^{(k)} - D_t^{-1}R^{-1}D_t^{-1}D_t^{(i)}D_t^{-1}\underline{\varepsilon}_t\underline{\varepsilon}_t'D_t^{-1}D_t^{(j)}D_t^{-1}D_t^{(k)}\\
& \hspace{1.5cm} - D_t^{-1}R^{-1}D_t^{-1}\underline{\varepsilon}_t\underline{\varepsilon}_t'D_t^{-1}D_t^{(i)}D_t^{-1}D_t^{(j)}D_t^{-1}D_t^{(k)} + D_t^{-1}R^{-1}D_t^{-1}\underline{\varepsilon}_t\underline{\varepsilon}_t'D_t^{-1}D_t^{(i,j)}D_t^{-1}D_t^{(k)}\\
& \hspace{1.5cm} - D_t^{-1}R^{-1}D_t^{-1}\underline{\varepsilon}_t\underline{\varepsilon}_t'D_t^{-1}D_t^{(j)}D_t^{-1}D_t^{(i)}D_t^{-1}D_t^{(k)} + D_t^{-1}R^{-1}D_t^{-1}\underline{\varepsilon}_t\underline{\varepsilon}_t'D_t^{-1}D_t^{(j)}D_t^{-1}D_t^{(i,k)},\\
c_{31} &= -2\left(-D_t^{-1}D_t^{(i)}D_t^{-1}D_t^{(j)}D_t^{-1}D_t^{(k)} + D_t^{-1}D_t^{(i,j)}D_t^{-1}D_t^{(k)}\right.\\
	& \hspace{1.5cm} \left. -D_t^{-1}D_t^{(j)}D_t^{-1}D_t^{(i)}D_t^{-1}D_t^{(k)} + D_t^{-1}D_t^{(j)}D_t^{-1}D_t^{(j,k)}\right),\\
c_{32} &= 2\left(-D_t^{-1}D_t^{(i)}D_t^{-1}D_t^{(j,k)} + D_t^{-1}D_t^{(i,j,k)}\right).
\end{align*}
\item For $i,j,k= s_1+1,\ldots, s_0$ which means that we differentiate with respect to the parameter $\rho$, we differentiate \eqref{D2-3} and we obtain
\begin{equation}\label{Der-3-R-ijk}
\dfrac{\partial^3l_t(\theta)}{\partial \theta_i\partial \theta_j\partial \theta_k} = Tr(c_{61} + c_{62} + c_{63} + c_{64} + c_{65})
\end{equation}
with
\begin{equation*}
\begin{aligned}
c_{61} &= -R^{-1}R^{(i)}R^{-1}R^{(j)}D_t^{-1}\underline{\varepsilon}_t\underline{\varepsilon}_t'D_t^{-1}R^{-1}R^{(k)} + R^{-1}R^{(i,j)}D_t^{-1}\underline{\varepsilon}_t\underline{\varepsilon}_t'D_t^{-1}R^{-1}R^{(k)}  \\
	& \hspace{1.5cm} -R^{-1}R^{(j)}D_t^{-1}\underline{\varepsilon}_t\underline{\varepsilon}_t'D_t^{-1}R^{-1}R^{(i)}R^{-1}R^{(k)} + R^{-1}R^{(j)}D_t^{-1}\underline{\varepsilon}_t\underline{\varepsilon}_t'D_t^{-1}R^{-1}R^{(i,k)} ,\\
c_{62} &=  -R^{-1}R^{(i)}R^{-1}D_t^{-1}\underline{\varepsilon}_t\underline{\varepsilon}_t'D_t^{-1} R^{-1}R^{(j)}R^{-1}R^{(k)} - R^{-1}D_t^{-1}\underline{\varepsilon}_t\underline{\varepsilon}_t'D_t^{-1}R^{-1}R^{(i)}R^{-1}R^{(j)}R^{-1}R^{(k)}\\
	& \hspace{1.5cm} + R^{-1}D_t^{-1}\underline{\varepsilon}_t\underline{\varepsilon}_t'D_t^{-1}R^{-1}R^{(i,j)}R^{-1}R^{(k)} - R^{-1}D_t^{-1}\underline{\varepsilon}_t\underline{\varepsilon}_t'D_t^{-1}R^{-1}R^{(j)}R^{-1}R^{(i)}R^{-1}R^{(k)} \\
	& \hspace{1.5cm} + R^{-1}D_t^{-1}\underline{\varepsilon}_t\underline{\varepsilon}_t'D_t^{-1}R^{-1}R^{(j)}R^{-1}R^{(i,k)},\\
c_{63} & = R^{-1}R^{(i)}R^{-1}D_t^{-1}\underline{\varepsilon}_t\underline{\varepsilon}_t'R^{-1}R^{(j,k)} + R^{-1}D_t^{-1}\underline{\varepsilon}_t\underline{\varepsilon}_t'R^{-1}R^{(i)}R^{-1}R^{(j,k)} -R^{-1}D_t^{-1}\underline{\varepsilon}_t\underline{\varepsilon}_t'R^{-1}R^{(i,j,k)},\\
c_{64} &= R^{-1}R^{(i)}R^{-1}R^{(j)}R^{-1}R^{(k)} - R^{-1}R^{(i,j)}R^{-1}R^{(k)} + R^{-1}R^{(j)}R^{-1}R^{(i)}R^{-1}R^{(k)} - R^{-1}R^{(j)}R^{-1}R^{(i,k)}, \\
c_{65} &= -R^{-1}R^{(i)}R^{-1}R^{(j,k)} + R^{-1}R^{(i,j,k)}.
\end{aligned}
\end{equation*}
\item For $i, j = 1,\ldots, s_1$ (differentiation with respect to $\underline{\omega}'$, $\mbox{vec}(A_l^\pm)'$ and $\mbox{vec}(B_{l'})'$) and for $k=s_1+1,\ldots, s_0$ (differentiation with respect to $\rho$), we obtain from \eqref{D2-2}
\begin{equation}\label{Der-3-D-ij}
\dfrac{\partial^3l_t(\theta)}{\partial \theta_i\partial \theta_j\partial \theta_k} = Tr(c_{41} + c_{51})
\end{equation}
where
\begin{equation*}
\begin{aligned}
c_{41} &= -R^{-1}D_t^{-1}D_t^{(i)}D_t^{-1}D_t^{(j)}D_t^{-1}\underline{\varepsilon}_t\underline{\varepsilon}_t'D_t^{-1}R^{-1}R^{(k)} + R^{-1}D_t^{-1}D_t^{(i,j)}D_t^{-1}\underline{\varepsilon}_t\underline{\varepsilon}_t'D_t^{-1}R^{-1}R^{(k)}\\
	& \hspace{1cm} -R^{-1}D_t^{-1}D_t^{(j)}D_t^{-1}D_t^{(i)}D_t^{-1}\underline{\varepsilon}_t\underline{\varepsilon}_t'D_t^{-1}R^{-1}R^{(k)} - R^{-1}D_t^{-1}D_t^{(j)}D_t^{-1}\underline{\varepsilon}_t\underline{\varepsilon}_t'D_t^{-1}D_t^{(i)}D_t^{-1}R^{-1}R^{(k)},\\
c_{51} &= -R^{-1}D_t^{-1}D_t^{(i)}D_t^{-1}\underline{\varepsilon}_t\underline{\varepsilon}_t'D_t^{-1}D_t^{(j)}D_t^{-1}R^{-1}R^{(k)} - R^{-1}D_t^{-1}\underline{\varepsilon}_t\underline{\varepsilon}_t'D_t^{-1}D_t^{(i)}D_t^{-1}D_t^{(j)}D_t^{-1}R^{-1}R^{(k)} \\
	& \hspace{1.5cm} + R^{-1}D_t^{-1}\underline{\varepsilon}_t\underline{\varepsilon}_t'D_t^{-1}D_t^{(i,j)}D_t^{-1}R^{-1}R^{(k)} - R^{-1}D_t^{-1}\underline{\varepsilon}_t\underline{\varepsilon}_t'D_t^{-1}D_t^{(j)}D_t^{-1}D_t^{(i)}D_t^{-1}R^{-1}R^{(k)}.
\end{aligned}
\end{equation*}
\item In the same way, from \eqref{D2-3} we obtain for $i= 1,\ldots, s_1$ (derivatives with respect to $\underline{\omega}'$, $\mbox{vec}(A_l^\pm)'$ and $\mbox{vec}(B_{l'})'$) and for $j,k=s_1+1,\ldots, s_0$ (derivatives with respect to  $\rho$):
\begin{equation}\label{Der-3-D-i}
\dfrac{\partial^3l_t(\theta)}{\partial \theta_i\partial \theta_j\partial \theta_k} = Tr(c_{61}' + c_{62}' + c_{63}')
\end{equation}
with
\begin{equation*}
\begin{aligned}
c_{61}' &= -R^{-1}R^{(j)}D_t^{-1}D_t^{(i)}D_t^{-1}\underline{\varepsilon}_t\underline{\varepsilon}_t'D_t^{-1}R^{-1}R^{(k)} -R^{-1}D_t^{-1}\underline{\varepsilon}_t\underline{\varepsilon}_t'D_t^{-1}D_t^{(i)}D_t^{-1}R^{-1}R^{(j)}R^{-1}R^{(k)},\\
c_{62}' &= -R^{-1}D_t^{-1}D_t^{(i)}D_t^{-1}\underline{\varepsilon}_t\underline{\varepsilon}_t'D_t^{-1}R^{-1}R^{(j)}R^{-1}R^{(k)} - R^{-1}D_t^{-1}\underline{\varepsilon}_t\underline{\varepsilon}_t'D_t^{-1}D_t^{(i)}D_t^{-1}R^{-1}R^{(j)}R^{-1}R^{(k)}, \\
c_{63}' &= R^{-1}D_t^{-1}D_t^{(i)}D_t^{-1}\underline{\varepsilon}_t\underline{\varepsilon}_t'D_t^{-1}R^{-1}R^{(i,j)} + R^{-1}D_t^{-1}\underline{\varepsilon}_t\underline{\varepsilon}_t'D_t^{-1}D_t^{(i)}D_t^{-1}R^{-1}R^{(i,j)}.
\end{aligned}
\end{equation*}
We remark that the last two terms of $c_6$ in \eqref{D2-3}
are not composed with the matrix $D_t$ and thus their derivatives vanishes.
\end{enumerate}
In order to estimate the moments of order $r_0$ of the third order derivatives, we need to study the term  $D_t^{-1}D_{0t}$ which appears, for example, in the third term of $c_{11}$ namely  $c_{11}^{(3)}=D_t^{-1}\underline{\varepsilon}_t\underline{\varepsilon}_t'D_t^{-1}D_t^{(i,j)}D_t^{-1}R^{-1}D_t^{-1}D_t^{(k)}$. Indeed, using the fact that $\underline{\varepsilon}_t=D_{0t} \tilde{\eta}_t$, we may write
\begin{align}
\mathbb{E}\left[\sup\limits_{\theta\in V(\theta_0)}\left\vert Tr\left(c_{11}^{(3)} \right)\right\vert\right] &\leq K \mathbb{E}\Vert \tilde{\eta}_t\tilde{\eta}_t'\Vert \mathbb{E}\left[\sup\limits_{\theta\in V(\theta_0)}\Vert D_t^{-1}D_{0t}\Vert^2\Vert R^{-1}\Vert D_t^{-1}D_t^{(i,j)} \Vert \Vert D_t^{-1}D_t^{(k)}\Vert\right] \nonumber\\
& \leq K\mathbb{E}\left[\sup\limits_{\theta\in V(\theta_0)} \Vert D_t^{-1}D_{0t}\Vert^2 \Vert D_t^{-1}D_t^{(i,j)} \Vert \Vert D_t^{-1}D_t^{(k)}\Vert\right].
\end{align}
Consequently we have to prove that for any $r_0 \geq 1$,
\begin{equation}\label{DDesti}
\mathbb{E}\left[\sup\limits_{\theta \in V(\theta_0)} \Vert D_t^{-1}D_{0t}\Vert^{r_0}\right] < \infty.
\end{equation}
Letting $\mathbb{B}_0=\mathbb{B}(\theta_0)$, by \eqref{Hgrandinfini} the $ i_1^{\text{th}}$ component  of $\underline{h}_{t}^{\underline{\delta}/2}$ equals
\begin{align*}
{h}_{i_1,t}^{\delta_{i_1}/2}(\theta_0) &= c_0 + \sum\limits_{k=0}^\infty\sum\limits_{j_1=1}^m\sum\limits_{j_2=1}^m\sum\limits_{i=1}^q \mathbb{B}_0^k(i_1,j_1)\left\{A_{0i}^+(j_1,j_2)(\varepsilon_{j_2,t-k-i}^+)^{\delta_{j_2}}
+ A_{0i}^-(j_1,j_2)(-\varepsilon_{j_2,t-k-i}^-)^{\delta_{j_2}}\right\}
\end{align*}
where  $c_0$ is a strictly positive constant. For a sufficiently small neighborhood $V(\theta_0)$ of $\theta_0$, we have
\begin{equation*}
\sup\limits_{\theta \in V(\theta_0)}\dfrac{A_{0i}^+(j_1,j_2)}{A_i^+(j_1,j_2)}< K,\qquad \sup\limits_{\theta \in V(\theta_0)}\dfrac{A_{0i}^-(j_1,j_2)}{A_i^-(j_1,j_2)}< K,\qquad \sup\limits_{\theta \in V(\theta_0)}\dfrac{\mathbb{B}_0^k(i_1,j_1)}{\mathbb{B}^k(i_1,j_1)}\leq (1 + \xi)^k
\end{equation*}
for any $i_1, j_1, j_2 \in \{1, \ldots, m\}$ and for all integer $k$ and all $\xi >0$. Moreover, the coefficients $A_{0i}^+(j_1,j_2)$ et $A_{0i}^-(j_1,j_2)$ are bounded below by a constant $c>0$ uniformly in $\theta$ on  $V(\theta_0)$. We thus obtain that
\begin{align*}
\sup\limits_{\theta\in V(\theta_0)} \dfrac{{h}_{i_1,t}^{\delta_{i_1}/2}(\theta_0)}{{h}_{i_1,t}^{\delta_{i_1}/2}(\theta)} &\leq K + K \sum\limits_{k=0}^\infty \sum\limits_{j_1=1}^m\sum\limits_{j_2=1}^m\sum\limits_{i=1}^q \left\{\dfrac{(1+\xi)^k\mathbb{B}^k(i_1,j_1)(\varepsilon_{j_2,t-k-i}^+)^{\delta_{j_2}}}{\omega + c\mathbb{B}^k(i_1,j_1)(\varepsilon_{j_2,t-k-i}^+)^{\delta_{j_2}}}
\right.\\
 &\left. \hspace{4cm}\qquad\qquad
 + \dfrac{(1+\xi)^k\mathbb{B}^k(i_1,j_1)(-\varepsilon_{j_2,t-k-i}^-)^{\delta_{j_2}}}{\omega + c\mathbb{B}^k(i_1,j_1)(-\varepsilon_{j_2,t-k-i}^-)^{\delta_{j_2}}}\right\}\\
 &\leq K + K \sum\limits_{j_2=1}^m\sum\limits_{i=1}^q\sum\limits_{k=0}^\infty (1+\xi)^k\rho^{ks}\left|\varepsilon_{j_2,t-k-i}\right|^{\delta_{j_2}s},
\end{align*}
for some $\rho\in [0,1[$, all $\xi>0$ and all $s\in [0,1]$. We then deduce that
\begin{align*}
\sup\limits_{\theta\in V(\theta_0)} \dfrac{{h}_{i_1,t}^{1/2}(\theta_0)}{{h}_{i_1,t}^{1/2}(\theta)} &=\sup\limits_{\theta\in V(\theta_0)} \left(\dfrac{{h}_{i_1,t}^{\delta_{i_1}/2}(\theta_0)}{{h}_{i_1,t}^{\delta_{i_1}/2}(\theta)}\right)^{1/\delta_{i_1}}\leq\left( K + K \sum\limits_{j_2=1}^m\sum\limits_{i=1}^q\sum\limits_{k=0}^\infty (1+\xi)^k\rho^{ks}\left|\varepsilon_{j_2,t-k-i}\right|^{\delta_{j_2}s}\right)^{1/\delta_{i_1}}.
\end{align*}
We distinguish two cases. If $\delta_{i_1}>1$, the concavity of the function $x\mapsto x^{\delta_{i_1}}$ on $[0,\infty[$ implies, by the Jensen inequality, that
\begin{align}\label{raph0}
\mathbb{E}\left[\sup\limits_{\theta\in V(\theta_0)} \dfrac{{h}_{i_1,t}^{1/2}(\theta_0)}{{h}_{i_1,t}^{1/2}(\theta)} \right] \leq\left( K + K \sum\limits_{j_2=1}^m\sum\limits_{i=1}^q\sum\limits_{k=0}^\infty (1+\xi)^k\rho^{ks}\mathbb{E}\left|\varepsilon_{j_2,t-k-i}\right|^{\delta_{j_2}s}\right)^{1/\delta_{i_1}}<\infty,
\end{align}
by using Corollary \ref{cor2}. Now, when $\delta_{i_1}=1$, Corollary \ref{cor2} entails that
\begin{align}\label{raph0bis}
\mathbb{E}\left[\sup\limits_{\theta\in V(\theta_0)} \dfrac{{h}_{i_1,t}^{1/2}(\theta_0)}{{h}_{i_1,t}^{1/2}(\theta)} \right] \leq K + K \sum\limits_{j_2=1}^m
\sum\limits_{i=1}^q\sum\limits_{k=0}^\infty (1+\xi)^k\rho^{ks}\mathbb{E}\left|\varepsilon_{j_2,t-k-i}\right|^{\delta_{j_2}s}<\infty.
\end{align}
In view of \eqref{raph0} and \eqref{raph0bis}, for all $r_0\ge 1$, we deduce that
\begin{equation}\label{h0-h1}
\mathbb{E}\left[\sup\limits_{\theta\in V(\theta_0)} \left\vert \dfrac{{h}_{i_1,t}^{1/2}(\theta_0)}{{h}_{i_1,t}^{1/2}(\theta)}\right\vert^{r_0} \right]<\infty
\end{equation}
and \eqref{DDesti} is true. \\
We also need to control the term $\Vert D_t^{-1}D_t^{(i,j,k)}\Vert$. There are many other terms that involve derivatives of order one and two that we may control thanks our previous estimations. Moreover $R^{(i,j,k)}=0$ because $R^{(i,j)}=0$.\\
The third order derivative of the matrix $D_t$ with respect to the parameters  $\theta_i$, $\theta_j$ and $\theta_k$ has the following expression for $i,j,k=1,\ldots, s_1$
\begin{align*}
\dfrac{\partial^3 D_t(i_1,i_1)}{\partial \theta_i \partial \theta_j\partial \theta_k} &= \dfrac{1}{\delta_{i_1}^2}h_{i_1,t}^{1/2}\dfrac{1}{h_{i_1,t}^{\delta_{i_1}/2}}\dfrac{\partial h_{i_1,t}^{\delta_{i_1}/2}}{\partial\theta_k}\left[\dfrac{1}{h_{i_1,t}^{\delta_{i_1}/2}}
\dfrac{\partial h_{i_1,t}^{\delta_{i_1}/2}}{\partial \theta_{j}} \dfrac{1}{h_{i_1,t}^{\delta_{i_1}/2}} \dfrac{\partial h_{i_1,t}^{\delta_{i_1}/2}}{\partial \theta_{i}} \left(\dfrac{1}{\delta_{i_1}} - 1 \right) + \dfrac{1}{h_{i_1,t}^{\delta_{i_1}/2}} \dfrac{\partial^2 h_{i_1,t}^{\delta_{i_1}/2}}{\partial \theta_{i} \partial \theta_{j}} \right] \\
&\qquad + \dfrac{1}{\delta_{i_1}}h_{i_1,t}^{1/2} \left[ -\dfrac{1}{h_{i_1,t}^{\delta_{i_1}/2}}\dfrac{\partial h_{i_1,t}^{\delta_{i_1}/2}}{\partial \theta_k}\dfrac{1}{h_{i_1,t}^{\delta_{i_1}/2}}\dfrac{\partial h_{i_1,t}^{\delta_{i_1}/2}}{\partial \theta_j}\dfrac{1}{h_{i_1,t}^{\delta_{i_1}/2}}\dfrac{\partial h_{i_1,t}^{\delta_{i,1}/2}}{\partial \theta_i}\left(\dfrac{1}{\delta_{i_1}} -1 \right) \right.\\
& + \dfrac{1}{h_{i_1,t}^{\delta_{i_1}/2}}\dfrac{\partial^2h_{i_1,t}^{\delta_{i_1}/2}}{\partial\theta_j\partial\theta_k}\dfrac{1}{h_{i_1,t}^{\delta_{i_1}/2}}\dfrac{\partial h_{i_1,t}^{\delta_{i_1}/2}}{\partial \theta_i}\left(\dfrac{1}{\delta_{i_1}}-1\right) + \dfrac{1}{h_{i_1,t}^{\delta_{i_1}/2}}\dfrac{\partial h_{i_1,t}^{\delta_{i_1}/2}}{\partial \theta_j} \dfrac{1}{h_{i_1,t}^{\delta_{i_1}/2}} \dfrac{\partial^2 h_{i_1,t}^{\delta_{i_1}/2}}{\partial \theta_i \partial\theta_k}\left(\dfrac{1}{\delta_{i_1}}-1\right)\\
&\left.\qquad - \dfrac{1}{h_{i_1,t}^{\delta_{i_1}/2}} \dfrac{\partial h_{i_1,t}{\delta_{i_1}/2}}{\partial \theta_k}\dfrac{1}{h_{i_1,t}^{\delta_{i_1}/2}} \dfrac{\partial^2 h_{i_1,t}^{\delta_{i_1}/2}}{\partial \theta_i\partial\theta_j} + \dfrac{1}{h_{i_1,t}^{\delta_{i_1}/2}}\dfrac{\partial^3h_{i_1,t}^{\delta_{i_1}/2}}{\partial \theta_i\partial\theta_j\partial\theta_k}\right].
\end{align*}
The terms in which the first and second order derivatives of  $\underline{h}_t^{\delta_{i_1}/2}$ are involved are already controlled thanks to \eqref{unif1} and \eqref{unif2}. Thus it remains to prove that
\begin{equation}\label{ss1}
\mathbb{E}\left\vert \dfrac{1}{h_{i_1,t}^{\delta_{i_1}/2}} \dfrac{\partial^3 h_{i_1,t}^{\delta_{i_1}/2}}{\partial \theta_i\partial \theta_j\partial \theta_k}(\theta)\right\vert^{r_0}<\infty.
\end{equation}
Starting from \eqref{Hgrandinfini}
\begin{align*}
 \mathbb{H}_t(i_1)= \sum\limits_{k=0}^{\infty} \sum\limits_{j_1=1}^m\mathbb{B}^k(i_1,j_1)\underline{c}_{t-k}(j_1),
\end{align*}
one may express the derivatives with respect to the different parameters.
We only have to treat the derivatives
$$\frac{\partial^3 \mathbb{H}_t(i_1)}{\partial \theta_i\partial \beta_j\partial \beta_k}$$
when  $\theta_i \neq \beta_j$ and $\theta_i\neq \beta_k$ because the other derivatives vanish.

There are three cases.
\begin{enumerate}[(i)]
\item For $i=1,\ldots, m$ (this means that we differentiate with respect to the parameter $\underline{\omega}$) and for fixed $j$ and $k$, it holds
\begin{align}\label{Der-3-Omega}
\dfrac{\partial^3\mathbb{H}_t}{\partial \omega_i \partial\beta_j \partial\beta_k } &= \sum\limits_{k'=2}^\infty\left[\sum\limits_{l=2}^{k'} \left\{\left(\sum\limits_{r=1}^{l-1} \mathbb{B}^{r-1}\dfrac{\partial \mathbb{B}}{\partial\beta_k}\mathbb{B}^{l-1-r}\right)\dfrac{\partial \mathbb{B}}{\partial\beta_j}\mathbb{B}^{k'-l}\right\}\right. \nonumber \\
&\hspace{3cm} \left.+ \sum\limits_{l = 1}^{k'-1}\left\{\mathbb{B}^{l-1}\dfrac{\partial \mathbb{B}}{\partial\beta_j} \left(\sum\limits_{r=1}^{k-l}\mathbb{B}^{r-1}\dfrac{\partial \mathbb{B}}{\partial\beta_k}\mathbb{B}^{k'-l-r}\right)\right\}\right]\dfrac{\partial \underline{c}_{t-k'}}{\partial\omega_i}.
\end{align}
Arguing as we did for the second order derivatives, we obtain
\small
\begin{align*}
\omega_i\beta_j\beta_k\dfrac{\partial^3\mathbb{H}_t(i_1)}{\partial \omega_i\partial\beta_j \partial\beta_k } &= \sum_{j_1=1}^m\sum\limits_{k'=2}^\infty\left[\sum\limits_{l=2}^{k'} \left\{\left(\sum\limits_{r=1}^{l-1} \mathbb{B}^{r-1}(i_1,j_1)\beta_k\dfrac{\partial \mathbb{B}(i_1,j_1)}{\partial\beta_k}\mathbb{B}^{l-1-r}(i_1,j_1)\right)\beta_j\dfrac{\partial \mathbb{B}(i_1,j_1)}{\partial\beta_j}\mathbb{B}^{k'-l}(i_1,j_1)\right\}\right. \\
&\hspace{-1cm} \left.+ \sum\limits_{l = 1}^{k'-1}\left\{\mathbb{B}^{l-1}(i_1,j_1)\beta_j\dfrac{\partial \mathbb{B}(i_1,j_1)}{\partial\beta_j} \left(\sum\limits_{r=1}^{k'-l}\mathbb{B}^{r-1}(i_1,j_1)\beta_k\dfrac{\partial \mathbb{B}(i_1,j_1)}{\partial\beta_k}\mathbb{B}^{k'-l-r}(i_1,j_1)\right)\right\}\right]\omega_i\dfrac{\partial \underline{c}_{t-k'}(j_1)}{\partial\omega_i}. \\
& \leq \sum\limits_{j_1=1}^m\sum\limits_{k'=2}^\infty \left[\sum\limits_{l=2}^{k'}(l-1)\mathbb{B}^{k'}(i_1,j_1) + \sum\limits_{l=1}^{k'-1}(k'-i)\mathbb{B}^{k'}(i_1,j_1)\right]\underline{c}_{t-k'}(j_1)\\
&\le  \sum\limits_{j_1=1}^m\sum\limits_{k'=2}^\infty k'(k'-1)\mathbb{B}^{k'}(i_1,j_1)\underline{c}_{t-k'}(j_1).
\end{align*}
\normalsize
Consequently, it holds
\begin{align}\label{Maj-Der-3-Omega}
\mathbb E \left |\dfrac{\omega_i\beta_j\beta_k}{\mathbb{H}_t(i_1)} \dfrac{\partial^3\mathbb{H}_t(i_1)}{\partial \omega_i\partial\beta_j\partial\beta_k} \right |^{r_0} <\infty.
\end{align}
\item For $i=m+1,\ldots,s_2$ (corresponding to differentiation with respect to  $\mbox{vec}(A^\pm_l)$) and for fixed $j$ and $k$, we have
\small
\begin{align}\label{Der-3-Alpha}
\dfrac{\partial^3\mathbb{H}_t(i_1)}{\partial \alpha_i^\pm\partial\beta_j\partial\beta_k} &=\sum\limits_{l'=1}^q\sum\limits_{k'=2}^\infty\sum\limits_{j_1=1}^m\sum\limits_{j_2=1}^m\left[\sum\limits_{l=2}^{k'}\left\{\left(\sum\limits_{r=1}^{l-1} \mathbb{B}^{r-1}(i_1,j_1)\dfrac{\partial \mathbb{B}(i_1,j_1)}{\partial\beta_k}\mathbb{B}^{l-1-r}(i_1,j_1)\right)\dfrac{\partial \mathbb{B}(i_1,j_1)}{\partial\beta_j}\mathbb{B}^{k-l}(i_1,j_1)\right\}\right. \nonumber\\
&\hspace{-1.5cm} \left.+ \sum\limits_{l = 1}^{k-1}\left\{\mathbb{B}^{l-1}(i_1,j_1)\dfrac{\partial \mathbb{B}(i_1,j_1)}{\partial\beta_j} \left(\sum\limits_{r=1}^{k'-l}\mathbb{B}^{r-1}(i_1,j_1)\dfrac{\partial \mathbb{B}(i_1,j_1)}{\partial\beta_k}\mathbb{B}^{k'-l-r}(i_1,j_1)\right)\right\}\right]\dfrac{\partial A_{i}^\pm(j_1,j_2)}{\partial \alpha_i^\pm}(\pm\varepsilon_{j_2,t-l'}^\pm)^{\delta_{j_2}},
\end{align}
\normalsize
and we write that
\footnotesize
\begin{align*}
\alpha_i^\pm\beta_j\beta_k\dfrac{\partial^3\mathbb{H}_t(i_1)}{\partial \alpha_i^\pm\partial\beta_j\partial\beta_k} &=\sum\limits_{l'=1}^q\sum\limits_{k'=2}^\infty\sum\limits_{j_1=1}^m\sum\limits_{j_2=1}^m\left[\sum\limits_{l=2}^{k'}\left\{\left(\sum\limits_{r=1}^{l-1} \mathbb{B}^{r-1}(i_1,j_1)\beta_k\dfrac{\partial \mathbb{B}(i_1,j_1)}{\partial\beta_k}\mathbb{B}^{l-1-r}(i_1,j_1)\right)\beta_j\dfrac{\partial \mathbb{B}(i_1,j_1)}{\partial\beta_j}\mathbb{B}^{k-l}(i_1,j_1)\right\}\right. \nonumber\\
&\hspace{-1.7cm} \left.+ \sum\limits_{l = 1}^{k-1}\left\{\mathbb{B}^{l-1}(i_1,j_1)\beta_j\dfrac{\partial \mathbb{B}(i_1,j_1)}{\partial\beta_j} \left(\sum\limits_{r=1}^{k'-l}\mathbb{B}^{r-1}(i_1,j_1)\beta_k\dfrac{\partial \mathbb{B}(i_1,j_1)}{\partial\beta_k}\mathbb{B}^{k'-l-r}(i_1,j_1)\right)\right\}\right]\alpha_i^\pm\dfrac{\partial A_{i}^\pm(j_1,j_2)}{\partial \alpha_i^\pm}(\pm\varepsilon_{j_2,t-l'}^\pm)^{\delta_{j_2}}\\
&\leq \sum\limits_{l'=1}^q\sum\limits_{k'=2}^\infty \sum\limits_{j_1=1}^m\sum\limits_{j_2=1}^m\left[\sum\limits_{l=2}^{k'}(l-1)\mathbb{B}^{k'}(i_1,j_1) + \sum\limits_{l=1}^{k'-1}(k'-1)\mathbb{B}^{k'}(i_1,j_1)\right]A_{i}^\pm(j_1,j_2)(\pm\varepsilon_{j_2,t-l'}^\pm)^{\delta_{j_2}}\\
&\leq \sum\limits_{j_1=1}^m\sum\limits_{k'=2}^\infty k'(k'-1)\mathbb{B}^{k'}(i_1,j_1)\underline{c}_{t-k'}(j_1).
\end{align*}
\normalsize
Hence we deduce that
\begin{align}\label{Maj-Der-3-Alpha}
\mathbb E \left |\dfrac{\alpha_i^\pm\beta_j\beta_k}{\mathbb{H}_t(i_1)} \dfrac{\partial^3\mathbb{H}_t(i_1)}{\partial \alpha_i^\pm\partial\beta_j\partial\beta_k} \right |^{r_0} <\infty.
\end{align}
\item For $i=s_2+1,\ldots, s_1$ (that corresponds to the parameters $\mbox{vec}(B_\ell)$) and for fixed $j$ and $k$ such that  $\beta_i \neq \beta_j \neq \beta_k$, we have
\begin{align}\label{Der-3-Beta}
\nonumber\dfrac{\partial^3\mathbb{H}_t(i_1)}{\partial \beta_i\partial\beta_j\partial\beta_k} &= \sum_{j_1=1}^m\sum\limits_{k'=3}^\infty \left[ \sum\limits_{l=3}^{k'} \left[ \sum\limits_{r=2}^{l-1}\left\{\left( \sum\limits_{a=1}^{r-1}\mathbb{B}^{a-1}_{i_1,j_1}\mathbb{B}^{(k)}_{i_1,j_1}\mathbb{B}^{r-1-a}_{i_1,j_1}\right)\mathbb{B}^{(j)}_{i_1,j_1}\mathbb{B}^{l-1-r}_{i_1,j_1}\right\}\mathbb{B}^{(i)}_{i_1,j_1}\mathbb{B}^{k'-l}_{i_1,j_1}\right.\right.\\\nonumber
&\left.\qquad + \sum\limits_{r=1}^{l-2}\left\{\mathbb{B}^{r-1}_{i_1,j_1}\mathbb{B}^{(j)}_{i_1,j_1}\left(\sum\limits_{a=1}^{l-1-r}\mathbb{B}^{a-1}_{i_1,j_1}\mathbb{B}^{(k)}_{i_1,j_1}\mathbb{B}^{l-1-r-a}_{i_1,j_1}\right)\right\}\mathbb{B}^{(i)}_{i_1,j_1}\mathbb{B}^{k'-l}_{i_1,j_1}\right] \\\nonumber
&\qquad + \sum\limits_{l=2}^{k'-1}\left\{\left(\sum\limits_{r=1}^{l-1}\mathbb{B}^{r-1}_{i_1,j_1}\mathbb{B}^{(j)}_{i_1,j_1}\mathbb{B}^{l-1-r}_{i_1,j_1}\right)\mathbb{B}^{(i)}_{i_1,j_1}\left(\sum\limits_{a=1}^{k'-l}\mathbb{B}^{a-1}_{i_1,j_1}\mathbb{B}^{(k)}_{i_1,j_1}\mathbb{B}^{k'-l-a}_{i_1,j_1}\right)\right\}\\\nonumber
& + \sum\limits_{l=3}^{k'-1}\left[\sum\limits_{r=2}^{l}\left\{\left(\sum\limits_{a=1}^{r-1}\mathbb{B}^{a-1}_{i_1,j_1}\mathbb{B}^{(k)}_{i_1,j_1}\mathbb{B}^{r-1-a}_{i_1,j_1}\right)\mathbb{B}^{(i)}_{i_1,j_1}\left(\sum\limits_{m=1}^{k'-r}\mathbb{B}^{m-1}_{i_1,j_1}\mathbb{B}^{(j)}_{i_1,j_1}\mathbb{B}^{k'-r-m}_{i_1,j_1}\right)\right\}\right]\\\nonumber
&\qquad + \sum\limits_{l=1}^{k'-2}\left[\mathbb{B}^{l-1}_{i_1,j_1}\mathbb{B}^{(i)}_{i_1,j_1}\sum\limits_{r=2}^{k'-l}\left\{\left(\sum\limits_{a=1}^{r-1}\mathbb{B}^{a-1}_{i_1,j_1}\mathbb{B}^{(k)}_{i_1,j_1}\mathbb{B}^{r-1-a}_{i_1,j_1}\right)\mathbb{B}^{(j)}_{i_1,j_1}\mathbb{B}^{k'-l-r}_{i_1,j_1}\right\}\right. \\
&\left.\left.\qquad + \mathbb{B}^{l-1}_{i_1,j_1}\mathbb{B}^{(i)}_{i_1,j_1}\sum\limits_{r=1}^{k'-l-1}\left\{\mathbb{B}^{r-1}_{i_1,j_1}\mathbb{B}^{(j)}_{i_1,j_1}\left(\sum\limits_{a=1}^{k'-l-r}\mathbb{B}^{a-1}_{i_1,j_1}\mathbb{B}^{(k)}_{i_1,j_1}\mathbb{B}^{k'-l-r-a}_{i_1,j_1}\right)\right\}\right]\right]\underline{c}_{t-k}(j_1),
\end{align}
where $\mathbb{B}_{i_1,j_1}=\mathbb{B}(i_1,j_1)$ denotes the $(i_1,j_1)-$th component   of the matrix $\mathbb{B}$.
By multiplying by $\beta_i$, $\beta_j$ and $\beta_k$ we obtain

\begin{align*}
\beta_i\beta_j\beta_k\dfrac{\partial^3\mathbb{H}_t(i_1)}{\partial \beta_i\partial\beta_j\partial\beta_k} &= \sum\limits_{j_1=1}^m\sum\limits_{k'=3}^\infty \left[ \sum\limits_{l=3}^{k'} \left[ \sum\limits_{r=2}^{l-1}\left\{\left( \sum\limits_{a=1}^{r-1}\mathbb{B}_{i_1,j_1}^{a-1}\beta_k\mathbb{B}_{i_1,j_1}^{(k)}\mathbb{B}_{i_1,j_1}^{r-1-a}\right)\beta_j\mathbb{B}_{i_1,j_1}^{(j)}\mathbb{B}_{i_1,j_1}^{l-1-r}\right\}\beta_i\mathbb{B}_{i_1,j_1}^{(i)}\mathbb{B}_{i_1,j_1}^{k'-l}\right.\right.\\\nonumber
&\left.\qquad + \sum\limits_{r=1}^{l-2}\left\{\mathbb{B}_{i_1,j_1}^{r-1}\beta_j\mathbb{B}_{i_1,j_1}^{(j)}\left(\sum\limits_{a=1}^{l-1-r}\mathbb{B}_{i_1,j_1}^{a-1}\beta_k\mathbb{B}_{i_1,j_1}^{(k)}\mathbb{B}_{i_1,j_1}^{l-1-r-a}\right)\right\}\beta_i\mathbb{B}_{i_1,j_1}^{(i)}\mathbb{B}_{i_1,j_1}^{k'-l}\right] \\\nonumber
&\qquad + \sum\limits_{l=2}^{k'-1}\left\{\left(\sum\limits_{r=1}^{l-1}\mathbb{B}_{i_1,j_1}^{r-1}\beta_j\mathbb{B}_{i_1,j_1}^{(j)}\mathbb{B}_{i_1,j_1}^{l-1-r}\right)\beta_i\mathbb{B}_{i_1,j_1}^{(i)}\left(\sum\limits_{a=1}^{k'-l}\mathbb{B}_{i_1,j_1}^{a-1}\beta_k\mathbb{B}_{i_1,j_1}^{(k)}\mathbb{B}_{i_1,j_1}^{k'-l-a}\right)\right\}\\\nonumber
& + \sum\limits_{l=3}^{k'-1}\left[\sum\limits_{r=2}^{l}\left\{\left(\sum\limits_{a=1}^{r-1}\mathbb{B}_{i_1,j_1}^{a-1}\beta_k\mathbb{B}_{i_1,j_1}^{(k)}\mathbb{B}_{i_1,j_1}^{r-1-a}\right)\beta_i\mathbb{B}_{i_1,j_1}^{(i)}\left(\sum\limits_{m=1}^{k'-r}\mathbb{B}_{i_1,j_1}^{m-1}\beta_j\mathbb{B}_{i_1,j_1}^{(j)}\mathbb{B}_{i_1,j_1}^{k'-r-m}\right)\right\}\right]\\\nonumber
&\qquad + \sum\limits_{l=1}^{k'-2}\left[\mathbb{B}_{i_1,j_1}^{l-1}\beta_i\mathbb{B}_{i_1,j_1}^{(i)}\sum\limits_{r=2}^{k'-l}\left\{\left(\sum\limits_{a=1}^{r-1}\mathbb{B}_{i_1,j_1}^{a-1}\beta_k\mathbb{B}_{i_1,j_1}^{(k)}\mathbb{B}_{i_1,j_1}^{r-1-a}\right)\beta_j\mathbb{B}_{i_1,j_1}^{(j)}\mathbb{B}_{i_1,j_1}^{k'-l-r}\right\}\right. \\
&\left.\left.\qquad + \mathbb{B}_{i_1,j_1}^{l-1}\beta_i\mathbb{B}_{i_1,j_1}^{(i)}\sum\limits_{r=1}^{k'-l-1}\left\{\mathbb{B}_{i_1,j_1}^{r-1}\beta_j\mathbb{B}_{i_1,j_1}^{(j)}\left(\sum\limits_{a=1}^{k'-l-r}\mathbb{B}_{i_1,j_1}^{a-1}\beta_k\mathbb{B}_{i_1,j_1}^{(k)}\mathbb{B}_{i_1,j_1}^{k'-l-r-a}\right)\right\}\right]\right]\underline{c}_{t-k}(j_1)\\
&\leq \sum\limits_{j_1=1}^m \sum\limits_{k'=3}^\infty \left[\sum\limits_{l=3}^{k'} (l-2)\mathbb{B}_{i_1,j_1}^{k'} + \sum\limits_{l=2}^{k'-1}(l-1)(k'-l)\mathbb{B}_{i_1,j_1}^{k'} + \sum\limits_{l=1}^{k'-2}(k'-l-1)\mathbb{B}_{i_1,j_1}^{k'}\right]\underline{c}_{t-k}(j_1)\\
&\leq \sum\limits_{j_1=1}^m\sum\limits_{k'=3}^\infty k'(k'-1)(k'-2)\mathbb{B}_{i_1,j_1}^{k'}\underline{c}_{t-k}(j_1).
\end{align*}
So we deduce that
\begin{align}\label{Maj-Der-3-Beta}
\mathbb E \left |\dfrac{\beta_i\beta_j\beta_k}{\mathbb{H}_t(i_1)} \dfrac{\partial^3\mathbb{H}_t(i_1)}{\partial \beta_i\partial\beta_j\partial\beta_k} \right |^{r_0} <\infty.
\end{align}
\end{enumerate}
The three above cases prove the existence of the moments of the third order derivatives. As before, our estimations are in fact uniform and we may write that on a neighborhood $V(\theta_0)$ of $\theta_0\in\overset{\circ}{\Theta}$, for all $i_1=1,\ldots, m$ and $i,j,k=1,\ldots, s_1$, we have
\begin{align}\label{unif3}
\mathbb{E}\sup\limits_{\theta\in V(\theta_0)}\left\vert \dfrac{1}{h_{i_1,t}^{\delta_{i_1}/2}}\dfrac{\partial^3h_{i_1,t}^{\delta_{i_1}/2}}{\partial \theta_i \partial\theta_j\partial\theta_k}(\theta)\right\vert^{r_0}<\infty.
\end{align}
These estimations imply that
\begin{equation}\label{unif3ellt}
\mathbb{E}\left[\sup\limits_{\theta\in V(\theta_0)}\left\vert \dfrac{\partial^3 l_t(\theta)}{\partial \theta_i \partial \theta_j \partial\theta_k}\right\vert \right]<\infty.
\end{equation}
By a Taylor expansion around $\theta_0$, for all $i$ and $j$ it holds that
\begin{equation}\label{Preuve-Conv-J-Connue}
\dfrac{1}{n}\sum\limits_{t=1}^n \dfrac{\partial^2 l_t}{\partial \theta_i \partial\theta_j} (\theta_{ij}^\ast) = \dfrac{1}{n}\sum\limits_{t=1}^n\dfrac{\partial^2 l_t}{\partial \theta_i\partial \theta_j} (\theta_0) + \dfrac{1}{n}\sum\limits_{t=1}^n\dfrac{\partial}{\partial \theta'}\left\{\dfrac{\partial^2 l_t}{\partial \theta_i\partial\theta_j}(\tilde{\theta}_{ij})\right\}(\theta_{ij}^\ast - \theta_0),
\end{equation}
where $\tilde{\theta}_{ij}$ lies between $\theta_{ij}^\ast$ et $\theta_0$. Using the almost sure convergence of $\tilde{\theta}_{ij}$ to $\theta_0$, the ergodic
theorem and \eqref{unif3ellt}, we imply that almost-surely
\begin{align*}
\underset{n\to +\infty}{\lim\sup}\left\Vert \dfrac{1}{n}\sum\limits_{t=1}^n\dfrac{\partial}{\partial \theta'}\left\{\dfrac{\partial^2l_t}{\partial\theta_i\partial\theta_j}(\tilde{\theta}_{ij})\right\}\right\Vert &\leq \underset{n\to +\infty}{\lim\sup} \dfrac{1}{n}\sum\limits_{t=1}^n \sup\limits_{\theta\in V(\theta_0)}\left\Vert \dfrac{\partial}{\partial \theta'}\left\{\dfrac{\partial^2l_t}{\partial\theta_i\partial\theta_j}(\theta_{ij})\right\}\right\Vert \\
& = \mathbb{E}_{\theta_0}\sup\limits_{\theta\in V(\theta_0)} \left\Vert \dfrac{\partial}{\partial \theta'}\left\{\dfrac{\partial^2l_t}{\partial\theta_i\partial\theta_j}(\theta_{ij})\right\}\right\Vert < \infty.
\end{align*}
Since $\Vert \theta_{ij}^\ast - \theta_0 \Vert \underset{n\to \infty}{\longrightarrow}0$, the second term in the right hand side of  \eqref{Preuve-Conv-J-Connue} converges to $0$ almost-surely. By the ergodic theorem and using the same arguments than in the proof of theorem 2.2 in \cite{FZ-bernoulli} it follows that
\begin{equation}\label{Preuve-J-Connue}
\dfrac{1}{n}\sum\limits_{t=1}^n \dfrac{\partial^2l_t}{\partial \theta_i\partial\theta_j}(\theta_{ij}^\ast) \longrightarrow J(i,j) \ \text{in probability}.
\end{equation}

\subsubsection{Invertibility of the matrix $J$}%
To prove the invertibility of the matrix $J$ we calculate the derivatives of the criterion $\{\partial l_t(\theta)\}\{\partial\theta_i\}$ and $\{\partial^2 l_t(\theta)\}\{\partial\theta_i\partial\theta_j\}$ as functions of $H_t$.
We start from
\[l_t(\theta) = \underline{\varepsilon}_t'H_t^{-1}\underline{\varepsilon}_t + \log(\det(H_t))\]
and we have the first derivative
\[\dfrac{\partial l_t(\theta)}{\partial \theta_i} = Tr\left[(H_t^{-1} - H_t^{-1}\underline{\varepsilon}_t\underline{\varepsilon}_t'H_t^{-1})\dfrac{\partial H_t}{\partial\theta_i}\right],\]
and the second derivative
\begin{equation*}
\begin{aligned}
\dfrac{\partial^2 l_t(\theta)}{\partial \theta_i\partial\theta_j} &= Tr \left[-H_t^{-1}\dfrac{\partial H_t}{\partial\theta_j}H_t^{-1}\dfrac{\partial H_t}{\partial\theta_i} + H_t^{-1}\dfrac{\partial H_t}{\partial\theta_i\partial\theta_j} + H_t^{-1}\dfrac{\partial H_t}{\partial\theta_j}H_t^{-1}\underline{\varepsilon}_t\underline{\varepsilon}_t'H_t^{-1}\dfrac{\partial H_t}{\partial \theta_i} \right.\\
&\qquad \left. + H_t^{-1}\underline{\varepsilon}_t\underline{\varepsilon}_t'H_t^{-1}\dfrac{\partial H_t}{\partial \theta_j}H_t^{-1}\dfrac{\partial H_t}{\partial\theta_i}  -H_t^{-1}\underline{\varepsilon}_t\underline{\varepsilon}_t'H_t^{-1}\dfrac{\partial^2H_t}{\partial\theta_i\partial\theta_j} \right].
\end{aligned}
\end{equation*}
Since
$$
J = \mathbb{E}\left[\dfrac{\partial^2 l_t(\theta_0)}{\partial \theta_i\partial\theta_j} \right ] = \mathbb E \left (
\mathbb{E}\left[\left.\dfrac{\partial^2 l_t(\theta_0) }{\partial \theta_i\partial\theta_j}\right\vert \mathcal{F}_{t-1}\right] \right  ) ,  $$
we compute the conditional expectation as follows (with the convention that $H_{0t}=H_t(\theta_0)$):
\begin{equation*}
\begin{aligned}
\mathbb{E}\left[\left.\dfrac{\partial^2 l_t}{\partial \theta_i\partial\theta_j}(\theta_0) \right\vert \mathcal{F}_{t-1}\right]    &= \mathbb{E}\left[ - Tr\left(H_{0t}^{-1} \dfrac{\partial H_{0t}}{\partial\theta_j}H_{0t}^{-1} \dfrac{\partial H_{0t}}{\partial\theta_i} \right) - Tr\left( H_{0t}^{-1} \underline{\varepsilon}_t\underline{\varepsilon}_t' H_{0t}^{-1} \dfrac{\partial^2 H_{0t}}{\partial\theta_i\partial\theta_j}\right) \right.\\
& \hspace{3cm} + Tr\left(H_{0t}^{-1} \dfrac{\partial H_{0t}}{\partial\theta_j}H_{0t}^{-1} \underline{\varepsilon}_t\underline{\varepsilon}_t' H_{0t} \dfrac{\partial H_{0t}}{\partial\theta_i} \right) + Tr\left( H_{0t}^{-1} \dfrac{\partial^2 H_{0t}}{\partial\theta_i\partial\theta_j}\right) \\
&\hspace{3cm} \left. \left.+ Tr\left( H_{0t}^{-1} \underline{\varepsilon}_t\underline{\varepsilon}_t'H_{0t}^{-1} \dfrac{\partial H_{0t}}{\partial\theta_i} H_{0t} \dfrac{\partial H_{0t}}{\partial\theta_j} \right) \right\vert \mathcal{F}_{t-1} \right]\\
&= Tr\left( H_{0t}^{-1} \dfrac{\partial^2 H_{0t}}{\partial\theta_i\partial\theta_j}\right) - Tr\left(H_{0t}^{-1} \dfrac{\partial H_{0t}}{\partial\theta_j}H_{0t}^{-1} \dfrac{\partial H_{0t}}{\partial\theta_i} \right)\\
&\hspace{3cm} + Tr\left(H_{0t}^{-1} \dfrac{\partial H_{0t}}{\partial\theta_j}H_{0t}^{-1}  \dfrac{\partial H_{0t}}{\partial\theta_i} \right)\mathbb{E}\left[\eta_t\eta_t' \big\vert \mathcal{F}_{t-1}\right]\\
&\hspace{3cm} + Tr\left(H_{0t}^{-1} \dfrac{\partial H_{0t}}{\partial\theta_i}H_{0t}^{-1}  \dfrac{\partial H_{0t}}{\partial\theta_j} \right)\mathbb{E}\left[\eta_t\eta_t' \big\vert \mathcal{F}_{t-1}\right] \\
&\hspace{3cm} -Tr\left( H_{0t}^{-1} \dfrac{\partial^2 H_{0t}}{\partial\theta_i\partial\theta_j}\right)\mathbb{E}\left[\eta_t\eta_t' \big\vert \mathcal{F}_{t-1}\right]\\
&= Tr\left(H_{0t}^{-1} \dfrac{\partial H_{0t}}{\partial\theta_j}H_{0t}^{-1} \dfrac{\partial H_{0t}}{\partial\theta_i} \right)
\end{aligned}
\end{equation*}
By the relation $Tr(A'B) = (\Vec A)'\Vec B$ we have
\begin{equation*}
Tr\left(H_{0t}H_{0t}^{(i)}H_{0t}H_{0t}^{(j)}\right) = \textbf{h}_i'\textbf{h}_j,
\end{equation*}
where $\textbf{h}_i = \Vec(H_{0t}^{-1/2}H_{0t}^{(i)}H_{0t}^{-1/2})$ and $\textbf{h}_j = \Vec(H_{0t}^{-1/2}H_{0t}^{(j)}H_{0t}^{-1/2})$. In view of $\Vec(ABC) = (C' \otimes A)\Vec B$ we have $\textbf{h}_i =((H_{0t}^{-1/2})'\otimes H_{0t}^{-1/2}) \textbf{d}_i$ with $\textbf{d}_i = \Vec(H_{0t}^{(i)})$.

We define the $m^2\times s_0$ matrices $\textbf{h} = (\textbf{h}_1 \vert \ldots \vert \textbf{h}_{s_0})$ and $\textbf{d} = (\textbf{d}_1 \vert \ldots \vert \textbf{d}_{s_0})$,
we have $\textbf{h} = \textbf{H}\textbf{d}$ with $\textbf{H} = (H_{0t}^{-1/2})' \otimes H_{0t}^{-1/2}$.\\
Reasoning by contradiction, we suppose that $J = \mathbb{E}\left[\textbf{h}'\textbf{h}\right]$ is singular. There exists a non-zero vector $\textbf{c} \in \mathbb{R}^{s_0}$, such that $\textbf{c}'J\textbf{c} = \mathbb{E}\left[\textbf{c}'\textbf{h}'\textbf{h}\textbf{c}\right] = 0$. Since $\textbf{c}'\textbf{h}'\textbf{h}\textbf{c} \geq 0$ almost surely, it means that $\textbf{c}'\textbf{h}'\textbf{h}\textbf{c} = \textbf{c}'\textbf{d}'\textbf{H}^2\textbf{d}\textbf{c} = 0$, almost surely.\\
The matrix $H_{0t}$ is definite-positive, then $H_{0t}^{-1/2}$ is too. This entails that $(H_{0t}^{-1/2})'\otimes H_{0t}^{-1/2}$ is definite positive. This implies that $\textbf{H}^2$ is a definite-positive matrix with probability 1, and consequently $\textbf{d}\textbf{c} = 0_{m^2}$ with probability 1.

We write $\textbf{c}$ as $\textbf{c} = (\textbf{c}_1',\textbf{c}_2')'$ where $\textbf{c}_1\in \mathbb{R}^{s_1}$ and $\textbf{c}_2\in \mathbb{R}^{s_4}$ where $s_4 = s_0 - s_1 = m(m-1)/2$ (which is the dimension of the parameters $\rho$). The rows $1,m+1,\ldots,m^2$ of the following equations
\begin{align}
\textbf{d}\textbf{c} &= \sum\limits_{i=1}^{s_0}c_i\dfrac{\partial}{\partial \theta_i}\text{vec}(H_{0t})= \sum\limits_{i=1}^{s_0}c_i\dfrac{\partial}{\partial \theta_i}\text{vec}(D_{0t}R_0D_{0t})= \sum\limits_{i=1}^{s_0}c_i\dfrac{\partial}{\partial \theta_i}\left[(D_{0t}\otimes D_{0t})\text{vec}(R_0)\right] \nonumber
\\
& = \sum\limits_{i=1}^{s_1}c_i\dfrac{\partial}{\partial \theta_i}\left[(D_{0t}\otimes D_{0t})\right]\text{vec}(R_0)+
\sum\limits_{i=s_1+1}^{s_0}c_i(D_{0t}\otimes D_{0t})\dfrac{\partial}{\partial \theta_i}\left[\text{vec}(R_0)\right] = 0_{m^2}\label{C.44}
\end{align}
yield
\begin{align}\label{C.45}
\sum\limits_{i=1}^{s_1}c_i\dfrac{\partial}{\partial \theta_i}\underline{h}_t^{\underline{\delta}/2}(\theta_0) = 0_m,\quad a.s.
\end{align}
by using \eqref{DD}. 
Differentiating the equation  \eqref{MAPGARCH}, we obtain that
\begin{align*}
\sum\limits_{k=1}^{s_1} c_k\dfrac{\partial}{\partial \theta_k}\underline{h}_t^{\underline{\delta}/2}(\theta_0) &= \sum\limits_{k=1}^{s_1}c_k\dfrac{\partial \underline{\omega}_0}{\partial \theta_k} + \sum\limits_{i=1}^q \sum\limits_{k=1}^{s_1}c_k\dfrac{\partial A_{0i}^+}{\partial \theta_k}(\underline{\varepsilon}_{t-i}^+)^{\underline{\delta}/2} + \sum\limits_{i=1}^q\sum\limits_{k=1}^{s_1}c_k\dfrac{\partial A_{0i}^-}{\partial \theta_k}(\underline{\varepsilon}_{t-i}^-)^{\underline{\delta}/2}\\
&\qquad+ \sum\limits_{j=1}^p\sum\limits_{k=1}^{s_1}c_k\dfrac{\partial B_{0j}}{\partial \theta_k} \underline{h}_{t-j}^{\underline{\delta}/2}(\theta_0) + \sum\limits_{j=1}^p\sum\limits_{k=1}^{s_1}c_kB_{0j}\dfrac{\partial \underline{h}_{t-j}^{\underline{\delta}/2}}{\partial \theta_k}(\theta_0)\\
& = \underline{\omega}_0^\ast + \sum\limits_{i=1}^qA_{0i}^{+,\ast}(\underline{\varepsilon}_{t-i}^+)^{\underline{\delta}/2} + \sum\limits_{i=1}^qA_{0i}^{-,\ast}(\underline{\varepsilon}_{t-i}^-)^{\underline{\delta}/2} + \sum\limits_{j=1}^pB_{0j}^\ast \underline{h}_{t-j}^{\underline{\delta}/2}\\&\qquad + \sum\limits_{j=1}^pB_{0j}\sum\limits_{k=1}^{s_1}c_k  \dfrac{\partial \underline{h}_{t-j}^{\underline{\delta}/2}}{\partial \theta_k}(\theta_0)
\end{align*}
where
\begin{equation*}
\underline{\omega}_0^\ast = \sum\limits_{k=1}^{s_1} c_k\dfrac{\partial\underline{\omega}_0}{\partial \theta_k}, \qquad A_{0i}^{+,\ast} = \sum\limits_{k=1}^{s_1} c_k \dfrac{\partial A_{0i}^+}{\partial \theta_k},\qquad A_{0i}^{-,\ast} = \sum\limits_{k=1}^{s_1} c_k \dfrac{\partial A_{0i}^-}{\partial \theta_k},\quad B_{0j}^\ast = \sum\limits_{k=1}^{s_1} c_k \dfrac{\partial B_{0j}}{\partial \theta_k}.
\end{equation*}
Because Equation \eqref{C.45} is satisfied for all $t$, we have
\begin{equation*}
\underline{\omega}_0^\ast + \sum\limits_{i=1}^q A_{0i}^{+,\ast}(\underline{\varepsilon}_{t-i}^+)^{\underline{\delta}/2} + \sum\limits_{i=1}^q A_{0i}^{-,\ast}(\underline{\varepsilon}_{t-i}^-)^{\underline{\delta}/2} + \sum\limits_{j=1}^pB_{0j}^\ast \underline{h}_{t-j}^{\underline{\delta}/2} = 0.
\end{equation*}
It follows that
\begin{align*}
\underline{h}_t^{\underline{\delta}/2}(\theta_0) &= (\underline{\omega}_0 - \underline{\omega}_0^\ast) + \sum\limits_{i=1}^q(A_{0i}^+ - A_{0i}^{+,\ast})(\underline{\varepsilon}_{t-i}^+)^{\underline{\delta}/2} + \sum\limits_{i=1}^q(A_{0i}^- - A_{0i}^{-,\ast})(\underline{\varepsilon}_{t-i}^-)^{\underline{\delta}/2}\\&\qquad + \sum\limits_{j=1}^p(B_{0j} - B_{0j}^\ast)\underline{h}_{t-j}^{\underline{\delta}/2}(\theta_0).
\end{align*}
Finally, we introduce the vector $\theta_1$ for which the  first $s_1$ components are
\begin{equation*}
\text{vec}(\underline{\omega}_0 - \underline{\omega}_0^\ast \vert A_{01}^+ - A_{01}^{+,\ast}\vert \ldots \vert A_{01}^- - A_{01}^{-,\ast} \vert \ldots \vert B_{01} - B_{01}^\ast\vert \ldots).
\end{equation*}
One may obtain  $\underline{h}_t^{\underline{\delta}/2}(\theta_0) = \underline{h}_t^{\underline{\delta}/2}(\theta_1)$ by choosing $\textbf{c}_1$ small enough in such a way that $\theta_1\in \Theta$. If  $\textbf{c}_1\neq 0$ then $\theta_1 \neq \theta_0$. This is in contradiction with the identifiability assumption and thus $\textbf{c}_1 = 0$. Consequently, Equation \eqref{C.44} becomes
\begin{equation*}
(D_{0t} \otimes D_{0t})\sum\limits_{k = s_1+1}^{s_0}c_k\dfrac{\partial}{\partial \theta_k}\text{vec}R_0 = 0_{m^2},\qquad a.s.
\end{equation*}
and then
\begin{equation*}
\sum_{i=s_1+1}^{s_0} c_k\dfrac{\partial}{\partial \theta_k} \text{vec}R_0 = 0_{m^2}.
\end{equation*}
Since the vectors $\partial \text{vec}R_0/\partial \theta_k$, $k = s_1+1,\ldots, s_0$ are linearly independent, the vector  $\textbf{c}_2 = (c_{s_1+1}, \ldots, c_{s_0})'$ is null and thus
$\textbf{c} = 0$. This is in contradiction with $\textbf{c}'\textbf{h}'\textbf{h}\textbf{c} = \textbf{c}'\textbf{d}'\textbf{H}^2\textbf{d}\textbf{c} = 0$ almost-surely. Therefore the assumption that  $J$ is not singular is absurd.

\subsubsection{Asymptotic irrelevance  of the initial values}\label{null-effect}
To conclude the proof, we have to deduce \eqref{NA1} and \eqref{ConvJ1} from \eqref{NA1} and \eqref{ConvJ1}. For this, we must show that the initial values have asymptotically no effect on the derivatives of the quasi likelihood. More precisely we may prove that
\begin{equation}\label{ecartGradient}
\left\Vert \dfrac{1}{\sqrt{n}}\sum\limits_{t=1}^n\left[\dfrac{\partial l_t(\theta_0)}{\partial\theta} - \dfrac{\partial \tilde{l}_t(\theta_0)}{\partial\theta'}\right] \right\Vert = \mathrm{o}_\mathbb P (1)
\end{equation}
and
\begin{equation}\label{ecartGradient2}
\dfrac{1}{n} \sum\limits_{t=1}^n\sup\limits_{\theta\in V(\theta_0)}\left\Vert \dfrac{\partial^2l_t(\theta)}{\partial\theta\partial\theta'} - \dfrac{\partial^2\tilde{l}_t(\theta)}{\partial\theta\partial\theta'} \right\Vert = \mathrm{o}_\mathbb P(1),
\end{equation}
for some neighbourhood $V(\theta_0)$.

The arguments are the same than in \cite{FZ-MAPGARCH}.

By \eqref{DL1}, it holds
\begin{equation*}
\dfrac{\partial l_t(\theta)}{\partial \theta_i} - \dfrac{\partial \tilde{l}_t(\theta)}{\partial \theta_i} = Tr(d_1 + d_2 + d_3),
\end{equation*}
with
\begin{align*}
d_1=&-D_t^{-1}\underline{\varepsilon}_t\underline{\varepsilon}_t'D_t^{-1}D_t(D_t^{-1} - \tilde{D}_t^{-1})R^{-1}D_t^{-1}D_t^{(i)}\\
d_2  =&  -D_t^{-1}\underline{\varepsilon}_t\underline{\varepsilon}_t'\tilde{D}_t^{-1}R^{-1}D_t^{-1}\left(D_t^{(i)} - \tilde{D}_t^{(i)}\right)\\
d_3  =&  -\left(D_t^{-1} - \tilde{D}_t^{-1}\right) \underline{\varepsilon}_t\underline{\varepsilon}_t'\tilde{D}_t^{-1}R^{-1}\tilde{D}_t^{-1} \tilde{D}_t^{(i)} - D_t^{-1}\underline{\varepsilon}_t\underline{\varepsilon}_t'\tilde{D}_t^{-1}R^{-1}\left(D_t^{-1} - \tilde{D}_t^{-1}\right)\tilde{D}_t^{(i)}\\
&  -D_t^{-1}R^{-1}D_t^{-1}D_t(D_t^{-1} - \tilde{D}_t^{-1})\underline{\varepsilon}_t\underline{\varepsilon}_t'D_t^{-1}D_t^{(i)} -D_t^{-1}R^{-1}\tilde{D}_t^{-1}\underline{\varepsilon}_t\underline{\varepsilon}_t'D_t^{-1}\left(D_t^{(i)} - \tilde{D}_t^{(i)}\right)\\
& -\left(D_t^{-1} - \tilde{D}_t^{-1}\right) R^{-1}\tilde{D}_t^{-1}\underline{\varepsilon}_t\underline{\varepsilon}_t'\tilde{D}_t^{-1} \tilde{D}_t^{(i)} - D_t^{-1}R^{-1}\tilde{D}_t^{-1}\underline{\varepsilon}_t\underline{\varepsilon}_t'\left(D_t^{-1} - \tilde{D}_t^{-1}\right)\tilde{D}_t^{(i)}\\
& + 2\left[D_t^{-1}\left(D_t^{(i)} - \tilde{D}_t^{(i)}\right) + \left(D_t^{-1} - \tilde{D}_t^{-1}\right)\tilde{D}_t^{(i)} \right].
\end{align*}
The term $d_3$ is a sum of term which can be handled as $d_1$ and $d_2$. Thus we need to prove that
$\sup_{\theta \in V(\theta_0)}\Vert D_t^{-1}\underline{\varepsilon}_t\Vert < \infty$,  $\sup_{\theta \in V(\theta_0)}\Vert D_t(D_t^{-1}-\tilde{D}_t^{-1})\Vert < \infty$, $\sup_{\theta \in V(\theta_0)}\Vert D_t^{-1}D_t^{(i)}\Vert < \infty$ and $\sup_{\theta \in V(\theta_0)}\Vert D_t^{-1}(D_t^{(i)} - \tilde{D}_t^{(i)})\Vert < \infty$.
From \eqref{MajH2.1}, \eqref{MajH-H1}, \eqref{MajDR1} and \eqref{MajH1.1}, we deduce that for any $t$
\begin{equation}\label{MajDinvD}
\begin{aligned}
\sup\limits_{\theta\in\Theta} \Vert D_t - \tilde{D}_t\Vert& \leq K\rho^t, \quad
\sup\limits_{\theta\in\Theta} \Vert D_t^{-1} \Vert \leq K , \quad
\sup\limits_{\theta\in\Theta} \Vert \tilde{D}_t^{-1} \Vert \leq K
, \\
\sup\limits_{\theta\in\Theta} \left|\frac{\tilde{h}_{i_1,t}^{1/2}(\theta)}{{h}_{i_1,t}^{1/2}(\theta)} \right|& \leq 1+K\rho^t \qquad\text{and }\quad
\sup\limits_{\theta\in\Theta} \left|\frac{{h}_{i_1,t}^{1/2}(\theta)}{\tilde{h}_{i_1,t}^{1/2}(\theta)} \right| \leq 1+K\rho^t,\, \text{ for } i_1=1,\dots,m.
\end{aligned}
\end{equation}
%
%
We remark that $D_t(D_t^{-1} - \tilde{D}_t^{-1}) = (\tilde{D}_t - D_t)\tilde{D}_t^{-1}$. Thus the above estimations yield
\begin{equation}\label{D.ecartD}
\sup\limits_{\theta\in\Theta}\Vert D_t(D_t^{-1} - \tilde{D}_t^{-1})\Vert= \sup\limits_{\theta\in\Theta}\Vert (\tilde{D}_t - D_t)\tilde{D}_t^{-1}\Vert \leq K\rho^t .
\end{equation}
%
By the matrix expressions \eqref{Hgrandinfini} and \eqref{ec2}, we have
\begin{equation*}
\mathbb{H}_t = \sum\limits_{k=0}^{t-r-1}\mathbb{B}^k\underline{c}_{t-k} + \mathbb{B}^{t-r}\mathbb{H}_r, \qquad \tilde{\mathbb{H}}_t = \sum\limits_{k=0}^{t-r-1}\mathbb{B}^k\tilde{\underline{c}}_{t-k} + \mathbb{B}^{t-r}\tilde{\mathbb{H}}_r,
\end{equation*}
where $r = \max\{p,q\}$. Since, $\underline{c}_t = \tilde{\underline{c}}_t$ for all $t>r$, we have
\begin{equation*}
\mathbb{H}_t - \tilde{\mathbb{H}}_t = \mathbb{B}^{t-r}(\mathbb{H}_r - \tilde{\mathbb{H}}_r),
\end{equation*}
and
\begin{equation*}
\dfrac{\partial}{\partial \theta_i} (\mathbb{H}_t - \tilde{\mathbb{H}}_t) = \mathbb{B}^{t-r} \dfrac{\partial}{\partial \theta_i} (\mathbb{H}_r - \tilde{\mathbb{H}}_r) + \sum\limits_{j=1}^{t-r} \mathbb{B}^{j-1}\mathbb{B}^{(i)}\mathbb{B}^{t-r-j}(\mathbb{H}_r - \tilde{\mathbb{H}}_r).
\end{equation*}
Since $\sup_{\theta\in\Theta}\rho(\mathbb{B}) < 1$ and \eqref{MajH1} or \eqref{Majh.delta}, we obtain
\begin{equation*}
\sup\limits_{\theta\in\Theta} \left\Vert \dfrac{\partial}{\partial \theta_i}(\mathbb{H}_t - \tilde{\mathbb{H}}_t)\right\Vert \leq K\rho^t
\end{equation*}
or equivalently
\begin{equation}\label{der.Majh.delta}
\sup\limits_{\theta\in\Theta} \left\Vert \dfrac{\partial}{\partial \theta_i}\left(\underline{h}_t^{\underline{\delta}/2}(\theta) - \tilde{\underline{h}}_t^{\underline{\delta}/2}(\theta)\right)\right\Vert \leq K\rho^t .
\end{equation}
%
By \eqref{DD} we have the expression of the derivative of $D_t$ (and analogously for the derivative of $\tilde{D}_t$). Thus we may write, for $i_1=1,...,m$ that
\begin{align*}
D_t^{-1}\dfrac{\partial}{\partial \theta_i}(D_t - \tilde{D}_t)(i_1,i_1) &= D_t^{-1} \big (  h_{i_1,t}^{1/2} - \tilde h_{i_1,t}^{1/2} \big )\dfrac{1}{\delta_{i_1}h_{i_1,t}^{\delta_{i_1}/2}}\dfrac{\partial h_{i_1,t}^{\delta_{i_1}/2}}{\partial\theta_i}(\theta)  \\
& + D_t^{-1} \left  (  h_{i_1,t}^{\delta_{i_1}/2} - \tilde h_{i_1,t}^{\delta_{i_1}/2} \right  )\dfrac{\tilde h_{i_1,t}^{1/2}}{\delta_{i_1} \, h_{i_1,t}^{\delta_{i_1}/2}\tilde h_{i_1,t}^{\delta_{i_1}/2}}\dfrac{\partial h_{i_1,t}^{\delta_{i_1}/2}}{\partial\theta_i}(\theta) \\
& + D_t^{-1} \dfrac{\tilde h_{i_1,t}^{1/2}}{\delta_{i_1} \tilde h_{i_1,t}^{\delta_{i_1}/2}}  \left ( \dfrac{\partial h_{i_1,t}^{\delta_{i_1}/2}}{\partial\theta_i}(\theta)-\dfrac{\partial \tilde h_{i_1,t}^{\delta_{i_1}/2}}{\partial\theta_i}(\theta)\right  ) \ .
\end{align*}
Using, for $i_1=1,...,m$, \eqref{Majh.delta}, \eqref{MajH2.1}, \eqref{MajDinvD}, \eqref{D.ecartD}, \eqref{der.Majh.delta} and in view of \eqref{unif1}, we obtain
\begin{equation}\label{MajAvecU_t}
\sup\limits_{\theta \in \Theta} \left\Vert D_t^{-1}\dfrac{\partial}{\partial \theta_i}(D_t - \tilde{D}_t)(i_1,i_1)\right\Vert \leq K\rho^t {u_t},
\end{equation}
where $u_t$ is a squared integrable variable.
Using \eqref{eta-tilde}, from \eqref{h0-h1} and \eqref{MajDinvD}, we deduce
\begin{align}\label{MajAvecV_t.1}
\sup\limits_{\theta\in V(\theta_0)} \Vert D_t^{-1}\underline{\varepsilon}_t \Vert  &= \sup\limits_{\theta\in V(\theta_0)}\Vert D_t^{-1}D_{0t}\tilde{\eta}_t\Vert \leq {v_t}\Vert \tilde{\eta}_t\Vert , \\
\sup\limits_{\theta\in V(\theta_0)} \Vert \tilde{D}_t^{-1}\underline{\varepsilon}_t \Vert  &= \sup\limits_{\theta\in V(\theta_0)}\Vert \tilde{D}_t^{-1}D_t\Vert\Vert D_t^{-1}\underline{\varepsilon}_t\Vert \leq (1 + K\rho^t){v_t}\Vert \tilde{\eta}_t\Vert,\label{MajAvecV_t.2}
\end{align}
where the random variable $v_t$ admits a fourth-order moment.
Now, using \eqref{MajDinvD}--\eqref{MajAvecV_t.2} and the
Cauchy-Schwarz inequality, we obtain
\begin{equation*}
\sup\limits_{\theta\in V(\theta_0)}\left\vert \dfrac{\partial l_t(\theta)}{\partial \theta_i} - \dfrac{\partial \tilde{l}_t(\theta)}{\partial \theta_i}\right\vert \leq K\rho^t{w_t},
\end{equation*}
where $w_t$ is an integrable variable. From the Markov inequality, we have
\begin{equation*}
\dfrac{1}{\sqrt{n}}\sum\limits_{t=1}^n \rho^tw_t =\mathrm{o}_{\mathbb{P}}(1),
\end{equation*}
which implies \eqref{ecartGradient}.
By exactly the same arguments, we obtain
\begin{equation*}
\sup\limits_{\theta\in V(\theta_0)}\left\vert \dfrac{\partial^2 l_t(\theta)}{\partial \theta_i\partial \theta_j} - \dfrac{\partial^2 \tilde{l}_t(\theta)}{\partial \theta_i\partial \theta_j}\right\vert \leq K\rho^t{w_t^\ast},
\end{equation*}
where $w_t^\ast$ is an integrable random variable.
Using the Borel-Cantelli lemma and the Markov inequality, we deduce that $\rho^t w_t^\ast$ goes to zero almost surely. Consequently, the Ces\'aro lemma implies that $n^{-1} \sum_{t=1}^n \rho^t w_t^\ast \to 0$ when $n$ goes to infinity, which entails \eqref{ecartGradient2}.

The proof of Theorem~\ref{AN-connu} is completed.\zak
\subsection{Proof of Theorem \ref{theorem31}}\label{A5}
In the sequel, we will use the version of the matrix representation \eqref{EcMat1} when the parameter $\delta$ is unknown. We write
\begin{equation}\label{EcMat2}
\mathbb{H}_t^{(\tau)} = \underline{c}_t^{(\tau)} + \mathbb{B}\cdot \mathbb{H}_{t-1}^{(\tau)},
\end{equation}
with
\[ \mathbb{H}_t^{(\tau)} = \begin{pmatrix} \underline{h}_t^{\underline{\tau}/2} \\ \underline{h}_{t-1}^{\underline{\tau}/2} \\ \vdots \\ \underline{h}_{t-p+1}^{\underline{\tau}/2} \end{pmatrix},\quad \underline{c}_t^{(\tau)} = \begin{pmatrix} \underline{\omega} + \sum\limits_{i=1}^q A_i^+(\underline{\varepsilon}_{t-i}^+)^{\underline{\tau}/2} + A_i^- (\underline{\varepsilon}_{t-1}^-)^{\underline{\tau}/2} \\ 0 \\ \vdots \\ 0 \end{pmatrix}\]
and
\[\mathbb{B} = \begin{pmatrix} B_1 & B_2 & & \ldots & & B_p\\ I_m & & & & & 0 \\ 0 & I_m & & & & \vdots\\ \vdots & & & \ddots & & \vdots \\ 0 & & & \ldots & I_m & 0 \end{pmatrix}, \]
and we can iterate the expression and we have
\begin{equation}\label{Hgrandinfinibis}
\mathbb{H}_t^{(\tau)} = \underline{c}_t^{(\tau)} + \mathbb{B}\underline{c}_{t-1}^{(\tau)} + \mathbb{B}^2\underline{c}_{t-2}^{(\tau)} + \ldots + \mathbb{B}^{t-1}\underline{c}_1^{(\tau)} + \mathbb{B}^t\mathbb{H}_0^{(\tau)} = \sum\limits_{k=0}^\infty \mathbb{B}^k \underline{c}_{t-k}^{(\tau)}.
\end{equation}

We prove our consistency statement when $\underline{\delta}_0$  is unknown. As in the case where $\underline{\delta}_0$ was known (see Section \ref{A2}), the proof is decomposed in the four following points which will be treated in separate subsections.
\begin{enumerate}[\footnotesize\bf\itshape\textit{\ref{A5}.}1.]
	\item Initial values do not influence quasi-likelihood: $\lim_{n\to\infty}\sup_{\nu \in \Delta} \vert \mathcal{L}_n(\nu) - \tilde{\mathcal{L}}_n(\nu) \vert  = 0$ a.s.
	\item Identifiability: If there exists $\in \mathbb{Z}$ such that $\underline{h}_t(\nu) = \underline{h}_t(\nu_0)$ almost surely and $R = R_0$, then $\nu = \nu_0$.
	\item Minimisation of the quasi log-likelihood on the true value: $\mathbb{E}_{\nu_0} \vert l_t(\nu_0) \vert < \infty,$ and if $\nu \neq \nu_0,\quad \mathbb{E}_{\nu_0}[l_t(\nu)] > \mathbb{E}_{\nu_0} [l_t(\nu_0)]$
	\item For any $\nu \neq \nu_0$ there exists a neighborhood $V(\nu)$ such that
\begin{align}\label{ivbis}	
\underset{n\to \infty}{\lim\inf}\inf\limits_{\nu^\ast \in V(\nu)} \tilde{\mathcal{L}}_n(\nu^\ast) > \mathbb{E}_{\nu_0}l_1(\nu_0),\quad a.s.
\end{align}
\end{enumerate}
There are many similarities with the proof of Theorem \ref{le-joli-label-de-yacouba}. We only indicates where the fact that the power is estimated has an importance is our reasoning.
%
%
\subsubsection{Initial values do not influence quasi-likelihood}
The proof is the same than the one done in  Subsection \ref{A21} when the power is assumed to be known.
\subsubsection{Identifiability}\label{new-proof1}
As regard to the proof of identifiability from Subsection \ref{A22}, it only remains to prove that if for $i_1=1,\dots,m$,  ${h}_{i_1,t}(\nu) / {h}_{i_1,t}(\nu_0) = 1$, a.s, then $\underline\tau = \underline{\delta}_0$. Let $\underline{\delta}_{0,i_1}$ (resp. $\underline\tau_{i_1}$) the $i_1^\text{th}$ element of $\underline{\delta}_0$ (resp. of $\underline\tau$).
We denote $Q^\pm(L)= \mathcal{B}(L)^{-1}\mathcal{A}^\pm(L)=\sum_{i\ge 1}Q_i^\pm L^i $ and $Q_0^\pm(L)= \mathcal{B}_0(L)^{-1}\mathcal{A}_0^\pm(L)=\sum_{i\ge 1}Q_{0i}^\pm L^i$.
Under  Assumption \textbf{A4}, by Proposition \ref{identif}, for any $i_2=1,\dots,m$, one may find $i_0\ge 1$ and $i_1\in\{1,\dots,m\}$ such that $Q_{0i_0}^++Q_{0i_0}^-\neq 0$. Since the coefficients of the matrix are positive, we denote by $(i_1,i_2)$ the position of a non zero element, for $i_2=1,\dots,m$.
By \eqref{upsilon} and \eqref{IC} we have
\begin{align*}
{h}_{i_1,t}^{\delta_{0,i_1}/2}(\nu_0) & = \sum_{j_1=1}^m\mathcal{B}_0(1)^{-1}(i_1,j_1)\underline{\omega}_0(j_1) +  \sum_{j_1=1}^m\sum_{l=1}^\infty  Q_{0l}^+(i_1,j_1) ({\varepsilon}_{j_1,t-l}^+)^{\delta_{0,j_1}} +\sum_{j_1=1}^m\sum_{l=1}^\infty  Q_{0l}^-(i_1,j_1) ({-\varepsilon}_{j_1,t-l}^-)^{\delta_{0,j_1}}  \\
 & = C_{i_1,t-i_0-1}(\nu_0) +     \sum_{j_1=1}^m Q_{0i_0}^+(i_1,j_1) ({\varepsilon}_{j_1,t-i_0}^+)^{\delta_{0,j_1}} +\sum_{j_1=1}^mQ_{0i_0}^-(i_1,j_1) (-{\varepsilon}_{j_1,t-i_0}^-)^{\delta_{0,j_1}}
 \end{align*}
where the quantities indexed by $t-i_0-1$ are $\mathcal F_{t-i_0-1}-$measurable. 
In the same way we have
\begin{align*}
{h}_{i_1,t}^{\tau_{i_1}/2}(\nu) & = \sum_{j_1=1}^m\mathcal{B}(1)^{-1}(i_1,j_1)\underline{\omega}(j_1) +  \sum_{j_1=1}^m\sum_{l=1}^\infty  Q_{l}^+(i_1,j_1) ({\varepsilon}_{j_1,t-l}^+)^{\tau_{j_1}} +\sum_{j_1=1}^m\sum_{l=1}^\infty  Q_{l}^-(i_1,j_1) (-{\varepsilon}_{j_1,t-l}^-)^{\tau_{j_1}}  
\\
 & = C_{i_1,t-i_0-1}(\nu) +    \sum_{j_1=1}^m Q_{i_0}^+(i_1,j_1) ({\varepsilon}_{j_1,t-i_0}^+)^{\tau_{j_1}} +\sum_{j_1=1}^mQ_{i_0}^-(i_1,j_1) (-{\varepsilon}_{j_1,t-i_0}^-)^{\tau_{j_1}} .
\end{align*}
Since ${h}_{i_1,t}(\nu) / {h}_{i_1,t}(\nu_0) = 1$, a.s, we have
\begin{align}\label{532}
& \frac{\left \{
C_{i_1,t-i_0-1}(\nu) +    \sum_{j_1=1}^m Q_{i_0}^+(i_1,j_1) ({\varepsilon}_{j_1,t-i_0}^+)^{\tau_{j_1}} +\sum_{j_1=1}^mQ_{i_0}^-(i_1,j_1) (-{\varepsilon}_{j_1,t-i_0}^-)^{\tau_{j_1}}
  \right \}^{\delta_{0,i_1}/\tau_{i_1}}}
  {
  C_{i_1,t-i_0-1}(\nu_0) +     \sum_{j_1=1}^m Q_{0i_0}^+(i_1,j_1) ({\varepsilon}_{j_1,t-i_0}^+)^{\delta_{0,j_1}} +\sum_{j_1=1}^mQ_{0i_0}^-(i_1,j_1) (-{\varepsilon}_{j_1,t-i_0}^-)^{\delta_{0,j_1}}
   } =1 \ \text{a.s.}
\end{align}
We denote $r_{j_1}= \tau_{j_1}/\delta_{0,j_1}$ and we introduce the function
$$ V(x_1,\dots,x_m) =
\frac{\left \{
 C_{i_1,t-i_0-1}(\nu) + \sum_{j_1=1}^m Q_{i_0}^+(i_1,j_1) (x_{j_1}^+)^{r_{j_1}\delta_{0,j_1}} +\sum_{j_1=1}^mQ_{i_0}^-(i_1,j_1) (-x_{j_1}^-)^{r_{j_1}\delta_{0,j_1}}
  \right \}^{1/r_{i_1}}}
  {
   C_{i_1,t-i_0-1}(\nu_0) + \sum_{j_1=1}^m Q_{0i_0}^+(i_1,j_1) (x_{j_1}^+)^{\delta_{0,j_1}} +\sum_{j_1=1}^mQ_{0i_0}^-(i_1,j_1) (-x_{j_1}^-)^{\delta_{0,j_1}}
   } .
$$
By \eqref{532}, $V(\eta_{t-1})=1$ almost-surely hence $V$ is almost surely constant on some neighborhood of zero (see Assumption \textbf{A8}). Hence
for any $y=x_{j_1}^{\delta_{0,j_1}}\in [a{,}b] \subset [0{,}+\infty[$:
$$ V(0,...,0,y,0,...,0) =
\frac{\left \{
 C_{i_1,t-i_0-1}(\nu) +    Q_{i_0}^+(i_1,i_2) (y)^{r_{i_2}}
  \right \}^{1/r_{i_1}}}
  {
   C_{i_1,t-i_0-1}(\nu_0) +      Q_{0i_0}^+(i_1,i_2) (y)
   } =1
$$
almost-surely.
Since the coefficients $C_{i_1,t-i_0-1}(\nu)$, $Q_{i_0}^+(i_1,i_2)$, $C_{i_1,t-i_0-1}(\nu_0)$ and $Q_{0i_0}^+(i_1,i_2)$  are positive, we deduce that $r_{i_2}=r_{i_1}:=r$ after differentiate twice the above equation.
Starting now from
$$ V(0,...,0,y,0,...,0) =
\frac{\left \{
 C_{i_1,t-i_0-1}(\nu) +    Q_{i_0}^+(i_1,i_2) (y)^{r}
  \right \}^{1/r}}
  {
   C_{i_1,t-i_0-1}(\nu_0) +      Q_{0i_0}^+(i_1,i_2) (y)
   } =1,
$$
we can deduce by differentiating twice again, as in \cite{HZ-APGARCH},  that $r=r_{i_1}=1$. Hence we have $\tau_{i_2}=\delta_{0,i_2}$. This is done for any $i_2=1,\dots,m$ so the result is proved.

\subsubsection{Minimisation of the likelihood on the true value}
Replacing $\theta_0$ by $\nu_0$, the proof is the same than the one 
when the power is assumed to be known.
\subsubsection{Proof of (\ref{ivbis})}
Once again, the proof is the same than the one 
when the power is assumed to be known.
\subsubsection{Conclusion}
The proof Theorem  \ref{theorem31} follows the argument from  Theorem  \ref{le-joli-label-de-yacouba}. 
\subsection{Proof of Theorem \ref{AN-inconnu-c-est-le-label-de-bruno}}\label{aninconnu}
Now we deal with the asymptotic normality result when $\underline{\delta}_0$ is unknown.
We follow the arguments and the different steps that we used in the proof of Theorem \ref{AN-connu} in Section \ref{anconnu}.
%
%
%
%
%
%
%
To establish the asymptotic normality result when the power is known, the proof is again decomposed in six intermediates points.
\begin{enumerate}[\footnotesize\bf\itshape\textit{\ref{aninconnu}.}1.]
	\item First derivative of the  quasi log-likelihood
	\item Existence of moments at any order of the score
	\item Asymptotic normality of the score vector:
	\begin{equation}\label{NAA}
\dfrac{1}{\sqrt{n}}\sum\limits_{t=1}^n \dfrac{\partial {l}_t(\nu_0)}{\partial\nu} \overset{\mathcal{L}}{\longrightarrow} \mathcal{N}(0,I).
\end{equation}
	\item Convergence to $J$:
	\begin{equation}\label{ConvJ2}
\dfrac1n \sum\limits_{t=1}^n \dfrac{\partial^2 {l}_t(\nu_{ij}^\ast)}{\partial \nu_i \partial\nu_j} \longrightarrow J(i,j) \mbox{ in probability,}
\end{equation}
	\item Invertibility of the matrix $J$
	\item Asymptotic irrelevance  of the initial values
\end{enumerate}
We introduce the following notations:
\begin{enumerate}[$\qquad \ast$]
	\item $s_0 = 2m + (p+2q)m^2 + m(m-1)/2$,
	\item $s_1 = 2m + (p+2q)m^2$,
	\item $s_2 = m + (p+2q)m^2$,
	\item $s_3 = m + 2qm^2$,
	\item $s_4 = m + qm^2$.
\end{enumerate}
\subsubsection{First derivative of the  quasi log-likelihood}
The aim of this subsection is to establish the expressions of the first order derivatives of the quasi log-likelihood. We may argue as in subsection \ref{A31}.

We denote $D_{0t} = D_t(\nu_0), R_0 = R(\nu_0)$,
\[D_{0t}^{(i)} = \dfrac{\partial D_t}{\partial \nu_i}(\nu_0), \quad R_0^{(i)} = \dfrac{\partial R}{\partial\nu_i}(\nu_0),\]
	\[D_{0t}^{(i,j)} = \dfrac{\partial^2D_t}{\partial\nu_i\partial\nu_j}(\nu_0),\quad R_0^{(i,j)} = \dfrac{\partial^2 R}{\partial\nu_i\partial\nu_j}(\nu_0),\]
and $\underline{\varepsilon}_t = D_{0t}\tilde{\eta}_t$, where $\tilde{\eta}_t(\nu) = R^{1/2}\eta_t(\nu)$ with   $\tilde{\eta}_t=\tilde{\eta}_t(\nu_0) = R_0^{1/2}\eta_t$.

When we  differentiate with respect to $\nu_i$ for $i = 1, \ldots, s_1$ (that is with respect to $(\underline{\omega}', {\alpha_{1}^{+}} ', \ldots, {\alpha_q^+}', {\alpha_1^-}', \ldots, {\alpha_q^-}', \beta'_{1}, \ldots, \beta'_{p}, \underline{\tau}')'$) we obtain:
\begin{align}\label{DL1nu}
\dfrac{\partial l_t(\nu)}{\partial \nu_i}  & = -Tr\left((\underline{\varepsilon}_t\underline{\varepsilon}_t'D_t^{-1}R^{-1} + R^{-1}D_t^{-1}\underline{\varepsilon}_t\underline{\varepsilon}_t') D_t^{-1}\dfrac{\partial D_t}{\partial \nu_i}D_t^{-1}\right) + 2Tr\left(D_t^{-1} \dfrac{\partial D_t}{\partial\nu_i}\right) \\
\dfrac{\partial l_t(\nu_0)}{\partial \nu_i}&  = Tr\left[ \left(I_m - R_0^{-1}\tilde{\eta}_t\tilde{\eta}_t'\right)D_{0t}^{(i)} D_{0t}^{-1} + \left(I_m - \tilde{\eta}_t\tilde{\eta}_t'R_0^{-1}\right)D_{0t}^{-1}D_{0t}^{(i)}\right]. \label{DL1nubis}
\end{align}
We differentiate with respect to $\nu_i$  for $i = s_1+1, \ldots, s_0$ (that is with respect to $\rho'$). We have
\begin{align}
\label{DL2nu}
\dfrac{\partial l_t(\nu)}{\partial \nu_i} &= -Tr\left( R^{-1} D_t^{-1} \underline{\varepsilon}_t\underline{\varepsilon}_t' D_t^{-1}R^{-1}\dfrac{\partial R}{\partial\nu_i} \right) + Tr\left(R^{-1} \dfrac{\partial R}{\partial\nu_i}\right) \\
\label{DL2nubis}
\dfrac{\partial l_t(\nu_0)}{\partial \nu_i}& = Tr\left[\left(I_m - R_0^{-1}\tilde{\eta}_t\tilde{\eta}_t'\right)R_0^{-1}R_0^{(i)}\right].
\end{align}
\subsubsection{Existence of moments at any order for the score}
Arguing as in the beginnig of Subsection \ref{A32} we have:
\begin{enumerate}[(i)]
\item for $i=1, \ldots, s_1$
\begin{equation*}
\begin{aligned}
 \left |\dfrac{\partial l_t(\nu_0)}{\partial \nu_i}\right | &
	& \leq K \left\Vert D_{0t}^{(i)} D_{0t}^{-1} \right\Vert ,
\end{aligned}
\end{equation*}
\item for $i=s_1+1,\ldots, s_0$
\begin{equation*}
\begin{aligned}
 \left |\dfrac{\partial l_t(\nu_0)}{\partial \nu_i}\right | &  \leq K ,
\end{aligned}
\end{equation*}
\item for $i,j = 1,\ldots, s_1$
\begin{equation*}
\begin{aligned}
\mathbb{E}\left\vert\dfrac{\partial l_t(\nu_0)}{\partial\nu_i} \dfrac{\partial l_t(\nu_0)}{\partial\nu_j}\right\vert  &\leq K\left(\mathbb{E}\Vert D_{0t}^{(i)}D_{0t}^{-1} \Vert^2\mathbb{E}\Vert D_{0t}^{(j)}D_{0t}^{-1} \Vert^2\right)^{1/2},
\end{aligned}
\end{equation*}
\item for $i= 1, \ldots, s_1$ and $j = s_1+1,\ldots, s_0$
\begin{equation*}
\begin{aligned}
\mathbb{E} \left\vert\dfrac{\partial l_t(\nu_0)}{\partial\nu_i} \dfrac{\partial l_t(\nu_0)}{\partial\nu_j}\right\vert  &\leq K \mathbb{E}\left\Vert D_{0t}^{(i)} D_{0t}^{-1} \right\Vert,
\end{aligned}
\end{equation*}
\item and finally for $i,j = s_1+1, \ldots, s_0$, we have
\begin{equation*}
\begin{aligned}
\mathbb{E}\left\vert\dfrac{\partial l_t(\nu_0)}{\partial\nu_i} \dfrac{\partial l_t(\nu_0)}{\partial\nu_j}\right\vert
	&\leq K.
\end{aligned}
\end{equation*}
\end{enumerate}
To have the finiteness of the moments of the first derivative of the quasi log-likelihood, it remains to treat the cases   (i), (iii) and (iv) above. Thus, we  have to control the term $\Vert D_{0t}^{(i)}D_{0t}^{-1}\Vert$.
Since
\[ D_{0t} = \mathrm{Diag}(\underline{h}_t^{1/2}(\nu_0)) = \mathrm{Diag}(\underline{h}_t^{\underline{\tau}/2}(\nu_0))^{1/\underline{\tau}}, \]
we can work component wise.

All the computations that we have done in Subsection \ref{A32} are valid. This means that we have the same estimations on the derivatives as long as we differentiate with respect to $\nu_i$ for $i\in \{1,\dots,s_2\}$ (that is when we do not  differentiate with respect to $\underline\tau_j$ for $j=1,\dots,m$).
Indeed, for $i_1=1,\dots,m$ and $i=1,\dots,s_2$, we have
\begin{align}\label{DDbis}
\dfrac{\partial D_{0t}(i_1,i_1)}{\partial \nu_{i}} & = \dfrac{1}{\tau_{i_1}} h_{i_1,t}^{1/2} \times \dfrac{1}{h_{i_1,t}^{\tau_{i_1}/2}} \dfrac{\partial h_{i_1,t}^{\tau_{i_1}/2}}{\partial\nu_i}(\nu_0)
\end{align}
and the reasonings are unchanged.

So we can focus ourselves on the derivatives with respect to $\underline{\tau}$:
\begin{align}\label{Dtau}
\dfrac{\partial D_t(i_1, i_1)}{\partial\tau_{j}}  & = h_{i_1,t}^{1/2}\left[ -\delta_{j,i_1} \dfrac{1}{\tau_{i_1}^2} \log\left(h_{i_1,t}^{\tau_{i_1}/2}\right) +\dfrac{1}{\tau_{i_1}} \dfrac{1}{h_{i_1,t}^{\tau_{i_1}/2}} \dfrac{\partial h_ {i_1,t}^{\tau_{i_1}/2}}{\partial \tau_{j}}(\nu_0)\right],
\end{align}
where $\delta_{j,i_1}$ denotes the Kronecker symbol.
Using the matrix expression \eqref{Hgrandinfinibis}, we calculate the derivatives  $\partial \mathbb{H}_t^{(\tau)} / \partial \nu_i$ for $i = s_2+1,\dots,s_2+m$ (with $s_1=s_2+m$)
and for $i_1=1,\dots,m$:
\[\dfrac{\partial\mathbb{H}_t^{(\tau)}(i_1)}{\partial\nu_i} = \sum\limits_{k = 0}^{\infty}\sum_{j_1=1}^m \mathbb{B}^k(i_1,j_1)\dfrac{\partial \underline{c}_{t-k}^{(\tau)}(j_1)}{\partial\nu_i},\]
with
$$\dfrac{\partial \underline{c}^{(\tau)}_{t-k}(j_1)}{\partial\nu_i} = \sum\limits_{l=1}^q\sum\limits_{j_2=1}^m \left \{  A_{l}^+(j_1,j_2)\dfrac{\partial({\varepsilon}_{j_2,t-l}^+)^{\tau_{j_2}}}{\partial\nu_i}
+  A_{l}^-(j_1,j_2)\dfrac{\partial(-{\varepsilon}_{j_2,t-l}^-)^{\tau_{j_2}}}{\partial\nu_i}\right \}. $$
Differentiating with respect to $\nu_i$ corresponds to a differentiation with respect to $\tau_{i_0}$ for an index $i_0\in\{1,\dots,m\}$. In the following computations, it is easy to work with an arbitrary order of derivation. So we write, for an order of derivation $\kappa\in\{1,2,3\}$, that
$$\dfrac{\partial^\kappa(\pm{\varepsilon}_{j_2,t-l}^\pm)^{\tau_{j_2}}}{\partial\nu^\kappa_i} =\dfrac{\partial^\kappa(\pm{\varepsilon}_{j_2,t-l}^\pm)^{\tau_{j_2}}}{\partial\tau_{i_0}^\kappa} = \left \{
\begin{array}{l}
 0 \ \text{if $j_2\neq i_0$} \\
\  \\
\log^\kappa(\pm{\varepsilon}_{i_0,t-l}^\pm) (\pm{\varepsilon}_{i_0,t-l}^\pm)^{\tau_{i_0}}  \ \text{if $j_2= i_0$} ,
\end{array}
\right .
$$
and we have
\begin{equation*}
\begin{aligned}
\dfrac{\partial^\kappa \underline{c}_{t-\kappa}^{(\tau)}(j_1)}{\partial\tau_{i_0}^\kappa} & =  \sum\limits_{l=1}^q A_l^+(j_1,i_0) \log^\kappa({\varepsilon}_{i_0,t-l}^+) ({\varepsilon}_{i_0,t-l}^+)^{\tau_{i_0}} + A_l^-(j_1,i_0)\log^\kappa(-{\varepsilon}_{i_0,t-l}^-) (-{\varepsilon}_{i_0,t-l}^-)^{\tau_{i_0}}.
\end{aligned}
\end{equation*}
By convention, we consider $\log(x^+) = 0$ when $x$ is negative and $\log(x^-) = 0$ when $x$ is positive.
\begin{equation*}
\begin{aligned}
\left\vert \dfrac{1}{\underline{c}^{(\tau)}_t(j_1)} \dfrac{\partial^\kappa \underline{c}_t^{(\tau)}(j_1)}{\partial \tau_{i_0}^\kappa} \right\vert &\leq  \dfrac{
\sum\limits_{l=1}^q A_l^+(j_1,i_0) |\log({\varepsilon}_{i_0,t-l}^+)|^\kappa ({\varepsilon}_{i_0,t-l}^+)^{\tau_{i_0}} + A_l^-(j_1,i_0)|\log(-{\varepsilon}_{i_0,t-l}^-)|^\kappa (-{\varepsilon}_{i_0,t-l}^-)^{\tau_{i_0}}
}
{\underline{\omega}(j_1)
+\sum\limits_{l=1}^q\sum\limits_{j_2=1}^m\left \{  A_{l}^+(j_1,j_2)({\varepsilon}_{j_2,t-l}^+)^{\tau_{j_2}}
+  A_{l}^-(j_1,j_2)(-{\varepsilon}_{j_2,t-l}^-)^{\tau_{j_2}}\right \}
}\\
& \le 	\dfrac{
\sum\limits_{l=1}^q A_l^+(j_1,i_0) |\log({\varepsilon}_{i_0,t-l}^+) |^\kappa({\varepsilon}_{i_0,t-l}^+)^{\tau_{i_0}} + A_l^-(j_1,i_0)|\log(-{\varepsilon}_{i_0,t-l}^-)|^\kappa (-{\varepsilon}_{i_0,t-l}^-)^{\tau_{i_0}}
}
{\underline{\omega}(j_1)
+\sum\limits_{l=1}^q\left \{  A_{l}^+(j_1,i_0)({\varepsilon}_{i_0,t-l}^+)^{\tau_{i_0}}
+  A_{l}^-(j_1,i_0)(-{\varepsilon}_{i_0,t-l}^-)^{\tau_{i_0}}\right \}
}\\	
	& \leq \sum\limits_{l=1}^q\left\vert \log({\varepsilon}_{i_0,t-l}^+) \right\vert^\kappa + \left\vert \log(-{\varepsilon}_{i_0,t-l}^-) \right\vert^\kappa .
\end{aligned}
\end{equation*}
Using the inequality
\[\mathbb{H}_t^{(\tau)}(i_1) \geq \omega + \sum\limits_{k'=1}^\infty\sum_{j_1=1}^m \mathbb{B}^{k'}(i_1,j_1)\underline{c}_{t-k'}^{(\tau)}(j_1)\ge \omega + \mathbb{B}^k(i_1,j_1)\underline{c}_{t-k}^{(\tau)},(j_1)\]
valid for any $k\ge 0$ where $\omega = \inf_{1\leq i\leq m} \underline{\omega}(i)$, and the fact that $x/(1+x)\le x^s$ for all $x\ge 0$, we obtain
\begin{equation*}
\begin{aligned}
 \sup\limits_{\nu \in \Delta} &\left\vert \dfrac{1}{\mathbb{H}_t^{(\tau)}(i_1)} \dfrac{\partial^\kappa \mathbb{H}_t^{(\tau)}(i_1)}{\partial \nu_i^\kappa} \right\vert \leq \sup\limits_{\nu \in \Delta} \sum\limits_{k=0}^\infty \sum\limits_{j_1=1}^m \dfrac{\mathbb{B}^k(i_1,j_1)}{\omega + \mathbb{B}^k(i_1, j_1)\underline{c}_{t-k}^{(\tau)}(j_1)} \left\vert \dfrac{\partial^\kappa \underline{c}_{t-k}^{(\tau)}(j_1)}{\partial \nu_i^\kappa} \right\vert \\
	&\leq \sup\limits_{\nu \in \Delta} \sum\limits_{k=0}^\infty \sum\limits_{j_1=1}^m \dfrac{\mathbb{B}^k(i_1,j_1)\underline{c}_{t-k}^{(\tau)}(j_1)}{\omega + \mathbb{B}^k(i_1, j_1)\underline{c}_{t-k}^{(\tau)}(j_1)} \left\vert \dfrac{1}{\underline{c}_{t-k}^{(\tau)}(j_1)}\dfrac{\partial^\kappa \underline{c}_{t-k}^{(\tau)}(j_1)}{\partial \nu_i^\kappa} \right\vert \\
	& \leq  \sup\limits_{\nu \in \Delta} \sum\limits_{k=0}^\infty \sum\limits_{j_1=1}^m \dfrac{\mathbb{B}^k(i_1,j_1)\underline{c}_{t-k}^{(\tau)}(j_1)}{\omega + \mathbb{B}^k(i_1, j_1)\underline{c}_{t-k}^{(\tau)}(j_1)} \left( \sum\limits_{j=1}^q \left\{ \left\vert \log({\varepsilon}_{i,t-k-j}^+)\right\vert^\kappa + \left\vert \log(-{\varepsilon}_{i,t-k-j}^-)\right\vert^\kappa \right\} \right)   \\
	& \leq   \sum\limits_{k=0}^\infty \sum\limits_{j_1=1}^m \sup\limits_{\nu \in \Delta} \left(\dfrac{\mathbb{B}^k(i_1,j_1)\underline{c}_{t-k}^{(\tau)}(j_1)}{\omega}\right)^s\left(\sum\limits_{j=1}^q  \left\{ \left\vert \log({\varepsilon}_{i,t-k-j}^+)\right\vert^\kappa + \left\vert \log(-{\varepsilon}_{i,t-k-j}^-)\right\vert^\kappa \right\}\right)  \\
	& \leq K \sum\limits_{k=0}^\infty \sum\limits_{j_1=1}^m \rho^{sk} \sup\limits_{\nu \in \Delta} \left(\underline{c}_{t-k}^{(\tau)}(j_1)\right)^s \left(\sum\limits_{j=1}^q\left\{ \left\vert \log({\varepsilon}_{i,t-k-j}^+)\right\vert^\kappa + \left\vert \log(-{\varepsilon}_{i,t-k-j}^-\right\vert^\kappa\right\} \right).
\end{aligned}
\end{equation*}
We have $\mathbb{E}_{\nu_0} \big (\sup_{\nu \in \Delta}(\underline{c}_{t-k}^{(\tau)}(j_1) )^{2s}\big )< \infty$, for all $s>0$ (see Corollary \ref{cor2}). So we obtain
\begin{align*}
\mathbb{E}_{\nu_0}\left[ \sup\limits_{\nu \in \Delta} \left\vert \dfrac{1}{\mathbb{H}_t^{(\tau)}(i_1)}\dfrac{\partial^\kappa \mathbb{H}_t^{(\tau)}(i_1)}{\partial\nu_i^\kappa} \right\vert \right] &\leq K \sum_{j=1}^q  (S^+_{i,j} + S^-_{i,j})
 \end{align*}
with for $j\ge 1$:
$$ S^\pm_{i,j}  = \sum\limits_{k=0}^\infty \rho^{sk} \left( \mathbb{E}_{\nu_0}\vert \log(\pm{\varepsilon}_{i,t-k-j}^\pm)\vert^{2\kappa}\right ) ^{1/2}  . $$
By stationarity, we treat only the terms $S^\pm_{i,1}$ and the computations are identical when one replaces $\underline{\varepsilon}^+$ by $\underline{\varepsilon}^-$ so we will only to treat $S^+_{i,1}$.
We have for any $A>0$
\begin{align*}
\mathbb{E}_{\nu_0}\vert \log({\varepsilon}_{i,t-k-1}^+)\vert^{2\kappa} & \le A +  \mathbb{E}_{\nu_0}\left \vert \log({\varepsilon}_{i,t-k-1}^+)\1_{ \vert \log({\varepsilon}_{i,t-k-1}^+)\vert^{2\kappa}  \ge A}\right\vert^{2\kappa} .
\end{align*}
It follows that
\begin{align}\label{decomp}
S_{i,1}^+ &\le  A + \sum\limits_{k=0}^\infty \rho^{sk}\left( \int_A^\infty \mathbb{P}\left(\left\vert \log({\varepsilon}_{i,t-k-1}^+)\right\vert^{2\kappa} > x\right)dx \right)^{1/2}\nonumber \\
	& \le  A+\sum\limits_{k=0}^\infty \rho^{sk}\left( \int_A^\infty \mathbb{P}\left(\left\vert \log({\varepsilon}_{i,t-k-1}^+) \right\vert > x^{1/2\kappa}\right)dx \right)^{1/2} \nonumber\\
	& \le  A+ \sum\limits_{k=0}^\infty\rho^{sk}\left(\int_A^\infty 2\kappa x^{2\kappa-1} \mathbb{P}\left(\left\vert \log({\varepsilon}_{i,t-k-1}^+(j_1))\right\vert > x\right)dx\right)^{1/2} \nonumber\\
	& \le  A+\sum\limits_{k=0}^\infty\rho^{sk} \left(\int_A^\infty 2\kappa x^{2\kappa-1}\left[ \mathbb{P}\left( \log({\varepsilon}_{i,t-k-1}^+) > x \right) + \mathbb{P}\left(\log({\varepsilon}_{i,t-k-1}^+) < -x\right)\right]dx\right)^{1/2} \nonumber\\
	&\le A+ \sum\limits_{k=0}^\infty\rho^{sk}\left( \int_A^\infty 2\kappa x^{2\kappa-1}\left[\mathbb{P}\left({\varepsilon}_{i,t-k-1}^+ \geq \exp\{x \}\right) + \mathbb{P}\left({\varepsilon}_{i,t-k-1}^+ < \exp\{-x\}\right)\right]dx\right)^{1/2} & \nonumber\\
	& \le  A+ \sum_{k=0}^\infty \rho^{sk} \left ( T_1^+ + T_2^+ \right ) ^{1/2} \ ,
\end{align}
with obvious notations. Using the Markov inequality one has
\begin{align}\label{s1}
T_1^+ & = \int_A^\infty 2\kappa x^{2\kappa-1}\mathbb{P}\left({\varepsilon}_{i,t-k-1}^+ \ge \exp\{x\}\right) dx  \le C \ .
\end{align}
The second term $T_2^+ $ is more difficult. One has to use the property
\begin{align}\label{hypA1}
\lim_{y\to 0^+} \frac{1}y\mathbb P({\varepsilon}_{i,t-k-1}^+ \le y ) &= 0 .
\end{align}
This is the Assumption A1 in \cite{pan} and in our case, it is a consequence of the fact that $\eta_t$ has a positive density on some neighborhood of zero (see Assumption \textbf{A8}).
We apply this property to $\mathbb{P}({\varepsilon}_{i,t-k-1}^+ < \exp\{-x\})$. Indeed one may find some $A>0$ such that $\exp\{-x\}$ is small. Thus for any $x\ge A$:
\begin{align}\label{pan1}
\mathbb{P}\left({\varepsilon}_{i,t-k-1}^+ < \exp\{-x\}\right) & \le  c \exp\{-x\} \ .
\end{align}
Hence, the second term $T_2^+$ satisfies
\begin{align}\label{s2}
T_2^+ & = \int_A^\infty 2\kappa x^{2\kappa-1}\mathbb{P}\left({\varepsilon}_{j_1,t-k}^+ < \exp\{-x \}\right) dx \le   c \int_A^\infty  2\kappa x^{2\kappa-1} \exp\{-x \}dx \le C<\infty \ .
\end{align}
We use \eqref{s1} and \eqref{s2} in \eqref{decomp} and we obtain that
$$ S_{i,1}^+  \le A +C  \sum_{k=0}^\infty  \rho^{sk}  <\infty \ .$$
This yields
\begin{align*}
\mathbb{E}_{\nu_0}\left[ \sup\limits_{\nu \in \Delta} \left\vert \dfrac{1}{\mathbb{H}_t^{(\tau)}(i_1)}\dfrac{\partial^\kappa \mathbb{H}_t^{(\tau)}(i_1)}{\partial\nu_i^\kappa} \right\vert \right] & < \infty
 \end{align*}
and it can be similarly be shown that, for any $r\ge 1$:
\begin{align}\label{esti-dtau}
\mathbb{E}_{\nu_0}\left[ \sup\limits_{\nu \in \Delta} \left\vert \dfrac{1}{\mathbb{H}_t^{(\tau)}(i_1)}\dfrac{\partial^\kappa \mathbb{H}_t^{(\tau)}(i_1)}{\partial\nu_i^\kappa} \right\vert^r \right] & < \infty
 \end{align}
or equivalently that for any $i_1=1,...,m$
 \begin{equation}\label{esti-dtau_bis}
 \mathbb{E}\sup\limits_{\nu\in V(\nu_0)}\left\vert \dfrac{1}{h_{i_1,t}^{\tau_{i_1}/2}}\dfrac{\partial^\kappa h_{i_1,t}^{\tau_{i_1}/2}}{\partial \nu_i^\kappa}(\nu)\right\vert^{r}<\infty.
 \end{equation}

\subsubsection{Asymptotic normality of the score vector}
 The arguments are the same than in the case of known power (see subsection \ref{A33}).
\subsubsection{Convergence to $J$}
\subsubsection*{$\rightsquigarrow$ Expression of the second order derivatives of the log-likelihood}
\addcontentsline{toc}{paragraph}{$\rightsquigarrow$\ Expression of the second order derivatives of the log-likelihood}
First one may remark that the algebraic expressions of the second order derivatives \eqref{D2-1}, \eqref{D2-2} and  \eqref{D2-3} are unchanged (even if the values of $s_0$ and $s_1$ take into account the parameter $\underline{\tau}$).
To control the derivates at the second order it is sufficient to control the term $\Vert D^{-1} D_{0t}^{(i,j)} \Vert$. For this, it necessary to calculate the second order derivatives of $D_{0t}$.
\subsubsection*{$\rightsquigarrow$ Existence of the moments of the second order derivatives of the log-likelihood}
\addcontentsline{toc}{paragraph}{$\rightsquigarrow$ Existence of the moments of the second order derivatives of the log-likelihood}
We give now the expressions of the second order derivatives of  $D_{0t}$.
\begin{enumerate}[(i)]
	\item For $i,j=1,\ldots, s_2$, we have
\begin{align*}
	 \dfrac{\partial^2 D_{0t}(i_1,j_1)}{\partial\nu_i\partial\nu_j} &= \dfrac{1}{\tau_{i_1}}h_{i_1,t}^{1/2}\left[ \dfrac{1}{h_{i_1,t}^{\tau_{i_1}/2}}\dfrac{\partial h_{i_1,t}^{\tau_{i_1}/2}}{\partial \nu_j} \dfrac{1}{h_{i_1,t}^{\tau_{i_1}/2}} \dfrac{\partial h_{i_1,t}^{\tau_{i_1}/2}}{\partial \nu_i} \left(\dfrac{1}{\tau_{i_1}} - 1 \right) + \dfrac{1}{h_{i_1,t}^{\tau_{i_1}/2}} \dfrac{\partial^2 h_{i_1,t}^{\tau_{i_1}/2}}{\partial \nu_i \partial \nu_j} \right]\\	
\end{align*}
	\item for $i_1 = 1,\ldots, m$ and $j=1,\ldots, s_2$, we obtain
\small{
\begin{align*}
\dfrac{\partial^2 D_{0t}(i_1,j_1)}{\partial\tau_{i_1}\partial\nu_j} &=  \dfrac{1}{\tau_{i_1}}h_{i_1,t}^{1/2}\left[ - \dfrac{1}{h_{i_1,t}^{\tau_{i_1}/2}}\dfrac{\partial h_{i_1,t}^{\tau_{i_1}/2}}{\partial\nu_j} \right.  + \dfrac{1}{\tau_{i_1}} \left(\dfrac{1}{\tau_{i_1}}\log\left(h_{i_1,t}^{\tau_{i_1}/2}\right) + \dfrac{1}{h_{i_1,t}^{\tau_{i_1}/2}}\dfrac{\partial h_{i_1,t}^{\tau_{i_1}/2}}{\partial \tau_{i_1}} \right)\dfrac{1}{h_{i_1,t}^{\tau_{i_1}/2}}\dfrac{\partial h_{i_1,t}^{\tau_{i_1}/2}}{\partial \nu_j} \\
	&\qquad \left. - \dfrac{1}{h_{i_1,t}^{\tau_{i_1}/2}}\dfrac{\partial h_{i_1,t}^{\tau_{i_1}/2}}{\partial \tau_{i_1}} \dfrac{1}{h_{i_1,t}^{\tau_{i_1}/2}}\dfrac{\partial h_{i_1,t}^{\tau_{i_1}/2}}{\partial \nu_i}  + \dfrac{1}{h_{i_1,t}^{\tau_{i_1}/2}}\dfrac{\partial^2 h_{i_1,t}^{\tau_{i_1}/2}}{\partial\tau_{i_1}\partial \nu_j}\right].
\end{align*}
}
\normalsize
	\item for $i_1=1,...,m$, $j = 1,\ldots, s_2$ and $i_0=1,\ldots, m$ such that $i_0\neq i_1$, we have
\begin{align*}
	 \dfrac{\partial^2 D_{0t}(i_1,j_1)}{\partial\tau_{i_0}\partial \nu_j}  &= \dfrac{1}{\tau_{i_1}}h_{i_1,t}^{1/2}\left[ \dfrac{1}{\tau_{i_0}}\left(-\dfrac{1}{\tau_{i_0}}\log\left(h_{i_1,t}^{\tau_{i_1}/2}\right) + \dfrac{1}{h_{i_1,t}^{\tau_{i_1}/2}}\dfrac{\partial h_{i_1,t}^{\tau_{i_1}/2}}{\partial \tau_{i_0}}\right) \dfrac{1}{h_{i_1,t}^{\tau_{i_1}/2}}\dfrac{\partial h_{i_1,t}^{\tau_{i_1}/2}}{\partial \nu_j}\right.\\
	 &\qquad\left. -  \dfrac{1}{h_{i_1,t}^{\tau_{i_1}/2}}\dfrac{\partial h_{i_1,t}^{\tau_{i_1}/2}}{\partial \tau_{i_0}}\dfrac{1}{h_{i_1,t}^{\tau_{i_1}/2}}\dfrac{\partial h_{i_1,t}^{\tau_{i_1}/2}}{\partial \nu_j} + \dfrac{1}{h_{i_1,t}^{\tau_{i_1}/2}}\dfrac{\partial^2 h_{i_1,t}^{\tau_{i_1}/2}}{\partial \tau_{i_0}\partial \nu_j}  \right]
\end{align*}
	\item for $i_1 = 1,\ldots, m$, we have
\small{	
\begin{align*}
\dfrac{\partial^2 D_{0t}(i_1,j_1)}{\partial\tau_{i_1}^2} & = \dfrac{1}{\tau_{i_1}}h_{i_1,t}^{1/2}\left[-\dfrac{1}{\tau_{i_1}}\left(-\dfrac{1}{\tau_{i_1}} \log\left(h_{i_1,t}^{\tau_{i_1}/2}\right) + \dfrac{1}{h_{i_1,t}^{\tau_{i_1}/2}} \dfrac{\partial h_ {i_1,t}^{\tau_{i_1}/2}}{\partial \tau_{i_1}} \right) \right.\\
&\qquad + \dfrac{1}{\tau_{i_1}}\left(-\dfrac{1}{\tau_{i_1}} \log\left(h_{i_1,t}^{\tau_{i_1}/2}\right) + \dfrac{1}{h_{i_1,t}^{\tau_{i_1}/2}} \dfrac{\partial h_ {i_1,t}^{\tau_{i_1}/2}}{\partial \tau_{i_1}} \right)^2 + \dfrac{1}{\tau_{i_1}^2}\log\left(h_{i_1,t}^{\tau_{i_1}/2}\right) \\
&\qquad  - \dfrac{1}{\tau_{i_1}}\dfrac{1}{h_{i_1,t}^{\tau_{i_1}/2}} \dfrac{\partial h_{i_1,t}^{\tau_{i_1}/2}}{\partial \tau_{i_1}} - \dfrac{1}{h_{i_1,t}^{\tau_{i_1}/2}} \dfrac{\partial h_{i_1,t}^{\tau_{i_1}/2}}{\partial \tau_{i_1}}\dfrac{1}{h_{i_1,t}^{\tau_{i_1}/2}} \dfrac{\partial h_{i_1,t}^{\tau_{i_1}/2}}{\partial \tau_{i_1}} \left.+ \dfrac{1}{h_{i_1,t}^{\tau_{i_1}/2}} \dfrac{\partial^2 h_{i_1,t}^{\tau_{i_1}/2}}{\partial \tau_{i_1}^2}\right]
\end{align*}
}
\normalsize
	\item finally for $i_1,i_0 = 1,\ldots, m$, we obtain
\small{
\begin{align*}
\dfrac{\partial^2 D_{0t}(i_1,j_1)}{\partial\tau_{i_1}\partial\tau_{i_0}} &= \dfrac{1}{\tau_{i_1}} h_{i_1,t}^{1/2}\left[\dfrac{1}{\tau_{i_0}}\left(-\dfrac{1}{\tau_{i_0}}\log\left(h_{i_1,t}^{\tau_{i_1}/2}\right) + \dfrac{1}{h_{i_1,t}^{\tau_{i_1}/2}}\dfrac{\partial h_{i_1,t}^{\tau_{i_1}/2}}{\partial \tau_{i_0}}\right)\left(-\dfrac{1}{\tau_{i_1}} \log\left(h_{i_1,t}^{\tau_{i_1}/2}\right) + \dfrac{1}{h_{i_1,t}^{\tau_{i_1}/2}} \dfrac{\partial h_ {i_1,t}^{\tau_{i_1}/2}}{\partial \tau_{i_1}} \right) \right.\\
&\qquad\left. - \dfrac{1}{\tau_{i_1}}\dfrac{1}{h_{i_1,t}^{\tau_{i_1}/2}} \dfrac{\partial h_{i_1,t}^{\tau_{i_1}/2}}{\partial \tau_{i_0}} + \dfrac{1}{h_{i_1,t}^{\tau_{i_1}/2}} \dfrac{\partial h_{i_1,t}^{\tau_{i_1}/2}}{\partial \tau_{i_0}}\dfrac{1}{h_{i_1,t}^{\tau_{i_1}/2}} \dfrac{\partial h_{i_1,t}^{\tau_{i_1}/2}}{\partial \tau_{i_1}} + \dfrac{1}{h_{i_1,t}^{\tau_{i_1}/2}} \dfrac{\partial^2 h_{i_1,t}^{\tau_{i_1}/2}}{\partial \tau_{i_1}\partial \tau_{i_0}}\right]
\end{align*}
}
\normalsize
\end{enumerate}
Since the first order derivatives are already controlled, and since the estimations done in the case with known power,  it remains to prove that
\begin{align}\label{DSTau}
\mathbb{E}\left\vert \dfrac{1}{h_{i_1,t}^{\tau_{i_1}/2}} \dfrac{\partial^2h_{i_1,t}^{\tau_{i_1}/2}}{\partial \tau_i \partial \nu_j}(\nu_0)\right\vert^{r_0} < \infty ,
\text{ for }i=1,...,m\text{ and }j=1,\ldots,s_2.
\end{align}
By \eqref{Hgrandinfinibis}, we have
\begin{equation*}
\mathbb{H}_t^{(\tau)}(i_1) = \sum\limits_{k=0}^\infty\sum\limits_{j_1=1}^d \mathbb{B}^k(i_1,j_1)\underline{c}_{t-k}^{(\tau)}(j_1).
\end{equation*}
It is easy to notice that
\[ \dfrac{\partial^2 \mathbb{H}_t^{(\tau)}(i_1)}{\partial \tau_i \partial \omega_j} = 0,\]
and for $i\neq j$
\[ \dfrac{\partial^2 \mathbb{H}_t^{(\tau)}(i_1)}{\partial \tau_{i} \partial \tau_{j}} = 0.\]
It remains to treat three cases to prove \eqref{DSTau}.
\begin{enumerate}[\qquad(a)]
	\item For $i=1,\ldots, m$  we have
	\[\dfrac{\partial^2\mathbb{H}_t^{(\tau)}(i_1)}{\partial\tau_i\partial\alpha^\pm_j} = \sum\limits_{k=0}^\infty\sum\limits_{j_1=1}^m\sum\limits_{l=1}^q \mathbb{B}^k(i_1,j_1)\dfrac{\partial A_l^\pm(i_1,j_1)}{\partial \alpha_j^\pm} \log(\pm{\varepsilon}_{i_1,t-l}^\pm)(\pm{\varepsilon}_{i_1,t-l}^\pm)^{\tau_{i_1}}, \]
where $\partial A_l^\pm(i_1,j_1)/\partial \alpha_j^\pm$ is a matrix with only one non null element equal 1 at the place of $\alpha_j^\pm$. It follows that
\begin{align*}
	\alpha_j^\pm\dfrac{\partial^2\mathbb{H}_t^{(\tau)}(i_1)}{\partial\tau_i\partial\alpha^\pm_j} &\leq \sum\limits_{k=0}^\infty\sum\limits_{j_1=1}^m\sum\limits_{l=1}^q \mathbb{B}^k (i_1,j_1) \alpha_j^\pm\log(\pm{\varepsilon}_{i_1,t-l}^\pm)(\pm{\varepsilon}_{i_1,t-l}^\pm)^{\tau_{i_1}}\\
	&\leq \sum\limits_{k=0}^\infty\sum\limits_{j_1=1}^m \mathbb{B}^k(i_1,j_1)\dfrac{\partial \underline{c}_{t-k}^{(\tau)}(j_1)}{\partial \tau_i}.
\end{align*}
Using the same techniques used to prove \eqref{esti-dtau_bis} with $\kappa=1$, we obtain that
\[\mathbb{E}\left\vert \dfrac{\alpha_j^\pm}{\mathbb{H}_t^{(\tau)}(i_1)} \dfrac{\partial^2\mathbb{H}_t^{(\tau)}(i_1)}{\partial \tau_i\partial \alpha_j^\pm}\right\vert^{r_0} < \infty.\]
	\item For $i=1, \ldots, m$ it holds
	\[\dfrac{\partial^2\mathbb{H}_t^{(\tau)}(i_1)}{\partial\tau_i\partial\beta_j} = \sum\limits_{k=1}^\infty\sum\limits_{j_1=1}^m\left\{\sum\limits_{l=1}^k\mathbb{B}^{l-1}(i_1,j_1)\mathbb{B}^{(j)}(i_1,j_1)\mathbb{B}^{k-l}(i_1,j_1)\right\}\dfrac{\partial\underline{c}_{t-k}^{(\tau)}(i_1)}{\partial\tau_i}. \]
Consequently
\begin{align*}
\beta_j \dfrac{\partial^2\mathbb{H}_t^{(\tau)}(i_1)}{\partial\tau_i\partial\beta_j} &\leq \sum\limits_{k=1}^\infty\sum\limits_{j_1=1}^m k \mathbb{B}^k(i_1,j_1)\dfrac{\partial\underline{c}_{t-k}^{(\tau)}(i_1)}{\partial\tau_i} ,
\end{align*}
and we proceed as in the previous case to conclude that
\[\mathbb{E}\left\vert \dfrac{\beta_j}{\mathbb{H}_t^{(\tau)}(i_1)}\dfrac{\partial^2 \mathbb{H}_t^{(\tau)}(i_1)}{\partial \tau_i \partial \beta_j}\right\vert^{r_0}< \infty.\]
	\item For $i=1,\ldots, m$
	\[\dfrac{\partial^2\mathbb{H}_t^{(\tau)}(i_1)}{\partial\tau_{i}^2} = \sum\limits_{k=0}^\infty\sum\limits_{j_1=1}^m \mathbb{B}^k(i_1,j_1) \dfrac{\partial^2 \underline{c}_{t-k}^{(\tau)}(i)}{\partial\tau_i^2} , \]
and we use \eqref{esti-dtau_bis} with $\kappa=2$ in order to obtain
\[\mathbb{E}\left\vert \dfrac{1}{\mathbb{H}_t^{(\tau)}(i_1)}\dfrac{\partial^2\mathbb{H}_t^{(\tau)}(i_1)}{\partial\tau_i^2}\right\vert^{r_0}<\infty.\]
\end{enumerate}

The above arguments can be generalized on a neighborhood of  $V(\nu_0)$ of $\nu_0$. So we have for all $i_1=1,\ldots, m$ and all $i,j=1,\ldots, s_1$:
\begin{equation*}
\mathbb{E}\sup\limits_{\nu \in V(\nu_0)}\left\vert \dfrac{1}{h_{i_1,t}^{\tau_{i_1}/2}} \dfrac{\partial^2 h_{i_1,t}^{\tau_{i_1}/2}}{\partial \nu_i\partial \nu_j}(\nu)\right\vert^{r_0} < \infty.
\end{equation*}
Using the same arguments as in Subsection \ref{A34} combined with the previous modification taking into account that the power is unknown (especially the estimation \eqref{esti-dtau_bis} with $\kappa =3$) we obtain
that for $i_1=1,\ldots, m$ and $i,j,k = 1,\ldots, s_1$, it holds
\begin{align*}
\mathbb{E}\sup\limits_{\nu \in V(\nu_0)}\left\vert \dfrac{1}{h_{i_1,t}^{\tau_{i_1}/2}}\dfrac{\partial^3h_{i_1,t}^{\tau_{i_1}/2}}{\partial \nu_i\partial \nu_j\partial\nu_k}(\nu)\right\vert^{r_0} < \infty.
\end{align*}

\subsubsection{Invertibility of the matrix $J$}\label{new-proof2}
Replacing $\theta$ by $\nu$ in the arguments that lead to \eqref{C.44}, one obtains that the rows $1,m+1,\ldots,m^2$ of the following equations
\begin{align}\label{C.44Inconnu}
\textbf{d}\textbf{c}
& = \sum\limits_{i=1}^{s_0}c_i\dfrac{\partial}{\partial \nu_i}\left[(D_{0t}\otimes D_{0t})\Vec(R_0)\right] = 0_{m^2}
\end{align}
yield
\begin{align}\label{C.45Inconnu}
\sum\limits_{i=1}^{s_1}c_i\dfrac{\partial}{\partial \nu_i}\underline{h}_t^{\underline\tau/2}(\nu_0)& = 0_m,\quad a.s.\\
\label{C.45bisInconnu}
(D_{0t} \otimes D_{0t})\sum\limits_{k = s_1+1}^{s_0}c_k\dfrac{\partial}{\partial \nu_k}\text{vec}R_0 &= 0_{m^2},\quad a.s.
\end{align}
Under \textbf{A8}, Equation \eqref{C.45bisInconnu} is equivalent to
\begin{equation*}
\sum_{i=s_1+1}^{s_0} c_k\dfrac{\partial}{\partial \nu_k} \text{vec}R_0 = 0_{m^2}.
\end{equation*}
Since the vectors $\partial \text{vec}R_0/\partial \nu_k$, $k = s_1+1,\ldots, s_0$ are linearly independent, the vector  $\textbf{c}_3 = ({c}_{s_1+1}, \dots,{c}_{s_0})'$ is null.
Consequently, Equation \eqref{C.44Inconnu} becomes
\begin{align*}
\sum\limits_{i=1}^{s_2}c_{i}\dfrac{\partial}{\partial \theta_i}\underline{h}_t^{\underline\tau/2}(\nu_0)
+\sum\limits_{i=s_2+1}^{s_2+m}c_{i}\dfrac{\partial}{\partial \nu_i}\underline{h}_t^{\underline\tau/2}(\nu_0) = 0_m,\quad a.s.
\end{align*}
or equivalently
\begin{align}\label{C.45Inconnubis}
\sum\limits_{i=1}^{s_2}c_{1,i}\dfrac{\partial}{\partial \theta_i}\underline{h}_t^{\underline\tau/2}(\nu_0)
+\sum\limits_{i=1}^{m}c_{2,i}\dfrac{\partial}{\partial \nu_i}\underline{h}_t^{\underline\tau/2}(\nu_0) = 0_m,\quad a.s.
\end{align}
In view of \eqref{MAPGARCH}, the $i_1^{\mbox{th}}$ component of $\underline{h}_{t}^{\underline{\tau}/2}(\nu_0)$ is
\begin{equation}\label{h.i_1}
\begin{aligned}
h_ {i_1,t}^{\tau_{i_1}/2} &= \underline{\omega}_0(i_1) +  \sum\limits_{i_2=1}^m \sum\limits_{i=1}^q  \left( A_{0i}^+(i_1,i_2)\left({\varepsilon}_{i_2,t-i}^+\right)^{\tau_{i_2}} + A_{0i}^-(i_1,i_2)\left(-{\varepsilon}_{i_2,t-i}^-\right)^{\tau_{i_2}}\right)\\
&\qquad +  \sum\limits_{i_2=1}^m \sum\limits_{i=1}^p B_{0i}(i_1,i_2)h_ {i_2,t-i}^{\tau_{i_2}/2}
\end{aligned}
\end{equation}
For $i_1=1,\dots,m$, from \eqref{DDbis} and \eqref{Dtau}, Equation \eqref{C.45Inconnubis}
reduces as
\begin{align}\label{eqC.45Inconnubis}
 h_{i_1,t}^{1/2}\dfrac{1}{\tau_{i_1}^2} \dfrac{1}{h_{i_1,t}^{\tau_{i_1}/2}} \left[\sum\limits_{i=1}^{s_2}c_{1,i}\tau_{i_1}\dfrac{\partial h_ {i_1,t}^{\tau_{i_1}/2}}{\partial \theta_{i}}(\nu_0)  +\sum\limits_{i=1}^{m}c_{2,i}\tau_{i_1}\dfrac{\partial h_ {i_1,t}^{\tau_{i_1}/2}}{\partial \tau_{i}}(\nu_0)-c_{2,i_1}h_{i_1,t}^{\tau_{i_1}/2}\log\left(h_{i_1,t}^{\tau_{i_1}/2}\right)\right]=0,\quad a.s.
\end{align}
From \eqref{h.i_1}, the derivatives are defined recursively by
\begin{equation*}
\begin{aligned}
\dfrac{\partial h_ {i_1,t}^{\tau_{i_1}/2}}{\partial \theta}(\nu) &= \underline{\tilde{c}}_t(\nu) +
 \sum\limits_{i_2=1}^m \sum\limits_{i=1}^p B_{i}(i_1,i_2)
\dfrac{\partial h_ {i_2,t-i}^{\tau_{i_2}/2}}{\partial \theta}(\nu),\\
\dfrac{\partial h_ {i_1,t}^{\tau_{i_1}/2}}{\partial \tau_j}(\nu) &=
\sum\limits_{i=1}^q  \left( A_{i}^+(i_1,j)\log\left({\varepsilon}_{j,t-i}^+\right)\left({\varepsilon}_{j,t-i}^+\right)^{\tau_{j}} + A_{i}^-(i_1,i_2)\log\left(-{\varepsilon}_{j,t-i}^-\right)\left(-{\varepsilon}_{j,t-i}^-\right)^{\tau_{j}}\right)\\
&\qquad +  \sum\limits_{i_2=1}^m \sum\limits_{i=1}^p B_{i}(i_1,i_2)
 \dfrac{\partial h_ {i_2,t-i}^{\tau_{i_2}/2}}{\partial \tau_j}(\nu)
\end{aligned}
\end{equation*}
for $j=1,\dots,m$, where $\underline{\tilde{c}}_t(\nu)$ is defined by
\begin{equation}\label{ctilde}
\begin{aligned}
\underline{\tilde{c}}_t(\nu) &= \left(0,\dots,1,0,\dots, \left({\varepsilon}_{i_1,t-1}^+\right)^{\tau_{i_1}},0, \dots,\left({\varepsilon}_{i_m,t-1}^+\right)^{\tau_{i_m}},0,\dots, \left({\varepsilon}_{i_1,t-q}^+\right)^{\tau_{i_1}},0, \dots,\left({\varepsilon}_{i_m,t-q}^+\right)^{\tau_{i_m}},\right.\\ &\qquad\left. 0,
,\dots, \left(-{\varepsilon}_{i_1,t-1}^-\right)^{\tau_{i_1}},0, \dots,\left(-{\varepsilon}_{i_m,t-1}^-\right)^{\tau_{i_m}},0,\dots, \left(-{\varepsilon}_{i_1,t-q}^-\right)^{\tau_{i_1}},0, \dots,\left(-{\varepsilon}_{i_m,t-q}^-\right)^{\tau_{i_m}},\right.\\ &\qquad\left. 0,
,\dots, h_ {i_1,t-1}^{\tau_{i_1}/2},0, \dots,h_ {i_m,t-1}^{\tau_{i_m}/2},0,\dots, h_ {i_1,t-p}^{\tau_{i_1}/2},0, \dots,h_ {i_m,t-p}^{\tau_{i_m}/2},\dots, 0\right)'.
\end{aligned}
\end{equation}
Let $R_t$ a random variable measurable with respect to $\sigma\{\eta_u, u \leq t\}$.

We restrict ourselves to the particular case of a CCC-APGARCH$(p,1)$ (see \eqref{MAPGARCH}). The general case can be easily deduced from the following arguments. We thus have
\begin{align*}
\dfrac{\partial h_ {i_1,t}^{\tau_{i_1}/2}}{\partial \tau_j}(\nu)  & = A_{1}^+(i_1,j)\log\left({\varepsilon}_{j,t-1}^+\right)\left({\varepsilon}_{j,t-1}^+\right)^{\tau_{j}} + A_{1}^-(i_1,j)\log\left(-{\varepsilon}_{j,t-1}^-\right)\left({-\varepsilon}_{j,t-1}^-\right)^{\tau_{j}} + R_{t-2}, \\
h_ {i_1,t}^{\tau_{i_1}/2}(\nu) &= \sum\limits_{i_2=1}^m \left( A_{1}^+(i_1,i_2)\left({\varepsilon}_{i_2,t-1}^+\right)^{\tau_{i_2}} + A_{1}^-(i_1,i_2)\left(-{\varepsilon}_{i_2,t-1}^-\right)^{\tau_{i_2}}\right) + R_{t-2}.
\end{align*}
Combining the above expressions and \eqref{ctilde}, under \textbf{A8}, Equation  \eqref{eqC.45Inconnubis} becomes
\begin{align*}
0&=\sum\limits_{i=1}^{m}\tau_{i_1}\left( c_{1,i_1+i\cdot m}\left({\varepsilon}_{i,t-1}^+\right)^{\tau_{i}}+
 c_{1,i_1+i\cdot m^2}\left(-{\varepsilon}_{i,t-1}^-\right)^{\tau_{i}}\right) + R_{t-2}^1  \\&\qquad+\sum\limits_{i=1}^{m}c_{2,i}\tau_{i_1}\left(A_{01}^+(i_1,i)\log\left({\varepsilon}_{i,t-1}^+\right)\left({\varepsilon}_{i,t-1}^+\right)^{\tau_{i}} + A_{01}^-(i_1,i)\log\left(-{\varepsilon}_{i,t-1}^-\right)\left(-{\varepsilon}_{i,t-1}^-\right)^{\tau_{i}} + R_{t-2}^2\right)\\&\qquad-c_{2,i_1}\left(\sum\limits_{i_2=1}^m \left( A_{01}^+(i_1,i_2)\left({\varepsilon}_{i_2,t-1}^+\right)^{\tau_{i_2}} + A_{01}^-(i_1,i_2)\left(-{\varepsilon}_{i_2,t-1}^-\right)^{\tau_{i_2}}\right) + R_{t-2}^3\right)\\&\qquad\times\log\left[\sum\limits_{i_2=1}^m \left( A_{01}^+(i_1,i_2)\left({\varepsilon}_{i_2,t-1}^+\right)^{\tau_{i_2}} + A_{01}^-(i_1,i_2)\left(-{\varepsilon}_{i_2,t-1}^-\right)^{\tau_{i_2}}\right) + R_{t-2}^3\right],\quad a.s.
\end{align*}
which is equivalent, almost surely, to the following two equations
\begin{equation}\label{eq+}
\begin{aligned}
\sum\limits_{i=1}^{m}\tau_{i_1} c_{1,i_1+i\cdot m}\left({\varepsilon}_{i,t-1}^+\right)^{\tau_{i}}
 + R_{t-2}^1 +\sum\limits_{i=1}^{m}c_{2,i}\tau_{i_1}&\left(A_{01}^+(i_1,i)\log\left({\varepsilon}_{i,t-1}^+\right)\left({\varepsilon}_{i,t-1}^+\right)^{\tau_{i}}  + R_{t-2}^2\right)\\ \qquad-c_{2,i_1}\left(\sum\limits_{i_2=1}^m  A_{01}^+(i_1,i_2)\left({\varepsilon}_{i_2,t-1}^+\right)^{\tau_{i_2}}  + R_{t-2}^3\right)&\log\left[\sum\limits_{i_2=1}^m  A_{01}^+(i_1,i_2)\left({\varepsilon}_{i_2,t-1}^+\right)^{\tau_{i_2}} + R_{t-2}^3\right]=0
\end{aligned}
\end{equation}
\begin{equation}\label{eq-}
\begin{aligned}
\sum\limits_{i=1}^{m}\tau_{i_1} c_{1,i_1+i\cdot m^2}\left(-{\varepsilon}_{i,t-1}^-\right)^{\tau_{i}}
 + R_{t-2}^1 +\sum\limits_{i=1}^{m}c_{2,i}\tau_{i_1}&\left(A_{01}^-(i_1,i)\log\left(-{\varepsilon}_{i,t-1}^-\right)\left(-{\varepsilon}_{i,t-1}^-\right)^{\tau_{i}}  + R_{t-2}^2\right)\\ \qquad-c_{2,i_1}\left(\sum\limits_{i_2=1}^m  A_{01}^-(i_1,i_2)\left(-{\varepsilon}_{i_2,t-1}^-\right)^{\tau_{i_2}}  + R_{t-2}^3\right)&\log\left[\sum\limits_{i_2=1}^m  A_{01}^-(i_1,i_2)\left(-{\varepsilon}_{i_2,t-1}^-\right)^{\tau_{i_2}} + R_{t-2}^3\right]=0.
\end{aligned}
\end{equation}
It follows from \textbf{A8} that, for some $\prod\limits_{i=1}^{m}[a_i,b_i]\subset[0,+\infty[^m$,
\begin{equation}\label{eq+bis}
\begin{aligned}
\sum\limits_{i=1}^{m}\tau_{i_1} c_{1,i_1+i\cdot m}x_i
 + R_{t-2}^1 +\sum\limits_{i=1}^{m} c_{2,i}&\left(A_{01}^+(i_1,i)x_i\log\left(x_i\right) + R_{t-2}^2\right)\\ \qquad-c_{2,i_1}\left(\sum\limits_{i_2=1}^m  A_{01}^+(i_1,i_2)x_{i_2} + R_{t-2}^3\right)&\log\left[\sum\limits_{i_2=1}^m  A_{01}^+(i_1,i_2)x_{i_2} + R_{t-2}^3\right]=0,\qquad a.s.
\end{aligned}
\end{equation}
\begin{equation}\label{eq-bis}
\begin{aligned}
\sum\limits_{i=1}^{m}\tau_{i_1} c_{1,i_1+i\cdot m^2}x_i
 + R_{t-2}^1 +\sum\limits_{i=1}^{m} c_{2,i}&\left(A_{01}^-(i_1,i)x_i\log\left(x_i\right) + R_{t-2}^2\right)\\ \qquad-c_{2,i_1}\left(\sum\limits_{i_2=1}^m  A_{01}^-(i_1,i_2)x_{i_2} + R_{t-2}^3\right)&\log\left[\sum\limits_{i_2=1}^m  A_{01}^-(i_1,i_2)x_{i_2} + R_{t-2}^3\right]=0,\qquad a.s.
\end{aligned}
\end{equation}
for any $(x_1,\dots,x_m)\in\prod\limits_{i=1}^{m}[a_i,b_i]$.
Differentiating three times Equations  \eqref{eq+bis} and \eqref{eq-bis} with respect to $x_i$, for $i=1,\dots,m$, we obtain that
$$\left(c_{2,i}-c_{2,i_1}\right)\left[\left(A_{01}^+(i_1,i)\right)^2+\left(A_{01}^-(i_1,i)\right)^2\right]=0.$$
Under Assumption $\textbf{A4}:$ if $p>0, \mathcal{A}_0^+(1) + \mathcal{A}_0^-(1) \neq 0,$
it is impossible to have $A_{01}^+(i_1,i)=A_{01}^-(i_1,i)=0$, for all $i=1,\dots,m$. Then there exists $i_0$ such that $\left(A_{01}^+(i_1,i_0)\right)^2+\left(A_{01}^-(i_1,i_0)\right)^2\neq 0$ and then we have $c_{2,i_0}-c_{2,i_1}=0$. Equations  \eqref{eq+bis} and \eqref{eq-bis} become
\begin{equation}\label{eq+bis1}
\begin{aligned}
\sum\limits_{i=1}^{m}\tau_{i_1} c_{1,i_1+i\cdot m}x_i
 + R_{t-2}^1 +\sum\limits_{i=1}^{m} c_{2,i}&\left(A_{01}^+(i_1,i)x_i\log\left(x_i\right) + R_{t-2}^2\right)\\ \qquad-\dfrac{1}{2}c_{2,i_0}\left(\sum\limits_{i_2=1}^m  A_{01}^+(i_1,i_2)x_{i_2} + R_{t-2}^3\right)&\log\left[\sum\limits_{i_2=1}^m  A_{01}^+(i_1,i_2)x_{i_2} + R_{t-2}^3\right]=0,\qquad a.s.
\end{aligned}
\end{equation}
\begin{equation}\label{eq-bis1}
\begin{aligned}
\sum\limits_{i=1}^{m}\tau_{i_1} c_{1,i_1+i\cdot m^2}x_i
 + R_{t-2}^1 +\sum\limits_{i=1}^{m} c_{2,i}&\left(A_{01}^-(i_1,i)x_i\log\left(x_i\right) + R_{t-2}^2\right)\\ \qquad-\dfrac{1}{2}c_{2,i_0}\left(\sum\limits_{i_2=1}^m  A_{01}^-(i_1,i_2)x_{i_2} + R_{t-2}^3\right)&\log\left[\sum\limits_{i_2=1}^m  A_{01}^-(i_1,i_2)x_{i_2} + R_{t-2}^3\right]=0,\qquad a.s.
\end{aligned}
\end{equation}
Differentiating twice again Equations  \eqref{eq+bis1} and \eqref{eq-bis1} with respect to $x_{i_0}$, we find that
$$c_{2,i_0}\left(A_{01}^+(i_1,i_0)+A_{01}^-(i_1,i_0)\right)\left[\sum\limits_{\substack{i_2=1\\i_2\neq i_0}}^m \left( A_{01}^+(i_1,i_2)\left({\varepsilon}_{i_2,t-1}^+\right)^{\tau_{i_2}} + A_{01}^-(i_1,i_2)\left(-{\varepsilon}_{i_2,t-1}^-\right)^{\tau_{i_2}}\right) + R_{t-2}^3\right]=0.$$
Since the law of $\eta_t$ is non degenerated (see Assumption \textbf{A3}), we deduce that $$\sum\limits_{\substack{i_2=1\\i_2\neq i_0}}^m \left( A_{01}^+(i_1,i_2)\left({\varepsilon}_{i_2,t-1}^+\right)^{\tau_{i_2}} + A_{01}^-(i_1,i_2)\left(-{\varepsilon}_{i_2,t-1}^-\right)^{\tau_{i_2}}\right) + R_{t-2}^3\neq 0 \ .$$
So $c_{2,i_0}=c_{2,i_1}=0$. Since $i_1$ is arbitrary, we deduce that $c_{2,i}=0$ for any $i=1,...,m$, or equivalently that $c_{1,i_1+i_0\cdot m}=c_{1,i_1+i_0\cdot m^2}=0$.

Thus the vector $\textbf{c}_2 = ({c}_{2,1}, \dots,{c}_{2,m})'= ({c}_{s_2+1}, \dots,{c}_{s_1})'$ is null.
We recall that  $\textbf{c}_3 = ({c}_{s_1+1}, \dots,{c}_{s_0})'$ is null,  the invertibility of the matrix $J$ is thus shown in this case of a CCC-APGARCH$(p,1)$.

In the general case of a CCC-APGARCH$(p,q)$, we show by induction that \eqref{C.45Inconnubis} entails necessarily
$$A_{01}^+(i_1,i_0)+A_{01}^-(i_1,i_0)=\dots=A_{0q}^+(i_1,i_0)+A_{0q}^-(i_1,i_0)=0,\qquad\forall,\, i_0,i_1=1,\dots,m,$$
which is impossible under Assumption $\textbf{A4}$ and thus $\textbf{c} = 0$.

This is in contradiction with $\textbf{c}'\textbf{h}'\textbf{h}\textbf{c} = \textbf{c}'\textbf{d}'\textbf{H}^2\textbf{d}\textbf{c} = 0$ almost-surely. Therefore the assumption that  $J$ is not singular is absurd.

\subsubsection{Asymptotic irrelevance  of the initial values}
It suffices to adapt the arguments  used in Subsection \ref{null-effect} when the power is known.

\bibliographystyle{apalike}
\bibliography{biblio-oth}
\addcontentsline{toc}{section}{References}
\clearpage
\pagenumbering{roman}

\begin{center}
{\Large
Estimation of multivariate asymmetric power GARCH models: \\[4mm] 
{\bf
Complementary results that are not submitted for publication}}
\end{center}
\vspace{2cm}

\section{Details on the proof of Theorem \ref{le-joli-label-de-yacouba}}

\subsection{Minimization of the quasi-likelihood on the true value}\label{A23a}
The criterion is not integrable on all point, but we first prove that $\mathbb{E}_{\theta_0}\left[l_t(\theta)\right]$ is well defined on $\mathbb{R} \cup \{+\infty\}$ for all $\theta$. Indeed we have
\begin{equation*}
\mathbb{E}_{\theta_0}\left[l_t^-(\theta)\right] \leq \mathbb{E}_{\theta_0}\left[\log^-\vert H_t \vert\right] = \mathbb{E}_{\theta_0}\left[ \log^-\vert D_t R D_t \vert\right] \leq \max\{0, \log(\vert R \vert \min\limits_{i} \underline{\omega}(i)^m)\} < \infty.
\end{equation*}
Now we show that  $\mathbb{E}_{\theta_0}\left[l_t(\theta_0)\right]$ is well defined on $\mathbb{R}$. We use $|\det(A)|\le \rho(A)^m\le  \|A\|^m$, the Jensen inequality and the  Corollary \ref{cor2}
\begin{equation*}
\begin{aligned}
\mathbb{E}_{\theta_0}\left[\log\vert H_t(\theta_0)\vert \right] &= \mathbb{E}_{\theta_0}\left[\dfrac{m}{s} \log\vert H_t(\theta_0)\vert^{\frac{s}{m}}\right] \\
& \leq \dfrac{m}{s}\log\left(\mathbb{E}_{\theta_0}\Vert H_t(\theta_0) \Vert ^{\frac{s}{m}}\right)\\
& \leq \dfrac{m}{s}\log\left(\mathbb{E}_{\theta_0}\Vert R \Vert^s \Vert D_{0t} \Vert ^{2s}\right) \\
& \leq K + \dfrac{m}{s} \log\left(\mathbb{E}_{\theta_0} \Vert D_{0t} \Vert^{2s}\right) \\
& = K + \dfrac{m}{s} \log\left(\mathbb{E}_{\theta_0}\left[\max\limits_{i}\left(h_{it}(\theta_0)\right)^s\right]\right)\\
& \leq K + \dfrac{m}{s} \log\left(\mathbb{E}_{\theta_0}\Vert \underline{h}_{t}(\theta_0)\Vert^s\right) <\infty\ .
\end{aligned}
\end{equation*}
Therefore we have
\begin{equation*}
\begin{aligned}
\mathbb{E}_{\theta_0}\left[l_t(\theta_0)\right] &= \mathbb{E}_{\theta_0}\left[ \varepsilon_t' H_t^{-1}(\theta_0) \varepsilon_t + \log\vert H_t(\theta_0) \vert \right]\\
& = \mathbb{E}_{\theta_0}\left[\eta_t '(H_t^{1/2})'(\theta_0)H_t^{-1}(\theta_0)H_t^{1/2}(\theta_0)\eta_t + \log\vert H_t(\theta_0)\vert \right] \\
& = \mathbb{E}_{\theta_0}\left[\mathrm{Tr}(\eta_t '(H_t^{1/2})'(\theta_0)H_t^{-1}(\theta_0)H_t^{1/2}(\theta_0)\eta_t) + \log\vert H_t(\theta_0)\vert \right] \\
& = \mathbb{E}_{\theta_0}\left[\mathrm{Tr}(\eta_t '\eta_t) + \log\vert H_t(\theta_0)\vert \right] \\
&= m + \mathbb{E}_{\theta_0}\left[\log\vert H_t(\theta_0)\vert \right] < +\infty.
\end{aligned}
\end{equation*}
Since $\mathbb{E}_{\theta_0}\left[l_t^-(\theta)\right] < \infty$ for any $\theta$, $\mathbb{E}_{\theta_0}\left[l_t^-(\theta_0)\right] < \infty$ and we deduce that $\mathbb{E}_{\theta_0}\left[l_t(\theta_0)\right] $ is well defined in $\mathbb{R}$. So, when one studies the function $\theta\mapsto\mathbb{E}_{\theta}\left[l_t(\theta)\right] $, we can restrict our study to the values of $\theta$ such that $\mathbb{E}_{\theta_0} \vert l_t(\theta)\vert < \infty$.

We denote $\lambda_{i,t}$ the positive eigenvalues of $H_t(\theta_0)H_t^{-1}(\theta)$. We have
\begin{equation*}
\begin{aligned}
\mathbb{E}_{\theta_0}\left[l_t(\theta)\right] - \mathbb{E}_{\theta_0}\left[l_t(\theta_0)\right] &= \mathbb{E}_{\theta_0}\left[\varepsilon_t'H_t^{-1}(\theta)\varepsilon_t + \log\vert H_t(\theta) \vert \right] -  \mathbb{E}_{\theta_0}\left[\varepsilon_t'H_t^{-1}(\theta_0)\varepsilon_t + \log\vert H_t(\theta_0) \vert \right] \\
&= \mathbb{E}_{\theta_0}\left[\varepsilon_t'H_t^{-1}(\theta)\varepsilon_t - \varepsilon_t'H_t^{-1}(\theta_0)\varepsilon_t \right] + \mathbb{E}_{\theta_0}\left[\log\vert H_t(\theta) \vert - \log\vert H_t(\theta_0) \vert \right]\\
&= \mathbb{E}_{\theta_0}\left[\varepsilon_t' \{H_t^{-1}(\theta) - H_t^{-1}(\theta_0)\}\varepsilon_t \right] + \mathbb{E}_{\theta_0}\left[\log\vert H_t(\theta)H_t^{-1}(\theta_0) \vert \right]\\
&= \mathbb{E}_{\theta_0}\left[\eta_t'(H_t^{1/2})'(\theta_0) \{H_t^{-1}(\theta) - H_t^{-1}(\theta_0)\}H_t^{1/2}(\theta_0)\eta_t \right]\\
& \qquad \quad+ \mathbb{E}_{\theta_0}\left[\log\vert H_t(\theta)H_t^{-1}(\theta_0) \vert \right]\\
&= \mathbb{E}_{\theta_0}\left[\eta_t'\{ (H_t^{1/2})'(\theta_0)H_t^{-1}(\theta) H_t^{1/2}(\theta_0)  - (H_t^{1/2})'(\theta_0)H_t^{-1}(\theta_0)H_t^{1/2}(\theta_0)\}\eta_t \right]\\
& \qquad \quad+ \mathbb{E}_{\theta_0}\left[\log\vert H_t(\theta)H_t^{-1}(\theta_0) \vert \right]\\
&= \mathbb{E}_{\theta_0}\left[\eta_t'\{ (H_t^{1/2})'(\theta_0)H_t^{-1}(\theta) H_t^{1/2}(\theta_0)  - I_m\}\eta_t \right] + \mathbb{E}_{\theta_0}\left[\log\vert H_t(\theta)H_t^{-1}(\theta_0) \vert \right]\\
&= \mathrm{Tr}\left\{\mathbb{E}_{\theta_0}\left[(H_t^{1/2})'(\theta_0)H_t^{-1}(\theta) H_t^{1/2}(\theta_0)  - I_m\right]\mathbb{E}\left[\eta_t\eta_t'\right]\right\}\\
&\qquad \quad+ \mathbb{E}_{\theta_0}\left[\log\vert H_t(\theta)H_t^{-1}(\theta_0) \vert \right]\\
&= \mathbb{E}_{\theta_0}\left[\mathrm{Tr}\left((H_t^{1/2})'(\theta_0)H_t^{-1}(\theta) H_t^{1/2}(\theta_0)  - I_m\right)\right] + \mathbb{E}_{\theta_0}\left[\log\vert H_t(\theta)H_t^{-1}(\theta_0) \vert \right]\\
&= \mathbb{E}_{\theta_0}\left[\sum\limits_{i=1}^m(\lambda_{i,t} - 1 - \log(\lambda_{i,t}))\right] \geq 0,
\end{aligned}
\end{equation*}
where we have used the inequality $\log(x) \leq x - 1$, for all $x>0$. But, if $x=1$, $\log(x) = x-1$, then the inequality is strict, except if for all $i$, $\lambda_{i,t} = 1$ $\mathbb{P}_{\theta_0}-$ a.s. This condition means that $H_t(\theta) = H_t(\theta_0)$ $\mathbb{P}_{\theta_0}-$a.s. It follows $\underline{h}_t(\theta) = \underline{h}_t(\theta_0)$ $\mathbb{P}_{\theta_0}-$ a.s. and $R(\theta) = R(\theta_0)$. By the identifiability proved in subsection \ref{A22}, we deduce that $\theta = \theta_0$.
\subsection{Proof of  \textup{(\ref{iv})}}\label{A24a}
We recall that
$${\mathcal{L}}_n(\theta) = \dfrac1n \sum\limits_{t=1}^n {l}_t, \qquad {l}_t = {l}_t(\theta) = \underline{\varepsilon}_t' {H}_t^{-1}\underline{\varepsilon}_t + \log \vert {H}_t\vert.$$
For $\theta \in \Theta$ and an integer $k>0$, we denote $V_k(\theta)$ the open ball of radius $1/k$ centered on $\theta$. We have
\begin{equation*}
\begin{aligned}
\inf\limits_{\theta^\ast \in V_k(\theta) \cap \Theta} \tilde{\mathcal{L}}_n(\theta^\ast) &= \inf\limits_{\theta^\ast \in V_k(\theta) \cap \Theta} \left(\mathcal{L}_n(\theta^\ast) + \tilde{\mathcal{L}}_n(\theta^\ast) - \mathcal{L}_n(\theta^\ast)\right) \\
& \geq \inf\limits_{\theta^\ast \in V_k(\theta) \cap \Theta} \mathcal{L}_n(\theta^\ast) + \inf\limits_{\theta^\ast \in V_k(\theta) \cap \Theta}\left(\tilde{\mathcal{L}}_n(\theta^\ast) - \mathcal{L}_n(\theta^\ast)\right)\\
& \geq \inf\limits_{\theta^\ast \in V_k(\theta) \cap \Theta} \mathcal{L}_n(\theta^\ast) - \sup\limits_{\theta^\ast \in V_k(\theta) \cap \Theta} \left\vert\tilde{\mathcal{L}}_n(\theta^\ast) - \mathcal{L}_n(\theta^\ast)\right\vert
\end{aligned}
\end{equation*}
Then
\begin{equation*}
\begin{aligned}
\underset{n\to+\infty}{\liminf}\inf\limits_{\theta^\ast \in V_k(\theta) \cap \Theta} \tilde{\mathcal{L}}_n(\theta^\ast) & \geq \underset{n\to+\infty}{\liminf} \inf\limits_{\theta^\ast \in V_k(\theta) \cap \Theta} \mathcal{L}_n(\theta^\ast) + \underset{n\to+\infty}{\liminf}\left(-\sup\limits_{\theta^\ast \in V_k(\theta) \cap \Theta} \left\vert\tilde{\mathcal{L}}_n(\theta^\ast) - \mathcal{L}_n(\theta^\ast)\right\vert\right) \\
& \geq \underset{n\to+\infty}{\liminf} \inf\limits_{\theta^\ast \in V_k(\theta) \cap \Theta} \mathcal{L}_n(\theta^\ast) - \underset{n\to+\infty}{\limsup}\sup\limits_{\theta\in\Theta} \left\vert\tilde{\mathcal{L}}_n(\theta) - \mathcal{L}_n(\theta)\right\vert ,
\end{aligned}
\end{equation*}
and by subsection \ref{A21} we deduce that
\begin{align}\label{truc}
\underset{n\to+\infty}{\liminf}\inf\limits_{\theta^\ast \in V_k(\theta) \cap \Theta} \tilde{\mathcal{L}}_n(\theta^\ast)  \geq \underset{n\to+\infty}{\liminf} \inf\limits_{\theta^\ast \in V_k(\theta) \cap \Theta} \mathcal{L}_n(\theta^\ast) & = \underset{n\to+\infty}{\liminf}\inf\limits_{\theta^\ast \in V_k(\theta) \cap \Theta} \dfrac1n\sum\limits_{t=1}^n l_t(\theta^\ast) \nonumber \\
& \geq \underset{n\to\infty}{\liminf} \dfrac1n\sum\limits_{t=1}^n \inf\limits_{\theta^\ast \in V_k(\theta) \cap \Theta} l_t(\theta^\ast).
\end{align}
We can not apply the classical ergodic theorem to $\{\inf\limits_{\theta^\ast \in V_k(\theta) \cap \Theta} l_t(\theta^\ast)\}_t$ because of the lack of integrability mentioned in Subsection \ref{A23a}. So we use an extension of the classical ergodic theorem (see \cite{billing} pages 284 and 495) and we have
\[ \underset{n\to+\infty}{\liminf}\  \dfrac1n\sum\limits_{t=1}^n \inf\limits_{\theta^\ast \in V_k(\theta) \cap \Theta} l_t(\theta^\ast) \longrightarrow \mathbb{E}  \left ( \inf\limits_{\theta^\ast \in V_k(\theta) \cap \Theta} l_1(\theta^\ast) \right )\]
By Beppo-Levi's theorem we have, when $k$ goes to infinity
\[\underset{k\to+\infty}{\liminf}\ \mathbb{E} \left ( \inf\limits_{\theta^\ast \in V_k(\theta) \cap \Theta} l_1(\theta^\ast)\right )  = \mathbb{E}\left ( \underset{k\to+\infty}{\liminf} \inf\limits_{\theta^\ast \in V_k(\theta) \cap \Theta} l_1(\theta^\ast) \right ) \longrightarrow \mathbb{E}\left[l_1(\theta)\right]. \]
Hence \eqref{truc} implies that
\begin{align*}
\underset{n\to+\infty}{\liminf}\inf\limits_{\theta^\ast \in V_k(\theta) \cap \Theta} \tilde{\mathcal{L}}_n(\theta^\ast) \ge   \mathbb{E}\left[l_1(\theta)\right]
\end{align*}
and \eqref{iv} follows from the result stated in Subsection \ref{A23a}.
\subsection{Conclusion: proof of Theorem \ref{le-joli-label-de-yacouba}}
By Subsection \ref{A21}, $\lim\limits_{n\to+\infty} \tilde{\mathcal{L}}_n(\theta_0) =  \lim\limits_{n\to+\infty}\mathcal{L}_n(\theta_0)$ and by the ergodic theorem, we have $\lim\limits_{n\to+\infty}\mathcal{L}_n(\theta_0) = \mathbb{E}_{\theta_0}\left[l_1(\theta_0)\right]$. Consequently, $\lim\limits_{n\to+\infty}\mathcal{L}_n(\theta_0)$ exists and is $\mathbb{E}_{\theta_0}\left[l_1(\theta_0)\right]$.
Since $\inf\limits_{\theta^\ast \in V_k(\theta_0) \cap \Theta} \tilde{\mathcal{L}}_n(\theta^\ast) \leq \tilde{\mathcal{L}}_n(\theta_0)$, we have
\begin{equation*}
\underset{n\to+\infty}{\limsup}\inf\limits_{\theta^\ast \in V_k(\theta_0) \cap \Theta}\tilde{\mathcal{L}}_n(\theta^\ast) \leq \underset{n\to+\infty}{\limsup}\, \tilde{\mathcal{L}}_n(\theta_0) = \lim\limits_{n\to+\infty}\tilde{\mathcal{L}}_n(\theta_0) = \lim\limits_{n\to+\infty}\mathcal{L}_n(\theta_0) = \mathbb{E}_{\theta_0}\left[l_1(\theta_0)\right]
\end{equation*}
Then for a small neighborhood ${V}(\theta_0)$ of $\theta_0$
\begin{equation}\label{consistance 1da}
\underset{n\to+\infty}{\limsup}\inf\limits_{\theta^\ast \in V(\theta_0) \cap \Theta}\tilde{\mathcal{L}}_n(\theta^\ast) \leq \underset{n\to+\infty}{\limsup}\, \tilde{\mathcal{L}}_n(\theta_0) = \mathbb{E}_{\theta_0}\left[l_1(\theta_0)\right] .
\end{equation}
The parameter space $\Theta$ can be covered as
\begin{equation*}
\Theta \subseteq {V}(\theta_0) + \bigcup\limits_{\theta \in \Theta }{V}(\theta)
\end{equation*}
with ${V}(\theta)$ is a neighborhood of $\theta$ verifying $(iv)$. 
By the compactness of $\Theta$, there exists a finite covering of $\Theta$
\begin{equation*}
\Theta \subseteq {V}(\theta_0) \cup {V}(\theta_1) \cup \ldots \cup {V}(\theta_k)
\end{equation*}
and then
\begin{equation}\label{consistance 1d2a}
\inf\limits_{\theta\in \Theta} \tilde{\mathcal{L}}_n(\theta) = \min\limits_{i=0,1,\ldots,k} \ \inf\limits_{\theta\in\Theta\cap{V}(\theta_i)} \tilde{\mathcal{L}}_n(\theta).
\end{equation}
Suppose that for all $N$, there exist $n \geq N$ such that $\hat{\theta}_n \in {V}(\theta_{i_0})$ with $i_0 = 1,\ldots, k$. Let $\varepsilon > 0$, we have by \eqref{consistance 1da} that there exists $N_1$ such that for all $n\geq N_1$,
\begin{equation*}
\inf\limits_{\theta^\ast \in \mathcal{V}(\theta_0)} \tilde{\mathcal{L}}_n (\theta^\ast) < \left[\underset{n\to+\infty}{\limsup} \inf\limits_{\theta^\ast \in V(\theta_0)}\tilde{\mathcal{L}}_n(\theta^\ast)\right] + \varepsilon \leq \mathbb{E}_{\theta_0}\left[l_1(\theta_0)\right] + \varepsilon,
\end{equation*}
and by \eqref{iv}, for $i_0 \neq 0$, we obtain the existence of $N_2$ such that for all $n\geq N_2$,
\begin{equation*}
\inf\limits_{\theta^\ast \in \mathcal{V}(\theta_{i_0})} \tilde{\mathcal{L}}_n (\theta^\ast) > \left[\underset{n\to+\infty}{\liminf} \inf\limits_{\theta^\ast \in V(\theta_{i_0})}\tilde{\mathcal{L}}_n(\theta^\ast)\right] - \varepsilon.
\end{equation*}
We can suppose that $N \geq N_1 \cup N_2$ and by the two latter relations and \eqref{consistance 1d2a}, we have for $i_0 \neq 0$
\begin{equation*}
\underset{n\to+\infty}{\liminf} \inf\limits_{\theta^\ast \in V(\theta_{i_0})}\tilde{\mathcal{L}}_n(\theta^\ast) - \varepsilon < \inf\limits_{\theta \in \mathcal{V}(\theta_{i_0})} \tilde{\mathcal{L}}_n(\theta) = \inf\limits_{\theta \in\Theta} \tilde{\mathcal{L}}_n(\theta) \leq \inf\limits_{\theta \in \mathcal{V}(\theta_0)} \tilde{\mathcal{L}}_n(\theta) < \mathbb{E}_{\theta_0}\left[l_1(\theta_0)\right] + \varepsilon.
\end{equation*}
Then we will have
\begin{equation*}
\underset{n\to+\infty}{\liminf} \inf\limits_{\theta^\ast \in V(\theta_{i_0})}\tilde{\mathcal{L}}_n(\theta^\ast) - \varepsilon < \mathbb{E}_{\theta_0}\left[l_1(\theta_0)\right] + \varepsilon,
\end{equation*}
then
\begin{equation}\label{consistance 1d3a}
\underset{n\to+\infty}{\liminf} \inf\limits_{\theta^\ast \in V(\theta_{i_0})}\tilde{\mathcal{L}}_n(\theta^\ast) \leq \mathbb{E}_{\theta_0}\left[l_1(\theta_0)\right].
\end{equation}
But by \eqref{iv} we have $\underset{n\to+\infty}{\liminf} \inf\limits_{\theta^\ast \in V(\theta_{i_0})}\tilde{\mathcal{L}}_n(\theta^\ast) > \mathbb{E}_{\theta_0}\left[l_1(\theta_0)\right]$ and this is in contradiction to \eqref{consistance 1d3a}.\\
For $n$ large enough, we conclude that $\hat{\theta}_n$ belongs to ${V}(\theta_0)$.
\zak

\section{Details on the proof of Theorem \ref{AN-connu}}
\subsection{First derivative of log-likelihood}
\noindent$\bullet$ Proof of \eqref{DL1} and \eqref{DL1bis}.

We differentiate with respect to $\theta_i$ for $i = 1, \ldots, s_1$ (that is with respect to $\underline{\omega}',  {\alpha_{1}^{+}} ', \ldots, {\alpha_q^+}', {\alpha_1^-}', \ldots, {\alpha_q^-}', \beta'_{1}, \ldots, \beta'_{p})'$). 
Indeed,
\begin{align*}
\dfrac{\partial \underline{\varepsilon}_t' D_t^{-1} R^{-1}D_t^{-1} \underline{\varepsilon}_t}{\partial \theta_i} &= Tr\left(\dfrac{\partial\underline{\varepsilon}_t'D_t^{-1}R^{-1}D_t^{-1}\underline{\varepsilon}_t}{D_t'} \dfrac{\partial D_t}{\partial \theta_i}\right) = Tr\left(\dfrac{\partial Tr(\underline{\varepsilon}_t\underline{\varepsilon}_t' D_t^{-1} R^{-1}D_t^{-1})}{D_t'}\dfrac{\partial D_t}{\partial \theta_i}\right)\\
	& = -Tr\left((\underline{\varepsilon}_t\underline{\varepsilon}_t'D_t^{-1}R^{-1} + R^{-1}D_t^{-1}\underline{\varepsilon}_t\underline{\varepsilon}_t') D_t^{-1}\dfrac{D_t}{\partial \theta_i}D_t^{-1}\right), \\
	 2\dfrac{\partial \log(\mathrm{det}(D_t)) + \log(\mathrm{det}(R))}{\partial \theta_i} &= 2Tr\left(\dfrac{\partial \log(\mathrm{det}(D_t))}{\partial D_t'} \dfrac{\partial D_t}{\partial \theta_i} + \underbrace{\dfrac{\partial \log(\det(R))}{\partial D_t'} \dfrac{\partial D_t}{\partial \theta_i}}_{ = 0}\right)\\
	&= 2Tr\left(D_t^{-1} \dfrac{\partial D_t}{\partial \theta_i} \right) ,
\end{align*}
and we obtain \eqref{DL1}.
Indeed, using the property $Tr(AB) = Tr(BA)$ in \eqref{DL1} yields
\begin{equation*}
\begin{aligned}
\dfrac{\partial l_t(\theta_0)}{\partial \theta_i} &= -Tr\left((\underline{\varepsilon}_{t}\underline{\varepsilon}_{t}'D_{0t}^{-1}R_0^{-1} + R_0^{-1}D_{0t}^{-1}\underline{\varepsilon}_{t}\underline{\varepsilon}_{t}') D_{0t}^{-1} D_{0t}^{(i)}D_t^{-1} - 2 D_t^{-1} D_{0t}^{(i)}\right)\\
	&=Tr\left(-D_{0t}\tilde{\eta}_t\tilde{\eta}_t'R_0^{-1}D_{0t}^{-1}D_{0t}^{(i)}D_{0t}^{-1} - R_0^{-1}\tilde{\eta}_t\tilde{\eta}_t' D_{0t}^{(i)}D_{0t}^{-1} + D_{0t}^{-1} D_{0t}^{(i)} + D_{0t}^{(i)}D_{0t}^{-1}\right)\\
	&= Tr\left(-D_{0t}^{-1}D_{0t}\tilde{\eta}_t\tilde{\eta}_t'R_0^{-1}D_{0t}^{-1}D_{0t}^{(i)} +(I_m -  R_0^{-1}\tilde{\eta}_t\tilde{\eta}_t' )D_{0t}^{(i)}D_{0t}^{-1} + D_{0t}^{-1} D_{0t}^{(i)} \right)\\
	&= Tr\left((I_m - \tilde{\eta}_t\tilde{\eta}_t'R_0^{-1})D_{0t}^{-1}D_{0t}^{(i)} +(I_m -  R_0^{-1}\tilde{\eta}_t\tilde{\eta}_t' )D_{0t}^{(i)}D_{0t}^{-1} \right)
\end{aligned}
\end{equation*}
and \eqref{DL1bis} is proved.  

\noindent$\bullet$  Proof of \eqref{DL2} and \eqref{DL2bis}.

We differentiate with respect to $\theta_i$  for $i = s_1+1, \ldots, s_0$ (that is with respect to $\rho'$). 
Indeed we have
\begin{equation*}
\begin{aligned}
 \dfrac{\partial \underline{\varepsilon}_t' D_t^{-1} R^{-1}D_t^{-1} \underline{\varepsilon}_t}{\partial \theta_i} &= Tr\left(\dfrac{\partial\underline{\varepsilon}_t'D_t^{-1}R^{-1}D_t^{-1}\underline{\varepsilon}_t}{R'} \dfrac{\partial R}{\partial \theta_i}\right) = Tr\left(\dfrac{\partial Tr(\underline{\varepsilon}_t\underline{\varepsilon}_t' D_t^{-1} R^{-1}D_t^{-1})}{R'}\dfrac{\partial R}{\partial \theta_i}\right)\\
	& = -Tr\left(R^{-1}D_t^{-1}\underline{\varepsilon}_t\underline{\varepsilon}_t'D_t^{-1}R^{-1} \dfrac{\partial R}{\partial \theta_i}\right),
\end{aligned}
\end{equation*}	
and
\begin{equation*}
\begin{aligned}
2\dfrac{\partial \log(\mathrm{det}(D_t)) + \log(\mathrm{det}(R))}{\partial \theta_i} &= Tr\left(\underbrace{\dfrac{\partial 2\log(\mathrm{det}(D_t))}{\partial R'} \dfrac{\partial R}{\partial \theta_i}}_{=0} + \dfrac{\partial \log(\det(R))}{\partial R'} \dfrac{\partial R}{\partial \theta_i}\right)\\
	&= Tr\left(R^{-1} \dfrac{\partial R}{\partial \theta_i} \right).
\end{aligned}
\end{equation*}	
Hence we obtain \eqref{DL2}. Now, we resume the above computations:
\begin{equation*}
\begin{aligned}
\dfrac{\partial l_t(\theta_0)}{\partial \theta_i} & = - Tr\left(R_0^{-1}D_{0t}^{-1} \underline{\varepsilon}_{t}\underline{\varepsilon}_{t}' D_{0t}^{-1}R_0^{-1} R_0^{(i)} - R_0^{-1}R_0^{(i)}\right)\\
	&= -Tr\left((-I_m+R_0^{-1}D_{0t}^{-1}\underline{\varepsilon}_{t}\underline{\varepsilon}_{t}'D_{0t}^{-1}) R_0^{-1}R_{0}^{(i)}\right)\\
	&= Tr\left((I_m-R_0^{-1}\right.D_{0t}^{-1}\underline{\varepsilon}_{t}\underline{\varepsilon}_{t}'D_{0t}^{-1}\left.) R_0^{-1}R_{0}^{(i)}\right) .
\end{aligned}
\end{equation*}
We remark that $ D_{0t}^{-1}\underline{\varepsilon}_{t}\underline{\varepsilon}_{t}'D_{0t}^{-1} = D_{0t}^{-1}D_{0t}R_0^{1/2}\eta_t \eta_t'R_0^{1/2} D_{0t}D_{0t}^{-1}= \tilde{\eta}_t\tilde{\eta}_t'$ and thus \eqref{DL2bis} is proved.

\newpage
\tableofcontents

\end{document}